\documentclass[10pt,a4paper]{amsart}
\usepackage{amssymb,amsthm,amsmath,mathrsfs}
\usepackage{hyperref}

\usepackage[T1]{fontenc}
\usepackage[utf8]{inputenc}
\usepackage[british]{babel}
\usepackage{geometry}
\usepackage{color}
\usepackage[all]{xy}

\newcommand{\R}{\mathbb{R}}

\newcommand{\cH}{\mathcal{H}}

\newcommand{\G}{\Gamma}
\newcommand{\g}{\gamma}

\newcommand{\fa}{\mathfrak{a}}

\newtheorem{theorem}{Theorem}[section]
\newtheorem{lem}[theorem]{Lemma}
\newtheorem{prop}[theorem]{Proposition}
\newtheorem{corol}[theorem]{Corollary}
\newtheorem{defin}[theorem]{Definition}
\newtheorem{fait}[theorem]{Fact}

\newtheorem{rem}[theorem]{Remark}

\newcommand{\hypG}{Let $G$ be a connected, real linear, semisimple Lie group of non-compact type. }

\title[Topological mixing]{Topological mixing of positive diagonal flows}
\author{Nguyen-Thi Dang}
\begin{document}
\begin{abstract}
Let $G$ be a semi-simple real Lie group without compact factors and $\G<G$ a Zariski dense, discrete subgroup. We study the topological dynamics of positive diagonal flows on $\G \backslash G$.
We extend Hopf coordinates to Bruhat-Hopf coordinates of $G$, which gives the framework to estimate the elliptic part of products of large generic loxodromic elements.
By rewriting results of Guivarc'h-Raugi into Bruhat-Hopf coordinates, we partition the preimage in $\Gamma \backslash G$ of the non-wandering set of mixing regular Weyl chamber flows, into finitely many dynamically conjugated subsets.
We prove a necessary condition for topological mixing, and when the connected component of the identity of the centralizer of the Cartan subgroup is abelian, we prove it is sufficient.
\end{abstract}

\maketitle

\section{Introduction}

Let $G$ be a connected, real semi-simple Lie group of non compact type i.e. without compact factors.
Let $A$ be a maximal split torus i.e. a maximal abelian subgroup whose Lie algebra $\fa$ is a Cartan subspace, 
denote by $\fa^+ \subset \fa$ a choice of closed positive Weyl chamber and by $\fa^{++}$ its interior, by $A^+=\exp(\fa^+)$ and $A^{++}:=\exp{\fa^{++}}$.
Let $\Gamma < G$ be a Zariski dense, discrete subgroup.
We study topological mixing of the right action by translation on $\Gamma \backslash G$ of one parameter subgroups of $A$ that are parametrized by non-trivial elements of $\fa^+$.

\subsection{Previous results}

In the case of lattices\footnote{lattices are Zariski dense subgroups} i.e. $\G \backslash G$ has finite volume for the Haar measure, topological mixing is a consequence of Howe-Moore \cite{HoweMoore79} Theorem.
Moore \cite{Moore87} even proved that it is exponentially mixing for the Haar measure.
\vspace*{.2 cm}

For the isometry group $\mathrm{SO}(n,1)^0$ of $\mathbb{H}^n$, the Cartan subspace $\fa$ is isomorphic to $\mathbb{R}$.
Assume that $\G$ is Zariski dense, discrete and torsion free.
Such right action corresponds to the geodesic frame flow of the hyperbolic orbifold $\G \backslash \mathbb{H}^n$. 
The geodesic frame flow factors the geodesic flow on the unit tangent bundle $T^1 \G \backslash \mathbb{H}^n$.
The latter identifies with the right action of $A$ on $\Gamma \backslash \mathrm{SO}(n,1)^0/\mathrm{SO}(n-1)$, where $\mathrm{SO}(n-1)$ is the stabilizer in $\mathrm{SO}(n)$ of a fixed unit vector in $T^1 \mathbb{H}^n$. 
The geodesic flow is topologically mixing on its non-wandering set\footnote{For example, topological mixing is equivalent to non-arithmeticity of the length spectrum by \cite{dalbo2000feuilletage}, which follow, for Zariski dense subgroup from  \cite{benoist2000proprietes} or \cite{kim2006length}. }.

Denote by $\Omega_G$ the preimage in $\Gamma \backslash \mathrm{SO}(n,1)^0$ of the non-wandering set of the geodesic flow. 
For convex cocompact subgroups, 
Winter \cite{winter2016exponential} and Sarkar-Winter \cite{sarkar2020exponential} proved exponential mixing for the push forward of the Bowen-Margulis-Sullivan (BMS) measure on the frame bundle.
Since this measure is supported in $\Omega_G$, these results imply topological mixing of the frame flow.

Under no other assumption for $\G$ than Zariski dense, Maucourant-Schapira \cite{MaucourantSchapira2019frameflows} proved that the frame flow is topological mixing on $\Omega_G$.
\vspace*{.2 cm}

For rank one (i.e. $\dim A=1$) locally symmetric spaces and discrete Zariski dense subgroup admitting a finite BMS measure,
Winter \cite{winter2015mixing} showed mixing for the frame flow.

\subsection{Main setting}
In this article, we focus on a higher rank semisimple Lie group $G$, meaning that $\dim A \geq 2$ and on an infinite covolume, discrete, Zariski dense subgroup $\G$ of $G$. 

Let $K$ be a maximal compact subgroup of $G$ for which the Cartan decomposition $KA^+K$ of elements of $G$ holds.
Denote by $M:=Z_K(A)$ the centralizer subgroup of $A$ in $K$.

For any $\theta \in \fa^{+}$, the \emph{nonnegative diagonal flow} $\phi_\theta^t$ corresponds to the right action by translation on $\Gamma \backslash G$ of $\exp (t \theta)$.
When $\theta \in \fa^{++} \setminus \lbrace 0\rbrace$, the flow $\phi_\theta^t$ is called \emph{positive diagonal}. 
Nonnegative diagonal flows $\phi_\theta^t$, where $\theta\in \fa^+$, induce right actions on $\G \backslash G/M$, so called \emph{Weyl chamber flows}.
They are called \emph{regular} when they are induced by positive diagonal flows.
The latter will play the same role in higher rank as the geodesic flow in the unit tangent bundle of the hyperbolic orbifold.
 
\subsection{Mixing of regular Weyl chamber flows} 
Conze-Guivarc'h \cite{conze2000limit} defined for $\mathrm{SL}(n,\R)$ and Zariski dense discrete subgroups a right $A$-invariant closed subset $\Omega \subset \Gamma \backslash G/M$ (cf. § 5.1 for a detailed construction). 
Their construction generalizes to all semisimple Lie groups without compact factors.
\begin{defin}\label{defin_omega}
We denote by $\Omega$ the smallest closed $A$-invariant subset of $\G \backslash G/M$ containing all periodic orbits of regular Weyl chamber flows and by $\Omega_G$ its preimage in $\G \backslash G$.
\end{defin}
The closed subset $\Omega$ is the analogue for Weyl chamber flows of the non-wandering set of the geodesic flow in the hyperbolic case.
With Glorieux \cite{dang_glorieux_2020}, we obtained the following necessary and sufficient mixing condition for regular Weyl chamber flows.

\begin{theorem}[ \cite{dang_glorieux_2020}]\label{theo-DG}
Let $G$ be a semisimple, connected, real linear Lie group, of non-compact type.
Let $\G$ be a Zariski dense, discrete subgroup of $G$.

A regular Weyl chamber flow $\phi_\theta^t$ is topologically mixing on $\Omega$ if and only if $\theta \in \overset{\circ}{\mathcal{C}}(\G)$.
\end{theorem}

The limit cone $\mathcal{C}(\G)$ was introduced by Benoist \cite{benoist1997proprietes}. For every Zariski dense $\G$, he proves that the limit cone is a closed, convex cone of $\fa^+$ of non-empty interior.
\begin{defin}\label{defin_limitcone}
Denote by $\lambda:G \rightarrow \fa^+$ the Jordan projection.
The \emph{limit cone} of $\G$ which is also called \emph{Benoist cone} $\mathcal{C}(\G)$, is the smallest closed cone of $\fa^+$ containing $\lambda(\G)$.
\end{defin}
Mixing ratio for regular Weyl chamber flow $\phi_\theta^t$, where $\theta$ lies in the interior of the limit cone, were obtained by Thirion \cite{thirion2009proprietes} for Ping-Pong groups,
Sambarino \cite{sambarino2015orbital} for Hitchin representations and Edwards-Lee-Oh \cite{edwards2020anosov} for Borel Anosov groups.

\subsection{Main result}
We study the topological dynamics of non-negative diagonal flows $(\Omega_G,\phi_\theta^t)$. We focus on its topological mixing properties.
Note that $\Omega_G$ is a right $AM$-invariant closed subset of $\G \backslash G$ and a principal $M$-bundle over $\Omega$, where $M$ is not necessarily connected.

Using a result of Guivarc'h-Raugi \cite{guivarc2007actions}, we partition $\Omega_G$ into finitely many $A$-invariant subsets that are dynamically conjugated to each other for nonnegative diagonal flows.

\begin{theorem}\label{theo_intro}
Let $G$ be a semisimple, connected, real linear Lie group, of non-compact type.
Let $\G$ be a Zariski dense, discrete subgroup of $G$.

Then there exists a normal subgroup of finite index $M_0 \lhd M_\Gamma \lhd M$ and
a partition $( \Omega_{[m]} )_{[m]\in M/M_\Gamma}$ of $\Omega_G$ such that 
\begin{itemize}
\item[(a)] every $\Omega_{[m]}$ is right $AM_\Gamma$-invariant and a principal $M_\G$-bundle over $\Omega$; 
\item[(b)] for all $\theta \in \fa^+$, the dynamical systems $\lbrace (\Omega_{[m]},\phi_\theta^t)  \rbrace_{[m]\in M/M_\Gamma}$ are conjugated to each other; 
\item[(c)] if $\theta \in \mathfrak{a}^{++}$ and $(\Omega_{[e_M]},\phi_\theta^t)$ is topologically mixing then
$\theta\in \overset{\circ}{ \mathcal{C}}(\Gamma)$ .
\end{itemize}
If furthermore $M_0$ is \emph{abelian} and
$\theta \in \mathfrak{a}^{++}$, then the converse of (c) is true:
\begin{itemize}
\item[(d)]$(\Omega_{[e_M]},\phi_\theta^t)$ is topologically mixing if and only if
$\theta\in \overset{\circ}{ \mathcal{C}}(\Gamma)$.
\end{itemize}
\end{theorem}
We expect that condition (d) holds in the general case, because Maucourant-Schapira \cite{MaucourantSchapira2019frameflows} proved topological mixing of the geodesic frame flow for $\mathrm{SO}(n,1)^0$ where $M=M_0=\mathrm{SO}(n-1)$.  
Condition (c) is a consequence of the joint work with Glorieux.

Observe that $M_0$ is abelian for example: split real semisimple Lie groups i.e. $\mathrm{SL}(n,\R)$, $\mathrm{Sp}(2n,\R)$, $\mathrm{SO}_0(p,p)$, $\mathrm{SO}_0(p,p+1)$;
and also for $\mathrm{SU}(p,p+1)$, $\mathrm{SU}(p,p)$, $\mathrm{SO}_0(p,p+2)$ and $\mathrm{SL}(n,\mathbb{C})$.
The closed, normal subgroup of finite index $M_\G$ of $M$ containing the connected component of the identity $M_0$ of $M$, is defined in 
Guivarc'h-Raugi \cite{guivarc2007actions} by
using the elliptic part of loxodromic elements of $\Gamma$.
It was also defined and studied in the appendix of  \cite{benoist2005convexes3}.
We call it the \emph{sign group} of $\G$.





Labourie \cite{labourie2006anosov} proved that $M_\G$ is trivial if $\G$ is the image of a Hitchin representations.
It thus follows from the above result that in this case, there are $2^{n-1}$ disjoint subsets in $\G \backslash \mathrm{PSL}(n,\mathbb{R})$ that share the same dynamical behavior for non-negative diagonal flows.
Consequently, positive diagonal flows are topologically mixing on any of these subsets if and only if they are parametrized by directions of the interior of the limit cone.

For Borel Anosov subgroups and independently, Lee-Oh \cite{lee2020ergodic} prove that there is an $A$-ergodic decomposition of every BMS measure into $AM_\Gamma$-semi-invariant and $A$-ergodic measures parametrized by $M/M_\Gamma$.
Any pair of such measures is the same up to right multiplication by elements of $M/M_\Gamma$, which concur with our result.

\subsection{Key ideas}

\subsection*{Bruhat-Hopf coordinates}

Denote by $\mathcal{F}^{(2)}$ the subset of transverse pairs in the Furstenberg boundary (cf. § 2.2) which identifies with $G/AM$ (cf. Proposition \ref{prop_transverse_pairs}).
Thirion \cite{thirion2007sous} generalized Hopf coordinates in higher rank by parametrizing point of $G/M$ with elements of $\mathcal{F}^{(2)} \times \fa$.
The left action of $G$ on $G/M$ reads using the Iwasawa cocycle $\sigma$ (cf. Definition \ref{defin-iwasawa-cocycle}) as follows
$$ g(\xi,\eta\; ; \; x) = (g \xi, g \eta \; ; \; \sigma(g,\xi)+ x).$$
The Weyl chamber flow reads by translating only the $\fa$ coordinate without changing the first two. 

Consider the set $ \lbrace G_s \rbrace_{s \in  \mathcal{S}}$ of maximal Bruhat cells of $G$.
For every $s \in \mathcal{S}$, we denote by $\mathcal{F}_s$ (resp. $\mathcal{F}_s^{(2)}$) the projection of $G_s$ in $\mathcal{F}$ (resp. $\mathcal{F}^{(2)}$).
 
In Section 3, we construct \emph{Bruhat-Hopf} coordinates $ \mathcal{H}_s: G_s \rightarrow \mathcal{F}_s^{(2)} \times AM$ that extend Hopf coordinates (cf. Definition \ref{defin_bh-coord}, \ref{defin-unip_bruhat-section}, Proposition \ref{prop_bh-coord}).
Note that they differ from coordinates coming from the unique Bruhat decomposition of $N^-MAN$ or their translate of the form $hN^-MAN$, where $h \in G$.
The projection $G \rightarrow G/M$ reads for all $s\in \mathcal{S}$ by preserving the coordinates in $\mathcal{F}^{(2)}$ and projecting the $AM$-coordinates to $\fa$.
The right translation by $AM$ on $G$ reads for all $(\check{\xi}, \xi ; u)_s \in \mathcal{F}_s^{(2)} \times AM$ and $x \in AM$ as $(\check{\xi}, \xi ; ux)_s$.

The left action of $G$ on itself reads in this family of Bruhat-Hopf coordinates $(\mathcal{H}_s)_{s\in \mathcal{S}}$ equivariantly in the coordinates in $\mathcal{F}^{(2)}$ and via left multiplication by the \emph{signed Iwasawa cocycles} $(\beta_{s',s})_{s,s' \in \mathcal{S}}$ (cf. Definition \ref{defin_cocycle}) of domains in $G \times \mathcal{F}$ and codomains in $AM$.  
They extend (cf. Proposition \ref{prop_bh-coord}) the Iwasawa cocycle in the sense that for all $\xi \in \mathcal{F}_s$ and $g \in G$ such that $g\xi \in \mathcal{F}_{s'}$, then 
$ \beta_{s',s}(g,\xi) \in \exp(\sigma(g,\xi)) M.$
We prove that the signed cocycles $(\beta_{s,s})_{s\in \mathcal{S}}$ are all cohomologous (cf. Fact \ref{fait_cocycle-transfer}) for the \emph{transition maps} $\mathscr{T}_{s,s'}: \mathcal{F}_s\cap \mathcal{F}_{s'} \rightarrow AM$ of Definition \ref{defin-transfer}. 

Furthermore, Bruhat-Hopf coordinates induce local coordinates of $K$ in $\mathcal{F} \times M$ by removing the second coordinate and projecting in $M$ the third one.

Likewise, the reader can check that Bruhat-Hopf coordinates induce local coordinates on $G/N$, $G/A$ and $G/MN$. 

\subsection*{The elliptic part of loxodromic elements}
Elements of $G$ whose Jordan projection is in the positive Weyl chamber are called \emph{loxodromic}.
Denote by $G^{lox}$ and $\G^{lox}$ the subset of loxodromic elements of the respective groups.
Loxodromic elements (see §4) have trivial unipotent parts and are conjugated to elements in $MA^{++}$.
The part in $A^{++}$, corresponding to the hyperbolic part, is given by the Jordan projection.
In \cite{benoist1996actions}, \cite{benoist1997proprietes} and \cite{benoist2000proprietes}, Benoist defines $(r,\varepsilon)$-loxodromic elements (see Definition \ref{defin-r-eps}) and obtains estimates for the Jordan projection of generic products of $(r,\varepsilon)$-loxodromic elements. 
We show that their elliptic part satisfy similar estimates.
\vspace*{.2 cm}

The elliptic part of a loxodromic element is conjugated to an element of $M$ which is defined up to conjugacy by $M$.
Therefore, the latter is only well defined when $M$ is abelian, in which case one can extend the Jordan projection from $G^{lox}$ to $\fa^{++} \times M$. 
Bruhat-Hopf coordinates gives a framework to solve this technical difficulty in the general case.

Fix a loxodromic element $g$ and denote by $g^+$ (resp. $g^-$) its attracting (resp. repelling) fixed point in $\mathcal{F}$ and by $\mathsf{b}(g^-)$ the basin of attraction of $g^+$ (cf. Proposition \ref{prop-lox-bassin}).
Starting from the formula $\sigma(g,g^+)=\lambda(g)$ satisfied by loxodromic elements, we define a multiplicative and signed Jordan projection for $g$.
For every $s\in \mathcal{S}$ such that $g^+ \in \mathcal{F}_s$, we set
$\mathscr{L}_s(g):= \beta_{s,s}(g,g^+).$
It is the unique element in $\exp(\lambda(g))M$ such that there is an element $h_s\in G_s$ unique up to right multiplication by $A$ such that $h_s^{-1} g h_s = \mathscr{L}_s(g)$.

Using the continuous maps $\mathscr{R}_{s',s}$ given in Definition \ref{defin-ratio}, we obtain an exact formula.  

\begin{prop}[\ref{prop-magic-cocycle} below]\label{prop-intro-cocycle}
\hypG

Then for all loxodromic element $g \in G^{lox}$, all integer $n\geq 1$ and $\xi \in \mathsf{b}(g^-)$, for any suitable $s_0,s_1,s_2 \in \mathcal{S}$  such that $(\xi,g^+,g^n\xi)\in \mathcal{F}_{s_0} \times \mathcal{F}_{s_1} \times \mathcal{F}_{s_2}$
$$
\beta_{s_2,s_0}(g^n,\xi)= \mathscr{R}_{s_1,s_2}(g;g^n\xi)^{-1} \mathscr{L}_{s_1}(g)^n \mathscr{R}_{s_1,s_0}(g;\xi).
$$
\end{prop}

We estimate the elliptic part of generic products of $(r,\varepsilon)$-loxodromic elements.
In order to do that, we introduce a family of constants $ \lbrace \delta_{r,\varepsilon} \; \vert \; 0 < \varepsilon \leq r \rbrace$ (cf. Definition \ref{defin_d-r-eps}) such that for all $r >0$, they satisfy $\lim_{\epsilon \rightarrow 0} \delta_{r,\epsilon} =0$ (cf. Proposition \ref{prop_equicontinu}).

\begin{prop}[\ref{prop_crucial} below]\label{prop_lox_ell}
\hypG
For all $r>0$ and $\varepsilon \in (0,r]$ and every family $g_1,...,g_l \in G$ of $(r,\varepsilon)$-loxodromic elements such that
\begin{itemize}
\item[$\star$] $r \leq \frac{1}{6} \mathrm{d} \big( \lbrace g_{i-1}^+,g_i^+ \rbrace, \partial \mathsf{b}(g_i^-)  \big) $ for all $1 \leq i \leq l$ with the convention $g_0=g_l$.
\end{itemize}
For all family $(s_i)_{0 \leq i \leq l} \subset \mathcal{S}$ such that 
\begin{itemize}
\item[$\star \star$]  $\mathcal{F}_{s_i} \supset \mathcal{V}_r(\partial \mathsf{b}(g_i^-))^\complement$ for every $1 \leq i \leq l$ and $\mathcal{F}_{s_0}\supset \mathcal{V}_{\varepsilon}(\partial \mathsf{b}(g_1^-))^\complement$.
\end{itemize}
Then for all integers $n_1,...,n_l \geq 1$,
the element $g_l^{n_l}...g_1^{n_1}$ is $(2r,2 \varepsilon)$-loxodromic with attracting (resp. repelling) point in $B(g_l^+,\varepsilon)$  (resp. $B(g_1^-,\varepsilon)$) and its extended Jordan projection satisfies
$$ \mathscr{L}_{s_l}(g_l^{n_l}... g_1^{n_1}) \in 
\mathscr{L}_{s_l}(g_l^{n_l} )
\mathscr{R}_{s_l,s_{l-1}}(g_l, g_{l-1}^+) ...
\mathscr{L}_{s_1}(g_1^{n_1})
\mathscr{R}_{s_1,s_l}(g_1, g_l^+)
 B( e_{AM}, 2l \delta_{r,\varepsilon} ) .$$
	\end{prop}

\subsection*{Decorrelation}
Denote by $M^{ab}$ the abelianization of $M$.
We define an abelianized Jordan projection for loxodromic elements $\mathscr{L}^{ab}: G^{lox} \rightarrow A^{++}M^{ab}$ using the previous local Jordan projections $\mathscr{L}_s$.
The number of connected components of $M^{ab}$ reached by the subset $\mathscr{L}^{ab}(\G^{lox}) $ suffices to understand $M_\G$.
Indeed, its abelianized $M_\G^{ab}$ is the
subgroup of $M^{ab}$ generated by the projection to $M^{ab}$ of $\mathscr{L}^{ab}(\G^{lox})$.
Thanks to Guivarc'h-Raugi \cite[Theorem 6.4]{guivarc2007actions} we deduce that the  
subgroup generated by $\mathscr{L}^{ab}(\G^{lox}) $ is dense in $AM_\Gamma^{ab}$.
Guivarc'h-Raugi also give a classification of $\G$-invariant minimal subsets of $K$.
We rewrite their result using Bruhat-Hopf coordinates of $K$ in Theorem \ref{theo_GR_Gamma_orbite_K} and define the invariant subsets $\Omega_{[m]}$ through their universal cover $\widetilde{\Omega}_{[m]}$ in $G$.
\vspace*{.2 cm}

Denote by $L(\G) \subset \mathcal{F}$ the limit set of $\G$ and by $L^{(2)}(\G):= L(\G) \times L(\G) \cap \mathcal{F}^{(2)}.$ 
The universal cover $\widetilde{\Omega}_G$ has Bruhat-Hopf coordinates $L^{(2)}(\G) \times AM$.

Without loss of generality, by using the joint work with Glorieux \cite{dang_glorieux_2020}, it suffices to prove the \emph{decorrelation Proposition} \ref{prop_decorrelation} i.e. that  there exists $(\xi_1,\check{\xi}_1 ) \in L^{(2)}(\G)$ such that for every $x \in AM$ and small $\delta>0$, the orbit $\G (\xi_1, \check{\xi}_1 \; ; \; x)_{\check{c}_1}$ is $\delta$-dense in an $M_\G$-orbit of the form $(\xi_1, \check{\xi}_1 \; ; \; y_{\delta}x M_\G )_{c_1}$ (for suitable $\check{c}_1,c_1 \in \mathcal{S}$).

The first step (Lemma \ref{lem_decorrelation_discret}) is to reach all connected components of $M_\G$ by the
left action of finitely many $(r,\varepsilon)$-loxodromic elements of $\G$ of attracting point close to $\xi_1$.
It does not use that $M$ abelian.

In the second step (Lemma \ref{lem_decorM}) we construct $(r,\varepsilon)$-loxodromic elements $\g_1,...,\g_l \in \G$ that satisfy the hypothesis of Proposition \ref{prop_lox_ell} and such that
$\mathscr{L}^{ab}( \lbrace \g_l^{n_l}...\g_1^{n_1} \; \vert \; n_1,..,n_l \geq 1 \rbrace )$
is $\delta$-dense in an $M_0$-invariant set that projects to $\log \pi_A( y_\delta x)+ \mathcal{C}_0$, where $\mathcal{C}_0 \subset \fa^{++}$ is a closed convex cone of non-empty interior.
We rely on density of squares in $M_0$, as well as density lemmata deduced from the assumption that $M_0$ is abelian.

Finally, we deduce the decorrelation with an overlapping cone argument.

\subsection{Organization of the paper}
In Section 2 we recall the classical Iwasawa, Bruhat decompositions of Lie groups and
characterize the transverse points in the Furstenberg boundary.
Section 3 is dedicated to the construction of Bruhat-Hopf coordinates.
In Section 4, using Bruhat-Hopf coordinates, we estimate the elliptic part of products of generic loxodromic elements.
In Section 5, we define the subgroup $M_\G$, the $\G$-invariant subsets of $G$ and prove Theorem \ref{theo_intro} (a)(b).
Section 6 is dedicated to the proof of decorrelation.
In Section 7 we prove the necessary and sufficient condition for topological mixing when $M_0$ abelian.
In the appendix, we prove the density lemmata.

\subsection*{Relation to other works}
Sections 2, 5, 6, 7 and the Appendix can be found in french in the author's PhD thesis \cite{dang2019these}.
Sections 3 and 4 improve the thesis's construction of Bruhat-Hopf coordinates and its estimates of the elliptic and hyperbolic parts of products of loxodromic elements.

Bruhat-Hopf coordinates were independently studied by Lee-Oh \cite{lee2020ergodic}.

\section*{Acknowledgments}
The author acknowledges funding by the Deutsche Forschungsgemeinschaft (DFG, German Research Foundation) – 281869850 (RTG 2229).
She would like to thank O. Glorieux, F. Maucourant, B. Schapira and B. Pozzetti for encouragements, fruitful discussions and comments.
 
\section{Background}
A classical reference for this section is \cite{helgason1978differential}.
Let $K$ be a maximal compact subgroup of $G$.
Denote by $\mathfrak{g}$ (resp. $\mathfrak{k}$) the Lie algebra of $G$ (resp. $K$).
Consider a Cartan decomposition $\mathfrak{g}=\mathfrak{k} \oplus \mathfrak{p}$.
Let $\mathfrak{a} \subset \mathfrak{p}$ be a \emph{Cartan subspace} i.e. a maximal abelian subspace of $\mathfrak{p}$ for which the adjoint endomorphism of every element is semisimple.
Denote by $\mathfrak{m}$ the centralizer of $\mathfrak{a}$ in $\mathfrak{k}$.

For every linear form $\alpha \in \mathfrak{a}^*$, set 
$\mathfrak{g}_\alpha := \lbrace v \in \mathfrak{g} \; \vert \; \forall u\in \mathfrak{a}, \; [u,v]=\alpha(u)v \rbrace.$
Note that $\mathfrak{g}_0= \mathfrak{m} \oplus \mathfrak{a}$.
The set of \emph{restricted roots} is given by
$\Sigma := \lbrace \alpha \in \mathfrak{a}^* \setminus 0 \; \vert \; \mathfrak{g}_\alpha \neq 0 \rbrace.$
By simultaneous diagonalisation over the abelian family of endomorphisms $ad(\mathfrak{a})$, we deduce the decomposition $\mathfrak{g}= \mathfrak{g}_0 \oplus_{\alpha \in \Sigma} \mathfrak{g}_\alpha$.
Note that $\Sigma$ is a finite set.
Let us now choose a \emph{positive Weyl chamber} of $\mathfrak{a}$ i.e. a connected component of $\mathfrak{a} \setminus \cup_{\alpha \in \Sigma} \ker(\alpha)$.
Denote the closed positive Weyl chamber by $\mathfrak{a}^+$ and $\mathfrak{a}^{++}$ its interior.
The set of \emph{positive roots}, denoted by $\Sigma^+$, is the subset of restricted roots which take positive values in the positive Weyl chamber.
This choice allows to define two particular nilpotent subalgebras $\mathfrak{n} =\oplus_{\alpha\in \Sigma^+} \mathfrak{g}_\alpha$ and $\mathfrak{n}_- =\oplus_{\alpha\in \Sigma^+} \mathfrak{g}_{-\alpha}$.

Finally, denote by $A:=\exp(\fa)$ the maximal split torus or Cartan subgroup, $A^+:= \exp(\fa^+)$ the closed positive Weyl chamber, $A^{++}:=\exp(\fa^{++})$ its interior, $N:=\exp(\mathfrak{n})$ (resp. $N^-:=\exp(\mathfrak{n}_-)$) the positive (resp. negative) maximal unipotent subgroups and $M$ the centralizer of $A$ in $K$, of Lie algebra $\mathfrak{m}$.
By definition, $A$ normalizes $N$ and $N^-$.
Furthermore, for all $a\in A^{++}$ and $h_\pm\in N^\pm$ the following convergences hold 
\begin{equation}\label{equ_unipotents}
a^{-n}h_{\pm}a^n  \underset{\pm \infty}{\longrightarrow} e_G .
\end{equation}

\subsection{Furstenberg boundary}
By Iwasawa decomposition ( cf.
\cite[Chapter IX, Thm 1.3 ]{helgason1978differential})
$G=KAN$ and $G=KAN^-$
and the maps (with the convention that $N^+=N$)
\begin{align*}
K \times A \times N^\pm &\longrightarrow G \\
(k,a,n) & \longmapsto kan
\end{align*}
are diffeomorphisms.
Denote by $g  \mapsto \big( k_{\mathcal{I}\pm}(g), a_{\mathcal{I}\pm}(g), u_{\mathcal{I}\pm}(g) \big) \in K \times A \times N^\pm $
the respective inverse diffeomorphisms.
Note\footnote{using Jacobi identity} that $[\mathfrak{g}_{0},\mathfrak{g}_\alpha] \subset \mathfrak{g}_{\alpha}$ for all $\alpha \in \Sigma_+$. 
Hence $\mathfrak{m}\oplus \mathfrak{a}\oplus \mathfrak{n}$ and $\mathfrak{m}\oplus \mathfrak{a}\oplus \mathfrak{n}_-$ are Lie subalgebras of $\mathfrak{g}$.
Consequently $MAN$ and $MAN^-$ are closed subgroups of $G$.
\begin{defin}\label{defin_furstenberg}
The \emph{Furstenberg boundary} is defined by $\mathcal{F}:= G/MAN$.
Denote by $k_\iota\in K$ a representative of the element in the Weyl group such that $\mathrm{Ad}(k_\iota)\mathfrak{a}^+ = -\mathfrak{a}^+$.  
Set $\eta_0:= MAN$ and $\check{\eta}_0:= k_\iota \eta_0$.
\end{defin}

The map
$k \in K \mapsto k \eta_0 \in \mathcal{F}$ is surjective and equivariant for the left action of $K$.
Furthermore, the stabilizer of $\eta_0$ is the closed subgroup $M$.
Therefore, we deduce an identification of $K/M$ with the Furstenberg boundary. 

Let us sketch the construction of a $K$-invariant Riemannian distance on $K$.
Start from a scalar product on $\mathfrak{k}$. 
Since $K$ is a compact subgroup, its Haar measure is finite.
By averaging the scalar product on $\mathfrak{k}$ along the Haar measure on $K$ for the adjoint action, we obtain an $Ad (K)$-invariant scalar product and norm on $\mathfrak{k}$.
Using the left action of $K$, we transport them on every tangent space and obtain a left $K$-invariant metric which is also invariant by conjugation.
Hence $K$ is endowed with an invariant Riemannian metric.
Its induced Riemannian distance is thus $K$-invariant. 
\begin{defin}\label{defin_dist-furstenberg}
Let $d_K$ be a $K$-invariant Riemannian distance on $K$. 
For every $\xi, \eta \in \mathcal{F}$ for any choice of representatives $k_\xi, k_\eta \in K$ such that $k_\xi \eta_0=\xi$ and $k_\eta \eta_0=\eta$, we consider the induced left $K$-invariant distance in $\mathcal{F}$
$$\mathrm{d}(\xi,\eta):= d_K(k_\xi M, k_\eta M).$$
\end{defin}
Let us define the Iwasawa cocycle. 
\begin{defin}\label{defin-iwasawa-cocycle}
For all $g \in G$ and $\xi \in \mathcal{F}$, we denote by $\sigma(g,\xi)$ the unique element\footnote{because $M$ normalises $N$, this element does not depend on the choice of the representative in $K$ of $\xi$.} in $\mathfrak{a}$ such that for all $k_\xi \in K$ such that $k_\xi \eta_0=\xi$,
$$g k_\xi \in K \exp(\sigma(g,\xi)) N. $$
The map $\sigma : G \times \mathcal{F} \rightarrow \mathfrak{a}$ is the \emph{Iwasawa cocycle}.
\end{defin}

\subsection{Transverse pairs in the Furstenberg boundary}
The following subset of $\mathcal{F} \times \mathcal{F}$ is a higher rank analogue to the set of pair of points in the geometric boundary of the hyperbolic plane $ \mathbb{H}^2 $ that parametrize oriented geodesics.
It also identifies for $\mathrm{SL}(n,\mathbb{R})$ with the space of transverse complete flags of $\R^n$.

\begin{defin}\label{defin_opposite_Furstenberg}
The subset of \emph{ordered transverse pairs} of $\mathcal{F} \times \mathcal{F}$ is defined by
$$\mathcal{F}^{(2)}:= \lbrace (g \eta_0, g \check{\eta}_0) \; \vert \; g \in G \rbrace.$$
Since $k_\iota$ is an involution, $(\check{\eta}_0, \eta_0)$ is also an ordered transverse pair.
Consequently, we say that $\xi , \eta \in \mathcal{F}$ are \emph{transverse} if any of the ordered pairs $(\xi, \eta)$ of $(\eta,\xi)$ are transverse.
\end{defin}
Denote by  $W:=N_K(A)/Z_K(A)$ the Weyl group of $G$.
We choose for every $w\in W$ a representative $k_w \in N_K(A)$. 
Then by Bruhat decomposition \cite[Chapter IX, Thm 1.4 ]{helgason1978differential},
$$G= \sqcup_{w \in W} B k_w B $$
where $B=MAN$.
Note that $N^-= k_\iota N k_\iota^{-1} $ and that 
$G= \sqcup_{w \in W} k_\iota B k_w B$, meaning that $N^- MAN$ is a cell in the Bruhat decomposition of $G$.
\begin{corol}[Chapter IX, Cor. 1.9 \cite{helgason1978differential}]\label{corol_bruhat}
\hypG
Then the map 
\begin{align*}
N^- &\longrightarrow N^- \eta_0 \\
n_- & \longmapsto n_- \eta_0
\end{align*}
is a diffeomorphism, its image is an open submanifold of $\mathcal{F}$ and its complement is a finite union of disjoint submanifolds of stricly smaller dimensions.
\end{corol}
It thus turns out that $N^- MAN$ is a maximal cell for the Bruhat decomposition.
We describe below the subset of transverse pairs in the Furstenberg boundary and include a proof for completeness.
\begin{prop}\label{prop_transverse_pairs}
\hypG
Then the following holds,
\begin{itemize}
\item[(i)] the set of transverse points to $\check{\eta}_0$ is $N^- \eta_0$,
\item[(ii)] for all $\eta, \xi \in \mathcal{F}$ and $k_\eta, \check{k}_\xi \in K$ such that $k_\eta \eta_0=\eta$ and $\check{k}_\xi \check{\eta}_0= \xi$,
$$ (\eta,\xi) \in \mathcal{F}^{(2)} \Longleftrightarrow \check{k}_\xi^{-1} k_\eta \in N^- MAN, $$
\item[(iii)] for all $\xi \in \mathcal{F}$ and $\check{k}_\xi \in K$ such that $\check{k}_\xi \check{\eta}_0= \xi$, the set of transverse points to $\xi$ is $\check{k}_\xi N^- \eta_0$.
\item[(iii')] for all $\xi \in \mathcal{F}$ and $k_\xi \in K$ such that $k_\xi \eta_0= \xi$, the set of transverse points to $\xi$ is $k_\xi N \check{\eta}_0$. 
\end{itemize}
Furthermore, the $G-$equivariant map
\begin{align*}
G/AM  & \longrightarrow \mathcal{F}^{(2)} \\
g AM & \longmapsto (g \eta_0, g \check{\eta}_0)
\end{align*}
is a diffeomorphism.
\end{prop}
\begin{proof}
(i) First remark that $N^-(\eta_0,\check{\eta}_0)= (N^- \eta_0, \check{\eta}_0)$. 
Let us now prove the converse i.e. that any point transverse to $\check{\eta}_0$ must be in $N^- \eta_0$.
Let $g\in G$ such that $(g\eta_0,\check{\eta}_0) \in \mathcal{F}^{(2)}.$
Then by definition, there exists $h\in G$ such that
$$(g \eta_0,\check{\eta}_0)=h(\eta_0,\check{\eta}_0).$$
On one hand $h\check{\eta}_0=\check{\eta}_0$,
hence $h\in \mathrm{Stab}(\check{\eta}_0)= k_\iota MAN k_\iota^{-1}$.
Since $N^-= k_\iota N k_\iota^{-1} $ and $MA$ is invariant by conjugation by $k_\iota$, we deduce that $$h\in MAN^-.$$
On the other hand $g \eta_0=h\eta_0$, hence $h^{-1}g\in \mathrm{Stab}(\eta_0)=MAN.$
Thus $$g\in h MAN \subset MAN^-MAN.$$
Since $MA$ normalizes $N^-$, we deduce that $g\in N^- MAN$.
Hence $g\eta_0 \in N^- \eta_0$.

(ii) It follows from (i) and by noticing that the pair
$ (k_\eta \eta_0,\check{k}_\xi \check{\eta}_0) \in \mathcal{F}^{(2)}$ if and only if $ (\check{k}_\xi^{-1}k_\eta \eta_0,\check{\eta}_0) \in \mathcal{F}^{(2)}.$ 

(iii) It follows from (ii) since $\check{k}_\eta( N^-\eta_0,\check{\eta}_0 ) \in \mathcal{F}^{(2)}$.

For the last statement, remark first that $G$ acts transitively on
 $\mathcal{F}^{(2)}$.
Furthermore 
$$\mathrm{Stab}_G(\eta_0, \check{\eta}_0)=MAN \cap MAN^- = AM.$$
We thus deduce the $G$-equivariance and bijectivity of the map 
\begin{align*}
G/AM & \longrightarrow \mathcal{F}^{(2)} \\
g AM & \longmapsto (g\eta_0, g \check{\eta}_0).
\end{align*}  
The left action of $G$ on the Furstenberg boundary $\mathcal{F}=G/MAN$ is differentiable
and so is its action on $\mathcal{F}\times \mathcal{F}$.
Thus, the map $g \mapsto (g\eta_0, g\check{\eta}_0)$ is differentiable.
The kernel of the differential in $e_G$ of the map $g \mapsto (g\eta_0, g\check{\eta}_0)$
contains $\mathfrak{m} \oplus \fa$.
Since the maps 
$N^- \rightarrow N^- \eta_0$ and $N \rightarrow N \check{ \eta}_0$ are diffeomorphisms, the 
differential in $e_G$ of $g \mapsto (g\eta_0, g\check{\eta}_0)$ is surjective from $\mathfrak{g}$ to $\mathfrak{n}_- \oplus \mathfrak{n}_+$.
By Bruhat decomposition in the Lie algebra
$\mathfrak{g}= \mathfrak{n}_- \oplus \mathfrak{m} \oplus \fa \oplus \mathfrak{n},$
we deduce that the kernel of the differential in 
$e_G$ of $g \mapsto (g\eta_0, g\check{\eta}_0)$ 
is equal to $\fa \oplus \mathfrak{m}$.
Thus, the map $G/AM \rightarrow \mathcal{F}^{(2)}$ is a local diffeomorphism in $AM$.
Finally, by transitivity of the left $G$ action on $G/AM$, we deduce that it is a diffeomorphism.
\end{proof}
We parametrize the maximal Bruhat cells of the Furstenberg boundary.
\begin{defin}\label{defin_bruhat-cell}
Let $\check{\xi} \in \mathcal{F}$, then for any representative $h(\check{\eta}) \in G$ such that $\check{\eta}=h(\check{\eta})\check{\eta}_0$,
we denote by $\mathsf{ b}( \check{\eta} ):= h(\check{\eta}) N^- \eta_0 $ the Bruhat cell opposite to $ \check{\eta}$.
\end{defin}
Thanks to the previous Proposition, the representative $h(\check{\eta})\in G$ is chosen up to right multiplication by $MAN^-$.
Remark that $\mathsf{b}(\eta_0)=N \check{\eta}_0$ and $\mathsf{b}(\check{\eta}_0)=N^- \eta_0$.
Using this notation, the set of Bruhat cells of $\mathcal{F}$ is naturally endowed with a left action of $G$ which satisfies $h \mathsf{b}(\check{\eta_0}):= \mathsf{b}(h\check{\eta}_0)$ for all $h \in G$.

\section{Bruhat-Hopf coordinates}
In his thesis, Thirion \cite[Chapter 8 §8.G.2]{thirion2007sous} introduced Hopf coordinates for $\mathrm{SL}(n,\mathbb{R})/M$. 
His construction generalizes to every semisimple Lie group without compact factors.
It is defined by 
\begin{align*}
G/M & \longrightarrow \mathcal{F}^{(2)} \times \fa \\
hM & \longmapsto (h \eta_0, h \check{\eta}_0\; ; \; \sigma(h,\eta_0)).
\end{align*}
The commuting left action of $G$ and right action of $A$ on $G/M$ read in those coordinates for all $(g,\theta,t) \in G \times \fa \times \mathbb{R}$ and $(\xi, \check{\xi} \; ; \; x) \in \mathcal{F}^{(2)} \times \fa$ as follows.
$$\tilde{\phi}_\theta^t \big( g(\xi, \check{\xi} \; ; \; x) \big)=(g \xi, g \check{\xi} \; ; \; \sigma(g,\xi)+ x+ t \theta).$$
The respective projections $G/M \rightarrow \mathcal{F}$ and $G/M \rightarrow G/AM$ read as the projection to the first coordinate in $\mathcal{F}$ and by removing the coordinate in $\fa$.

In this section, we extend locally and equivariantly (for the left action of $G$ and right action of $A$) Hopf coordinates to $G$.

Any local trivialisation of $G \rightarrow G/AM$ provides local coordinates in $\mathcal{F}^{(2)} \times AM$ that are equivariant for the right action of $A$.
The restricted left $G$-action provides a local $AM$-cocycle.
In general, neither these cocycles extend the Iwasawa cocycle nor will those coordinates locally extend Hopf coordinates.
We construct families of trivialisations of $G \rightarrow G/AM$ starting from those of $G \rightarrow \mathcal{F}$, for which these cocycles generalize the Iwasawa cocycle and obtain commutative diagrams with the Hopf coordinates.

In §3.1, using Bruhat decomposition, we construct from cross-sections of $G \rightarrow \mathcal{F}$, local coordinates in $\mathcal{F}^{(2)} \times AM$ of $G$ (cf.  Definition \ref{defin_bh-coord}).
We set notations for the rest of the article and in Definition \ref{defin_compatible-sections}, define covering families of cross-sections and of the same type (i.e. that are translates of one another by $G$-action).

In §3.2, we set notations for the transition functions between the trivializations in Definition \ref{defin-transfer}.
In Proposition \ref{prop_tranfer} we express classical properties of these functions with respect to our choice of notation.
Note that the cross-section parameters follow a Chasles relation.

In §3.3, for every family of differentiable cross-sections $(s_i)_{i\in I}$ of $G \rightarrow \mathcal{F}$ whose domain cover $\mathcal{F}$, we read in those coordinates the left action of $G$ on itself.
The behavior is the same as for Hopf coordinates for the first two coordinates in $\mathcal{F}^{(2)}$.
We define $AM$-valued functions in Definition \ref{defin_cocycle} of domain in $G \times \mathcal{F}$.
We prove in Proposition \ref{prop_cocycle} that those functions are cocycles that encode the information in $AM$ for the left action of $G$ on itself.
We deduce that the information contained in the second and third coordinate in $\mathcal{F}^{(2)} \times AM$ are not needed when one reads the left action of $G$.
In Fact \ref{fait_cocycle-transfer}, we extend to the cocyle the Chasles type relations we previously had with the transition functions.

In §3.4, we give a sufficient condition on a differential cross-section of $G \rightarrow \mathcal{F}$ for which the associated coordinates of $G$ extend Hopf coordinates.
We prove in Proposition \ref{prop_bh-coord} that when the cross-sections $(s_i)_{i \in I}$ take value in $K$, the signed multiplicative Iwasawa cocycles $(\beta_{s_i,s_j})_{i,j\in I}$ defined in the third paragraph generalize the Iwasawa cocycle.
We obtain an equivariant and commutative diagram with Hopf coordinates.

In §3.5, we prove in Proposition \ref{prop_BH_K} that local coordinates of $G$ that extends Hopf coordinates provide local coordinates of $K$ that take value in $\mathcal{F} \times M$.
Furthermore, the map $k_{\mathcal{I}}: G \rightarrow K$ reads in those coordinates by keeping the first coordinate in $\mathcal{F}$ and projecting the last one in $M$.

In the last paragraph, using Bruhat decomposition and Iwasawa decomposition, we construct two families of cross-sections of $G \rightarrow \mathcal{F}$ defined on Bruhat cells of $\mathcal{F}$: \emph{unipotent} and \emph{compact Bruhat sections} in 
Definition \ref{defin-unip_bruhat-section}
We define \emph{Bruhat-Hopf coordinates} as the local extensions of Hopf coordinates given by Proposition \ref{prop_bh-coord} with respect to the compact Bruhat sections.
In Proposition \ref{prop_transfer-bruhat} we
parametrize these cross-sections.


\subsection{Local trivialisations}
Let $s$ be a non-trivial cross-section of the $MAN$-bundle $G \rightarrow \mathcal{F}$, we denote by $\mathcal{F}_s$ its domain, $\mathcal{F}_s^{(2)}:= (\mathcal{F}_s \times \mathcal{F}) \cap \mathcal{F}^{(2)}$ the subsets of ordered transverse pairs of first coordinate in $\mathcal{F}_s$. 
We explicit the trivialization of $G$ given by this cross-section.
\begin{fait}\label{fait_bh-coord}
\hypG
Let $s$ be a differentiable cross-section of $G \rightarrow \mathcal{F}$.

Then the two maps below are diffeomorphisms. 
\begin{align*}
\mathcal{F}_s \times N \times AM 
& \longrightarrow 
s(\mathcal{F}_s) N AM \subset G \\
(\xi, u, x)
&\longmapsto
s(\xi)ux.
\end{align*}
\begin{align*}
\mathcal{F}_s \times N \times AM 
& \longrightarrow 
\mathcal{F}_s^{(2)} \times AM \\
(\xi, u, x)
&\longmapsto
(\xi, s(\xi)u \check{\eta}_0\; ; \; x)_s.
\end{align*}
\end{fait}
\begin{proof}
By hypothesis, the map $s: \mathcal{F}_s \rightarrow G$ is a cross-section and by Iwasawa decomposition in $NAM$, we deduce that the first map is a diffeomorphism.

By Proposition \ref{prop_transverse_pairs} (iii')
for every $\xi \in \mathcal{F}_s$, the set of transverse points to $\xi$ is $s(\xi)N \check{\eta}_0$.
Hence the map 
\begin{align*}
\mathcal{F}_s\times N & \longrightarrow \mathcal{F}_s^{(2)} \\
(\xi, u) & \longmapsto (\xi, s(\xi)u \check{\eta}_0)
\end{align*}
is a diffeomorphism. 
Consequently, the second map is a diffeomorphism.
\end{proof}

\begin{defin}\label{defin_bh-coord}
For any differentiable cross-section $s$ of $G \rightarrow \mathcal{F}$, we denote by $\mathcal{B}_s$ the associated differentiable trivialization
\begin{align*}
\mathcal{B}_s: s(\mathcal{F}_s)N AM \subset G & \longrightarrow
\mathcal{F}_s^{(2)} \times AM \\
g= s(\xi)ux & \longmapsto (g \eta_0, g \check{\eta}_0 \; ;\; x)_s.
\end{align*}

When $s$ is compact valued i.e. a cross-section of $K \rightarrow \mathcal{F}$, the same map is denoted $\mathcal{H}_s$. 
\end{defin}
In order to parametrize every element of $G$ in such coordinates, we construct families of differentiable cross-sections whose domain cover $\mathcal{F}$.
For all $g \in G$ and any cross-section $s : \mathcal{F}_s \rightarrow G$, we define the left translate by
\begin{align*}
g\cdot s : g \mathcal{F}_s & \longrightarrow G \\
\xi & \longmapsto g  s(g^{-1}\xi).
\end{align*} 
This provides a left $G$ action on the space of cross-sections of $G \rightarrow \mathcal{F}$.
For any $b \in MAN,$ we define the cross-section
\begin{align*}
s.b :  \mathcal{F}_s & \longrightarrow G \\
\xi & \longmapsto  s(\xi)b.
\end{align*} 
\begin{defin}\label{defin_compatible-sections}
A family of cross-section $(s_i)_{i\in I}$ of the bundle $G \rightarrow \mathcal{F}$ is 
\emph{covering} when the family of domains $\lbrace\mathcal{F}_{s_i}\rbrace_{i\in I}$ covers $\mathcal{F}$ i.e.
$$\mathcal{F} \subset \cup_{i\in I} \mathcal{F}_{s_i}.$$
The family $(s_i)_{i\in I}$ is of the same type if for any $i,j \in I$ there exists $g_{ij} \in K$ such that 
$$ s_i=g_{ij} \cdot s_j.$$
\end{defin}
Using that $K$ acts transitively on $\mathcal{F}$ and the compacity of the latter, one can construct finite families of differentiable cross-sections of the same type that are covering.
We provide two such families in Definition \ref{defin-unip_bruhat-section}. 
\subsection{Transition functions}
First we explicit the transition functions between $\mathcal{B}_s$ parametrizations. 
In Proposition \ref{prop_tranfer} we express classical properties of these functions with respect to our choice of notation.
Note that the cross-section parameters follow a Chasles relation.
\begin{fait}\label{fait_transfer}
\hypG
Let $s$ and $s'$ be two differentiable cross-sections of $G \rightarrow \mathcal{F}$ such that $\mathcal{F}_s \cap \mathcal{F}_{s'} \neq \emptyset$.
Then for all $\xi \in\mathcal{F}_s \cap \mathcal{F}_{s'},$  
$$ a_{\mathcal{I}} \big( s(\xi)^{-1}s'(\xi) \big) \; \;  k_{\mathcal{I}} \big(s(\xi)^{-1}s'(\xi)\big) \in AM.
$$
\end{fait}
\begin{proof}
For every $\xi \in  \mathcal{F}_s \cap \mathcal{F}_{s'}$, we denote by $\mathscr{T}_{s,s'}(\xi):= a_{\mathcal{I}} \big( s(\xi)^{-1}s'(\xi) \big) \; \;  k_{\mathcal{I}} \big(s(\xi)^{-1}s'(\xi)\big).$ 
Due to the hypothesis that $s$ and $s'$ are both cross-sections of $G \rightarrow \mathcal{F}$, we deduce that $s'(\xi)\in~s(\xi) MAN$.
Hence the compact part 
$k_{\mathcal{I}} \big(s(\xi)^{-1}s'(\xi)\big)$ is in $M$ and $\mathscr{T}_{s,s'}(\xi)$ is in $AM$.
\end{proof}
Let us explicit the notation of the transfert maps between $\mathcal{B}_s$ and $\mathcal{B}_{s'}$ parametrizations where $s$ and $s'$ are differentiable cross-sections of intersecting domains.
\begin{defin}\label{defin-transfer}
Let $s$ and $s'$ be two differentiable cross-sections of $G \rightarrow \mathcal{F}$ such that $\mathcal{F}_s \cap \mathcal{F}_{s'} \neq \emptyset$.
We define the \emph{transition map} 
\begin{align*}
\mathscr{T}_{s,s'}: \mathcal{F}_s \cap \mathcal{F}_{s'}  & \longrightarrow AM \\
\xi & \longmapsto 
 a_{\mathcal{I}} \big( s(\xi)^{-1}s'(\xi) \big) \; \;  k_{\mathcal{I}} \big(s(\xi)^{-1}s'(\xi)\big), 
\end{align*}
which associate to every $\xi \in \mathcal{F}_s \cap \mathcal{F}_{s'}$, the unique element in $AM$ such that
$$s'(\xi) \in s(\xi)N\mathscr{T}_{s,s'}(\xi).$$
\end{defin}
Remark that the transition functions between two compact valued cross-sections take value in $M$, however, it is not the case when one of them is not compact valued.

\begin{prop}\label{prop_tranfer}
\hypG
Let $s$ and $s'$ be differentiable cross-sections of $G \rightarrow \mathcal{F}$ such that $\mathcal{F}_s \cap \mathcal{F}_{s'} \neq \emptyset$.

Then the following holds.
\begin{itemize}
\item[(i)] The map $\mathscr{T}_{s,s'}$ is differentiable and the identity map of $s(\mathcal{F}_s \cap \mathcal{F}_{s'}) N AM$ reads in $\mathcal{B}_{s'}$ and $\mathcal{B}_s$ coordinates as follows:
\begin{align*}
\big(\mathcal{F}_{s'}^{(2)} \cap \mathcal{F}_{s}^{(2)}\big) \times AM
&\longrightarrow 
\big(\mathcal{F}_{s}^{(2)} \cap \mathcal{F}_{s'}^{(2)}\big) \times AM
\\
 (\xi, \check{\xi} \; ; \; x)_{s'} &\longmapsto (\xi, \check{\xi} \; ; \; \mathscr{T}_{s,s'}(\xi) x)_{s}.
\end{align*}

\item[(ii)] For all differentiable cross-section $s''$ such that $\mathcal{F}_{s''}\cap \mathcal{F}_{s'} \cap \mathcal{F}_s \neq \emptyset,$ for all $\xi$ in the latter,
$$ \mathscr{T}_{s'',s}(\xi) = \mathscr{T}_{s'',s'}(\xi) \mathscr{T}_{s',s}(\xi).$$
\item[(iii)] For all $\xi \in \mathcal{F}_{s'} \cap \mathcal{F}_s$, 
$$ \mathscr{T}_{s',s}(\xi)= \mathscr{T}_{s,s'}(\xi)^{-1} .$$
\item[(iv)] For all $x\in AM$ and $u\in N$, 
$$\mathscr{T}_{s,s.xu}=x=\mathscr{T}_{s,s.ux}.$$
\end{itemize}
\end{prop}
The first three points enforce the computational 'rule' that double cross-sections subscript cancel.
\begin{proof}
(i) note that $s'( \mathcal{F}_s\cap \mathcal{F}_{s'}) NAM= s( \mathcal{F}_s\cap \mathcal{F}_{s'}) NAM$ since $s$ and $s'$ are both cross-sections of $G \rightarrow \mathcal{F}$.
We want to write every element of this subset of $G$ in the two coordinates.
By Definition \ref{defin_bh-coord} of the coordinate maps, the first two coordinates in $\mathcal{F}^{(2)}$ do not depend on $s$ and $s'$. 
We only need to compute the change in the last coordinate.
Fix an element $g \in  s'( \mathcal{F}_s\cap \mathcal{F}_{s'}) NAM$ and denote by $(\xi , \check{\xi} \; ; \; x)_{s'} \in \mathcal{F}_{s'}^{(2)} \times AM$ its coordinates with respect to the section $s'$. 
Using Fact \ref{fait_bh-coord} on $g$ and $s'$, there exists a unique element  $u_{\check{\xi}} \in N$ such that $g$ admits the following decomposition 
$$g=s'(\xi) u_{\check{\xi}}x.$$
Let us deduce the last $\mathcal{B}_s$ coordinate of $g$ by finding its decomposition in $s(\mathcal{F}_s\cap \mathcal{F}_{s'}) N AM$. 
Since $\xi \in \mathcal{F}_{s'} \cap \mathcal{F}_s$, by Definition \ref{defin-transfer}, there exists a unique element $ u_{s',s}(\xi) \in N$ such that
$$s'(\xi) = s(\xi) u_{s',s}(\xi) \mathscr{T}_{s,s'}(\xi).$$
Then we replace it in $s'(\xi) u_{\check{\xi}}x$, 
$$ s'(\xi) u_{\check{\xi}}x = s(\xi) u_{s',s}(\xi) \mathscr{T}_{s,s'}(\xi) u_{\check{\xi}}x .$$
Since $AM$ normalizes $N$ and $\mathscr{T}_{s,s'}(\xi) \in AM$, we deduce the following $s(\mathcal{F}_s\cap \mathcal{F}_{s'}) NAM$ decomposition of $g$,
$$g= s'(\xi) u_{\check{\xi}}x = s(\xi) \;\;\big( u_{s',s}(\xi) \mathscr{T}_{s,s'}(\xi) u_{\check{\xi}} \mathscr{T}_{s,s'}(\xi)^{-1}\big) \;\; \mathscr{T}_{s,s'}(\xi)x .$$
Hence, the $\mathcal{B}_{s}$-coordinates of $g$ is $\big(\xi , \check{\xi}\; ; \;  \mathscr{T}_{s,s'}(\xi)x\big)_s.$

(ii) is a direct consequence of the relation $\mathcal{B}_{s}\mathcal{B}_{s''}^{-1} =\mathcal{B}_{s}\mathcal{B}_{s'}^{-1} \; \mathcal{B}_{s'}\mathcal{B}_{s''}^{-1}$ where each map is restricted to $s\big(\mathcal{F}_{s''}\cap \mathcal{F}_{s'} \cap \mathcal{F}_s\big)N AM$.

(iii) follows from (ii) since $e_{AM}=\mathscr{T}_{s,s} = \mathscr{T}_{s,s'}\mathscr{T}_{s',s}$.

(iv) we recall that for all $x \in AM$ and $u \in N$ the section $s.xu$ (resp. $s.ux$) is defined for every $\xi \in \mathcal{F}_s$ by  $s.xu(\xi)=s(\xi) xu$ (resp. $s.ux(\xi)=s(\xi)ux$).
Using that $AM$ normalises $N$, we deduce the unique decomposition in $s(\mathcal{F}_s) N AM,$
$$s.xu(\xi)=s(\xi) \; ( xux^{-1} ) \; x .$$
Hence the maps $\mathscr{T}_{s,s.xu}$ and $\mathscr{T}_{s,s.ux}$ are constant equal to $x$. 
\end{proof}
\subsection{Cocycle}
Fix a covering family of differentiable cross-sections $(s_i)_{i\in I}$ of $G \rightarrow \mathcal{F}$ and let us read in $(\mathcal{B}_{s_i})_{i\in I}$ coordinates the left action of $G$ on itself.
The left action of $G$ on the first two coordinates in $\mathcal{F}^{(2)}$ is given by $g( \xi, \check{\xi})=(g\xi, g \check{\xi}).$ 

In Proposition \ref{prop_cocycle}, we prove that the $AM$-valued function defined below, called \emph{signed Iwasawa cocycle}, contains the remaining information on the third coordinate.
Its domain is in $G \times \mathcal{F}$, meaning that the information contained in the second and third coordinate in $\mathcal{F}^{(2)} \times AM$ are not needed when one reads the left action of $G$.

In Fact \ref{fait_cocycle-transfer}, we extend to the cocyle the Chasles type relations we previously had with the transition functions.

\begin{defin}\label{defin_cocycle}
Let $s_0, s_1$ be differentiable cross-sections of $G \rightarrow \mathcal{F}$.

For every $g \in G$ and $\xi \in \mathcal{F}_{s_0}$ such that $g \xi \in \mathcal{F}_{s_1}$, we denote $\beta_{s_1,s_0}(g,\xi)$ the unique element in $AM$ such that
$$ g s_0(\xi) \in s_1(g\xi)   \beta_{s_1,s_0}(g,\xi)N.$$
When $s_1=s_0$, we set $\beta_{s_0}:=\beta_{s_0,s_0}.$

Whenever $s_0$ and $s_1$ take value in $K$, the cocycle $\beta_{s_1,s_0}$ is called \emph{signed (multiplicative) Iwasawa cocycle} or in a shorter way, \emph{signed cocycle}.
\end{defin}

\begin{prop}\label{prop_cocycle}
\hypG
Let $s_0, s_1$ be differentiable cross-sections of $G \rightarrow \mathcal{F}$.

For all $g \in G$ and every element in $s_0(\mathcal{F}_{s_0}) NAM$ of coordinates $(\xi, \check{\xi} \; ; \; x)_{s_0} \in \mathcal{F}_{s_0}^{(2)} \times AM$ such that $g \xi \in \mathcal{F}_{s_1}$, we denote by $g \big(\xi,\check{\xi} \; ; \; x \big)_{s_0}$ its left multiplication by $g$.
Then the latter's coordinates with respect to $s_1$ satisfy
\begin{equation}\label{eq_g-cocycle}
g \big(\xi, \check{\xi} \; ; \; x \big)_{s_0}
= \big( g \xi, g \check{\xi} \; ; \; \beta_{s_1,s_0}(g,\xi)x \big)_{s_1}.
\end{equation}

Therefore, for every covering family of smooth cross-sections $(s_i)_{i\in I}$ of $G \rightarrow \mathcal{F}$, the $AM$-valued functions $(\beta_{s_i,s_j} )_{i,j\in I}$ satisfy the following cocycle relations. \\
For every $i,j,k \in I$, all $\xi_i \in \mathcal{F}_{s_i}$ and $g_j,g_k \in G$ such that $g_j \xi_i \in \mathcal{F}_{s_j}$ and $g_k g_j \xi_i \in \mathcal{F}_{s_k}$ then
\begin{equation}\label{eq_rel-cocycle}
\beta_{s_k,s_i}(g_kg_j, \xi_i)=\beta_{s_k, s_j}(g_k,g_j\xi_i)\; \beta_{s_j,s_i}(g_j,\xi_i) .
\end{equation}

Furthermore, for all $y\in AM,$ the right multiplication by $y$ of the element of coordinate $(\xi, \check{\xi} \; ; \; x)_{s_0}$ denoted by $( \xi, \check{\xi} \; ; \; x)_{s_0}y$ satisfies
\begin{equation}\label{eq_AM-action}
(\xi, \check{\xi} \; ; \; x)_{s_0} y = ( \xi, \check{\xi} \; ; \; xy)_{s_0}  .
\end{equation}
\end{prop}
As in Proposition \ref{prop_tranfer} concerning the transition functions, \eqref{eq_g-cocycle} and \eqref{eq_rel-cocycle} enforce the computational 'rule' that double cross-sections subscript cancel.
\begin{proof}
Because $AM$ normalizes $N$, the following diagram is $G$-equivariant for the left action of $G$ and commutative.
$$  \xymatrix{  && g \in G \ar[ld]_{N} \ar[rd]^{AM} \\ 
&gN \in G/N  \ar[rd]_{AM}   && G/AM \simeq \mathcal{F}^{(2)} \ni g(\eta_0, \check{\eta}_0) \ar[ld]^{N} \\
  && g \eta_0 \in \mathcal{F} }$$
Thanks to the lower left side $G/N \rightarrow \mathcal{F}$ of the diagram we deduce that local trivializations of $G \rightarrow \mathcal{F}$ induces local trivializations of $G/N \rightarrow \mathcal{F}$, of fiber $AM$.
Indeed, for every differentiable cross-section $s:\mathcal{F}_s \rightarrow G$, the map
$\mathcal{F}_s \times AM \rightarrow G/N $ that associate to $(\xi\; ; \; x )_s \in \mathcal{F}_s \times AM$ the element $s(\xi)x N \in G/N$ is the inverse of a local coordinate system.

\noindent Let $(s_i)_{i\in I}$ be a covering family of cross-sections of $G \rightarrow \mathcal{F}$.
Then the cocycles $(\beta_{s_i,s_j})_{i,j\in I}$ of Definition \ref{defin_cocycle} and the left action of $G$ on $\mathcal{F}$ encode the left action of $G$ on $G/N$.
Indeed, let $hN \in G/N$ be an element of coordinates $(\xi \; ; \; x)_{s_i} \in \mathcal{F}_{s_i} \times AM $ and $g \in G$ such that $g\xi \in \mathcal{F}_{s_j}$.
By the restricted coordinates map, we write 
$hN=s_i(\xi) x N.$
Hence $$ghN=gs_i(\xi)xN.$$
By Definition \ref{defin_cocycle} of $\beta_{s_j,s_i}$, there exists a unique $u\in N$ such that $gs_i(\xi)=s_j(\xi)\beta_{s_j,s_i}(g,\xi) u. $ 
Replacing it in the expression of $ghN$ and using that $AM$ normalizes $N$ yields
$$ghN= s_j(g\xi) \beta_{s_j,s_i}(g,\xi) u x N= 
s_j(g\xi) \beta_{s_j,s_i}(g,\xi)x \; \big( x^{-1}u x N\big).$$
Hence $ghN$ has coordinates $(g\xi \; ; \; \beta_{s_j,s_i}(g,\xi)x)_{s_j}.$

\noindent Thanks to the higher right hand side of the diagram, the same cocycles $(\beta_{s_i,s_j})_{i,j\in I}$ combined with the left action of $G$ on $\mathcal{F}^{(2)}$ allow us to write in local trivialisations the left action of $G$ on itself.
Hence, equation \eqref{eq_g-cocycle} holds.

The cocycle relation given by equation \eqref{eq_rel-cocycle} follows from the equivariance of the diagram for the left action of $G$.

For equation \eqref{eq_AM-action}, note first that for every cross-section $s$ of $G \rightarrow \mathcal{F}$,  the subset $s(\mathcal{F}_s)N AM$ is invariant by right $AM$-translation.
Furthermore, right translating by $AM$ preserve the parts of the decomposition in $s(\mathcal{F}_s)N$.
Finally, this translates in $\mathcal{F}_s^{(2)} \times AM$ to a trivial action in the $\mathcal{F}^{(2)}$ coordinates and a translation in the third $AM$ coordinate.
\end{proof}
Lastly, let us combine the relations between transition functions and cocycles for the $(\mathcal{B}_{s_i})_{i\in I}$ coordinates.
\begin{fait}\label{fait_cocycle-transfer}
\hypG
Let $s_0, s_0', s_1,s_1'$ be differentiable cross-sections of $G \rightarrow \mathcal{F}$.
Then for all $g \in G$ and $\xi \in \mathcal{F}_{s_0} \cap \mathcal{F}_{s_0'}$ such that $g \xi \in \mathcal{F}_{s_1} \cap \mathcal{F}_{s_1'}$,
$$\beta_{s_1',s_0'}(g,\xi)=
\mathscr{T}_{s_1',s_1}(g\xi) \; \beta_{s_1,s_0}(g,\xi)\; \mathscr{T}_{s_0,s_0'}(\xi)
.$$
\end{fait}
Note that cross-section subscript that are doubled, cancel out with our notations.
\begin{proof}
Let $(\xi,\check{\xi}\; ; \; x)_{s_0'} \in \big(\mathcal{F}_{s_0'} \cap \mathcal{F}_{s_0} \big) ^{(2)} \times AM$ and $g \in G$ such that $g \xi \in \mathcal{F}_{s_1'} \cap \mathcal{F}_{s_1}$.
By equation \eqref{eq_g-cocycle} of the previous Proposition \ref{prop_cocycle} for the local coordinates given by $s_1'$ and $s_0'$,   
$$ g(\xi, \check{\xi} \; ; \; x )_{s_0'}=(g\xi,g \check{\xi} \; ; \; \beta_{s_1',s_0'}(g,\xi) x )_{s_1'} .$$
Then by the transition identity of Proposition \ref{prop_bh-coord} (i) between $s_0'$ and $s_0$ on the left side of the previous equation,
$$g(\xi, \check{\xi} \; ; \; x )_{s_0'}=g(\xi, \check{\xi} \; ; \;\mathscr{T}_{s_0,s_0'}(\xi)  x )_{s_0}.$$
Again by the cocycle identity on the right hand side between $s_0$ and $s_1$,
$$g(\xi, \check{\xi} \; ; \; x )_{s_0'}
= (g\xi,g \check{\xi} \; ; \; \beta_{s_1,s_0}(g,\xi) \mathscr{T}_{s_0,s_0'}(\xi) x )_{s_1}.$$
Lastly, the transition identity between $s_1$ and $s_1'$ on the right side of the equation yields
$$(g\xi,g \check{\xi} \; ; \; \beta_{s_1',s_0'}(g,\xi) x )_{s_1'}= 
(g\xi,g \check{\xi} \; ; \; \mathscr{T}_{s_1',s_1}(g\xi) \beta_{s_1,s_0}(g,\xi) \mathscr{T}_{s_0,s_0'}(\xi)  x )_{s_1'}.$$
\end{proof}
\subsection{Local extensions of Hopf coordinates}
Given a family of covering differentiable cross-sections
$(s_i)_{i\in I}$ of $G \rightarrow \mathcal{F}$, the associated cocycles do not extend the Iwasawa cocycle.
Hence, in a general setting, the maps $(\mathcal{B}_{s_i})_{i\in I}$ do not extend Hopf coordinates of $G/M$.

We prove that when the cross-sections $(s_i)_{i \in I}$ take value in $K$, the signed multiplicative Iwasawa cocycles $(\beta_{s_i,s_j})_{i,j\in I}$ generalize the Iwasawa cocycle.
We obtain an equivariant and commutative diagram with the Hopf coordinates.

\begin{prop}\label{prop_bh-coord}
\hypG
Let $s$ be a \emph{compact} valued, differentiable cross-section of $G \rightarrow \mathcal{F}$, then 
$\mathcal{H}_s$ extends the Hopf coordinates restricted to $s(\mathcal{F}_s) NAM$ i.e. 
the following diagram is commutative.
$$  \xymatrix{  
	 s(\mathcal{F}_s)NAM \ar[d]_{\pi_M} 
	\ar[r]
		&\mathcal{F}_s^{(2)} \times AM \ar[d] \ni ( \xi, \check{\xi} \; ; \; x)_s \\ 
	 G/M \ar[r]  
		&\mathcal{F}^{(2)} \times \fa \ni
(\xi,  \check{\xi} \; ; \; \log x_A)
		 }$$
Moreover, it is equivariant with the left action of $G$, i.e. for all $\xi \in \mathcal{F}_s$, for all $g \in G$ and all compact valued section $s'$ such that $g \xi \in \mathcal{F}_{s'}$, the element
$$ g(\xi, \check{\xi} \; ; \; x)_{s}= (g \xi, g \check{\xi}\; ; \; \beta_{s',s}(g,\xi)x )_{s'}$$
projects in $G/M$ to 
$$ g(\xi,\check{\xi} \; ; \; \log x_A)= (g \xi, g \check{\xi} \; ; \; \sigma(g,\xi)+  \log x_A) .$$
Similarly, it is equivariant with the right action of $A$ i.e.
for all $(\xi, \check{\xi} \; ; \; x ) \in \mathcal{F}_s^{(2)} \times AM$ and all 
$\theta \in \fa\setminus \lbrace 0\rbrace$, 
the element
$$ \tilde{\phi}_\theta^t(\xi, \check{\xi}\; ; \; x)_s = (  \xi, \check{\xi}\; ; \; x e^{t \theta} )_s $$
projects to 
$$\tilde{\phi}_\theta^t(\xi, \check{\xi}\; ; \; \log x_A )=   ( \xi, \check{\xi} \; ; \; \log x_A+ t \theta ) .$$
\end{prop}

\begin{proof}
Recall that the lower part of the diagram reads as $gM \mapsto (g \eta_0, g \check{\eta}_0 \; ; \; \sigma(g,\eta_0))$.
The upper part reads as $g \mapsto (g \eta_0, g \check{\eta}_0 \; ; \; x )_s$ where $x$ is the component in $MA$ given by the Iwasawa decomposition of $s(g \eta_0)^{-1} g$. 
Commutativity of the diagram then follows from the hypothesis $s(\mathcal{F}_s)\subset K$ and the Definition \ref{defin_furstenberg} of the Iwasawa cocycle $\exp (\sigma(g,\eta_0)) = a_{\mathcal{I}}(g)= a_{\mathcal{I}}(s(g \eta_0)^{-1} g) $.

Let us check left $G$-equivariance.
Let $s$ and $s'$ be compact valued differential cross-sections.
By Proposition \ref{prop_cocycle} \eqref{eq_g-cocycle} the left $G$-action starting from $s(\mathcal{F}_s) NAM$ landing in $s'(\mathcal{F}_{s'})NAM$ is expressed in the $AM$-coordinate thanks to the cocycle $\beta_{s',s}$.
Since $s$ and $s'$ take value in $K$, by the previous point $\mathcal{H}_s$ and $\mathcal{H}_{s'}$ are local extensions of the Hopf coordinates where the left $G$-action is given in the $\fa$-coordinate by the Iwasawa cocycle.
Hence the equivariance.

The last part follows from the commutativity of the diagram and Proposition \ref{prop_cocycle} \eqref{eq_AM-action} that describe how to read the right multiplication by elements of $AM$ in coordinates associated to cross-sections of $G\rightarrow \mathcal{F}$. 
\end{proof}

\subsection{Local coordinates of $K$}
We prove that every differentiable compact valued cross-section of $G \rightarrow \mathcal{F}$ also induces local coordinates of $K$ that take value in $\mathcal{F} \times M$.
We explicit in coordinates the map $k_{\mathcal{I}}$.

By endowing $K$ with the left $G$-action defined for every $g \in G$ and $k \in K$ by $g.k=k_{\mathcal{I}}(gk)$, we make the projection $G$-equivariant.

\begin{prop}\label{prop_BH_K}
\hypG
Let $s$ be a differentiable \emph{compact} valued cross-section of $G \rightarrow \mathcal{F}$.
Then the restriction to the first and last coordinates of $\mathcal{H}_s$ provide local coordinates of $K$ as follows.
\begin{align*}
\mathcal{F}_s \times M & \longrightarrow s(\mathcal{F}_s)M \subset K \\
(\xi \; ; \; c)_s & \longmapsto s(\xi)c.
\end{align*}
Furthermore, the map $k_{\mathcal{I}}:G \rightarrow K$ reads in coordinates as
\begin{align*}
\mathcal{F}_s^{(2)} \times AM & \longrightarrow \mathcal{F}_s\times M  \\
(\xi, \check{\xi} \; ; \; x)_s & \longmapsto (\xi \; ;\; x_M)_s
\end{align*}
and for every covering family of compact valued cross-sections $(s_i)_{i\in I}$ of $G\rightarrow \mathcal{F}$, 
the $M$-coordinate of the cocycles $(\beta_{s_i,s_j})_{i,j\in I}$ parametrize the left $G$ action on $K$.
\end{prop}
\begin{proof}
The map $k \mapsto k \eta_0$ allows to identify $\mathcal{F}$ with $K/M$.
Consequently, every compact valued differential cross-section induces a local trivialization.

Let $s$ be a compact valued differential cross-section of $G \rightarrow \mathcal{F}$.
Then $s(\mathcal{F}_s)M \subset K$.
Coordinates in $\mathcal{F}_s^{(2)} \times AM$ are the same as unique decompositions in $s(\mathcal{F}_s)NAM$ where the $N$ part is associated to the second coordinate in $\mathcal{F}$ and the $AM$ part the last coordinate.
Since $AM$ normalises $N$ and $M$ commutes with $A$, the compact part of the Iwasawa decomposition $KAN$ of every element in  $s(\mathcal{F}_s)NAM$ is given by the product of its elements in $s(\mathcal{F}_s)$ and $M$. 
Hence, the following diagram
$$\begin{array}{ccccc}
G & \overset{AN}{\longrightarrow}& K & \overset{M}{\longrightarrow}& \mathcal{F} \\
g & \longmapsto& k_{\mathcal{I}}(g) & \longmapsto& g \eta_0
\end{array}$$
reads in local $\mathcal{B}_s$ coordinates as
$$\begin{array}{ccccc}
\mathcal{F}_{s}^{(2)} \times AM & \longrightarrow& \mathcal{F}_{s} \times M & \longrightarrow& \mathcal{F}_{s} \\
(\xi,\check{\xi} \; ; \; x)_{s} & \longmapsto& (\xi\; ; \; x_M)_s & \longmapsto& \xi.
\end{array}$$
Let $(s_i)_{i\in I}$ be a family of covering  differentiable cross-sections of $K \rightarrow \mathcal{F}$.
That the left $G$ action in $K$ reads as the projection in $M$ of the cocycles $\big(\beta_{s_i,s_j}\big)_{i,j \in I}$ now follows from the equivariance of the second diagram in local coordinates.  
\end{proof}
\subsection{Bruhat-Hopf coordinates} 
We define two (covering) families of cross-sections of $G \rightarrow \mathcal{F}$ defined on Bruhat cells of $\mathcal{F}$: \emph{unipotent} and \emph{compact Bruhat sections}.
We define \emph{Bruhat-Hopf coordinates} as the local extensions of Hopf coordinates given by Proposition \ref{prop_bh-coord} with respect to the compact Bruhat sections.
In Proposition \ref{prop_transfer-bruhat} we prove that every unipotent (resp. compact) Bruhat section is parametrized by a point of $\mathcal{F}$ and an element in $AM$ (resp. $M$).

By Corollary \ref{corol_bruhat} of Bruhat decomposition, the map
\begin{align*}
N^- & \longrightarrow N^- \eta_0 = \mathsf{b}(\check{\eta}_0) \\
u & \longmapsto u \eta_0
\end{align*}
is a diffeomorphism. 
Denote by $[e]$ its inverse.
\begin{defin}\label{defin-unip_bruhat-section}
A \emph{unipotent Bruhat section} is a left translate by $G$ of the map $[e]$.
We denote them by $[h]:= h \cdot [e]$ where $h \in G$.
For every $h \in G$, the unipotent Bruhat section $[h]$ has domain $hN^- \eta_0= \mathsf{b}(h\check{\eta}_0)$, codomain $hN^-$ and is defined for all $\xi \in  \mathsf{b}(h \check{\eta}_0)$ by
$$[h](\xi)= h [e](h^{-1}\xi). $$
 
A \emph{compact Bruhat section} is the compact component in the $KAN$ decomposition of a unipotent Bruhat section, meaning that for every $h \in G$, the associated compact Bruhat section is defined by $k_{\mathcal{I}}\circ [h]$.  

\emph{Bruhat-Hopf coordinates} (resp. \emph{Bruhat coordinates}) are the families of coordinates of $G$ given by covering families of compact (resp. unipotent) Bruhat sections.
\end{defin} 

For every $h \in G$, the unipotent Bruhat section $[h]$ and the compact Bruhat section $k_{\mathcal{I}}\circ [h] $ share the same domain: 
the Bruhat cell $\mathsf{b}(h\check{\eta}_0)$ opposite to $h \check{\eta}_0$.
However, the latter's codomain is in $K$.
Remark that Bruhat cells $\mathsf{b}(h\check{\eta}_0)$ are open submanifolds of $\mathcal{F}$ therefore of non-empty interior.
Using compactness of $\mathcal{F}$, one can choose finite families of covering Bruhat sections of any type.

For every $\check{\xi} \in \mathcal{F}$, we pick a compact element $h_{\check{\xi}} \in K$ such that $h_{\check{\xi}}\check{\eta}_0 = \check{\xi}$.
The choice of this compact family 
$( h_{\check{\xi}} )_{\check{\xi} \in \mathcal{F}} \subset K$ determines a covering family of unipotent Bruhat sections.
Abusing notation, we denote each of them by
$$[\check{\xi}]:= [h_{\check{\xi}}].$$
Likewise, we determine a choice of compact Bruhat section for every domain $\mathsf{b}(\check{\xi})$ where $\check{\xi} \in \mathcal{F}$.
We denote them by 
$$k(\check{\xi}):= k_{\mathcal{I}} \circ [\check{\xi}].$$
\begin{rem}\label{rem_unip-com-section}
The Proposition below implies that for any $h \in G$ such that $h\check{\eta}_0= \check{\xi}$, there is a unique element $x_* \in AM$ such that $[h]=[\check{\xi}].x_*$.

Similarly, any compact Bruhat section $s$ is determined by its domain $\mathsf{b}(\check{\xi})$ with $\check{\xi} \in \mathcal{F}$ and an element $c \in M$ such that $s= k(\check{\xi}).c$.
\end{rem}

Recall that $k_{\mathcal{I}-}$ (resp. $a_{\mathcal{I}-}$) denotes the coordinate in $K$ (resp. $A$) in the Iwasawa decomposition $G= KAN^-$ and that for every cross-sections $s,s'$ of $G \rightarrow \mathcal{F}$, for all $\xi \in \mathcal{F}_s \cap \mathcal{F}_{s'}$, we defined $\mathscr{T}_{s,s'}(\xi)$ as the unique element in $AM$ such that 
$$s'(\xi) \in s(\xi ) N \mathscr{T}_{s,s'}(\xi). $$
\begin{prop}\label{prop_transfer-bruhat}
\hypG
The following holds.
\begin{itemize}
\item[(1)] For every $u_* \in N^-$, then $[u_*]=[e]$ i.e.
$$\mathscr{T}_{[e],[u_*]}=e_{AM}.$$
\item[(2)] For every $x_* \in AM$ and $u_* \in N^-$, then 
$[x_* u_*]=[e].x_*= [u_* x_*]$ i.e.
$$\mathscr{T}_{[e],[x_* u_*]}=x_*=\mathscr{T}_{[e],[u_*x_*]}.$$
\item[(3)] For every $h \in G$, then 
$[h]=[k_{\mathcal{I}-}(h)]. a_{\mathcal{I}-}(h)  $ i.e. 
$$\mathscr{T}_{[k_{\mathcal{I}-}(h)], [h]}= a_{\mathcal{I}-}(h).$$
\end{itemize}
\end{prop}
\begin{proof}
Note that for every $h \in N^-AM$, because $h \check{\eta}_0=\check{\eta}_0$, every unipotent Bruhat section of this form have the same domain i.e.
$$ \mathcal{F}_{[h]}=\mathcal{F}_{[e]}= \mathsf{b}(\check{\eta}_0).$$

Let us prove (1) i.e. that for all $u_* \in N^-$ and $\xi \in \mathcal{F}_{[u_*]}$
then $[u_*](\xi) \in [e](\xi ) N.$
Since $[u_*](\xi)=u_* [e](u_*^{-1}\xi)$ for every $\xi \in \mathsf{b}(\check{\eta}_0)$, then by Definition \ref{defin-unip_bruhat-section} of $[e]$, we deduce that $[u_*]$ takes value in $N^-$.
Hence $[e](\xi)^{-1}[u_*](\xi) \in N^-$. 
Furthermore, using that $[e]$ and $[u_*]$ are cross-sections of $G \rightarrow \mathcal{F}$, we deduce 
$[e](\xi)^{-1}[u_*](\xi) \in N^- \cap MAN$.
Therefore, by uniqueness of the Bruhat decomposition $[e](\xi)^{-1}[u_*](\xi)=e_G$ and $\mathscr{T}_{[e],[u_*]}=e_{AM}$.

For statement (2), for all $(u_*,x_*) \in N^- \times AM$ and $\xi \in \mathcal{F}_{[e]}$, then $[x_*u_*](\xi) =x_*u_* [e](u_*^{-1}x_*^{-1}\xi)$. 
Using that $AM$ normalizes $N^-$, we deduce that the map $\xi \mapsto [x_*u_*](\xi)x_*^{-1}$ is a differentiable cross-section of $G \mapsto \mathcal{F}$ taking value in $N^-$ and of domain $\mathsf{b}(\check{\eta}_0)$. 
Hence, by uniqueness of the Bruhat decomposition in $N^- NAM$, we deduce that $\mathscr{T}_{[e],[x_*u_*].x_*^{-1}}=e_{AM}.$
Now we apply Proposition \ref{prop_tranfer} (ii) on transition functions to deduce that
$$e_{AM}=\mathscr{T}_{[e],[x_*u_*].x_*^{-1}}=\mathscr{T}_{[e],[x_*u_*]}
\mathscr{T}_{[x_*u_*],[x_*u_*].x_*^{-1}} .$$
Then point (iv) of the same Proposition yields
 $\mathscr{T}_{[x_*u_*],[x_*u_*].x_*^{-1}}=x_*^{-1}$, hence
$$\mathscr{T}_{[e],[x_*u_*]}x_*^{-1}=e_{AM}.$$
For the second part of the equality, note that $[u_*x_*]=[x_* (x_*^{-1}u_*x_*) ]$.
Since $AM$ normalizes $N^-$, the conjugated term is in $N^-$ and the rest follows from the previous point. 
 
For statement (3),  we write the $KAN^-$ decomposition $h=k_{\mathcal{I}-}(h) a_{\mathcal{I}-}(h) u_{\mathcal{I}-}(h) $.
Then by properties of the left action of $G$ on $[e]$, we deduce that
$$h \cdot [e]= k_{\mathcal{I}-}(h)\cdot [a_{\mathcal{I}-}(h) u_{\mathcal{I}-}(h)].$$
Hence by statement (2), we deduce $[h]=[k_{\mathcal{I}-}(h)]. a_{\mathcal{I}-}(h).$
\end{proof}

\section{Products of loxodromic elements}

Recall that an element of $G$ is \emph{unipotent} (resp. \emph{elliptic}, \emph{hyperbolic}) if it is conjugated to an element in $N$ (resp. $K$, $A$).  
By semisimplicity of the Lie group, every element $g\in G$ admits a unique decomposition $g=g_eg_hg_u$, called the \emph{Jordan decomposition},
	where $g_e$, $g_h$ and $g_u$ commute and $g_e$ (resp. $g_h$, $g_u$) is called the \emph{elliptic part} (resp. \emph{hyperbolic part}, \emph{unipotent part}) of $g$.  
	
\begin{defin}\label{defin_jordan-proj}
	For any element $g\in G$, there is a unique element $\lambda(g)\in \fa^+$ such that the hyperbolic part of $g$ is conjugated to $\exp(\lambda(g))\in A^+$.
	The map $\lambda : G \rightarrow \fa^+$ is called the \emph{Jordan projection}. 
	\end{defin}

	An element $g\in G$ is \emph{loxodromic}  if $\lambda(g) \in \fa^{++}$.
	Denote by $G^{lox}$ (resp. $\G^{lox}$) the subset of loxodromic elements of $G$ (resp. $\G$).
Since any element of $N$ that commutes with $A^{++}$ is trivial, the unipotent part of loxodromic elements is trivial.
Furthermore, the only elements of $K$ that commute with $A^{++}$ are in $M$. 
We deduce that the elliptic part of loxodromic elements is conjugated to elements in $M$.
Therefore, $g$ is loxodromic if and only if there exists $h \in G$ such that $h^{-1}gh \in MA^{++}$.

	Hence, for every loxodromic element $g\in G$, there exists $h_g\in G$ and $m(g)\in M$ so that we can write $g=h_g m(g)e^{\lambda(g)} h_g^{-1}$.
	However, for every $m\in M$ we can also write $$g=(h_g m) (m^{-1} m(g)m ) e^{\lambda(g)} (h_gm)^{-1}.$$	 
Which means that the \emph{angular} part $m(g)$ is only well defined up to conjugacy by $M$.
We thus use specific cross-sections of $G \rightarrow G/AM$, to study the elliptic part of products of loxodromic elements.

For every loxodromic element $g\in G$, denote by $g^+:=h \eta_0$ and $g^-:= h \eta_0$.
The Iwasawa cocycle of $g$ on $g^+$ is equal to its Jordan projection (see for instance \cite[Fact 2.6]{dang_glorieux_2020})
$$\sigma(g,g^+)=\lambda(g).$$
In §4.1, by using differential cross-sections of $G \rightarrow \mathcal{F}$ that factor the projection $G \rightarrow G/AM$, we extend locally and to loxodromic elements the previous formula.

In §4.2, we recall the dynamical properties of the left action of loxodromic elements on the Furstenberg boundary. 
This leads us to another definition of $(r,\varepsilon)$-loxodromic elements, where $r$ is a positive number that measures the distance between the attracting point of the loxodromic element and the boundary of its basin of attraction and $\varepsilon$ measures how contracting it is. 
Using the Bruhat sections of $G \rightarrow \mathcal{F}$, we give another proof that every loxodromic element, iterated enough times, will become $(r,\varepsilon)$-loxodromic.

In §4.3, we compute the cocycle of loxodromic elements on points of their basin of attraction, with respect to the latter's unipotent Bruhat section.
By defining the Ratio maps that are similar to the error terms for products of generic loxodromic elements given in \cite{benoist2000proprietes}, we obtain an exact computation in Proposition \ref{prop-magic-cocycle}.

We define in §4.4 a family of equicontinuity constants $\delta_{r,\varepsilon}$ for compact Bruhat sections.
We claim the construction can be adapted for any family of covering $K$-valued cross-sections of $G \rightarrow \mathcal{F}$ of the same type.

In the last paragraph, we estimate simultaneously the elliptic and hyperbolic part of products of generic loxodromic elements in Proposition \ref{prop_crucial}.
The proof is based on a Ping-Pong argument.

\subsection{Extended Jordan projections for loxodromic elements}
For every loxodromic element $g\in G$, denote by $g^+:=h \eta_0$ and $g^-:= h \eta_0$.
By \cite[Fact 2.6]{dang_glorieux_2020}, its Iwasawa cocycle on $g^+$ is equal to the Jordan projection
$$\sigma(g,g^+)=\lambda(g).$$
Let us define a multiplicative and local extension to $MA^{++}$ of the Jordan projection of loxodromic elements.
\begin{defin}\label{defin-mul-jordan}
Let $s$ be a differentiable cross-section of $G \rightarrow \mathcal{F}$.
For every loxodromic element $g \in G$ such that $g^+ \in \mathcal{F}_s$, we denote by
$$\mathscr{L}_s(g):= \beta_{s}(g,g^+). $$ 
For compact or unipotent Bruhat sections, such a map is called a \emph{local extended Jordan projection} (for loxodromic elements). 
\end{defin}

\begin{fait}\label{fait_lox}
\hypG
Fix a family of unipotent Bruhat sections denoted by $([\xi])_{\xi \in \mathcal{F}}$ of respective domains $\mathsf{b}(\xi)$.
Let $g\in G$ be a loxodromic element.
Then the following holds.
\begin{itemize}
\item[(1)] For every $h \in G$ such that $h^{-1}gh \in MA^{++}$ then $\mathscr{L}_{[h]}(g)=h^{-1}gh .$
\item[(2)] For every loxodromic element $g$, the element $h_g$ of Bruhat coordinates $(g^+,g^-\; ; \; e_{AM})_{[g^-]}$ satisfies 
$h_g^{-1}gh_g=\mathscr{L}_{[g^-]}(g).$
\item[(3)] For every cross-section $s$ such that $g^+ \in \mathcal{F}_s$, then
$$ \mathscr{L}_s(g)= \mathscr{T}_{[g^-],s}(g^+)^{-1}\mathscr{L}_{[g^-]}(g) \mathscr{T}_{[g^-],s}(g^+) .$$
\item[(4)] Therefore
$ \mathscr{L}_s(g) \in M e^{\lambda(g)} .$
\end{itemize}

\end{fait}
\begin{proof}
(1) By Definition \ref{defin_cocycle} of the cocycle, $\beta_{[h]}(g,g^+)$ is the unique element in $AM$ such that
$$g[h](g^+) \in [h](g^+)\beta_{[h]}(g,g^+) N .$$
By Definition \ref{defin-unip_bruhat-section} of the unipotent Bruhat section,
$[h](g^+)=h [e](h^{-1}g^+).$
Since $g^+=h \eta_0$, we deduce that $[h](g^+)=h [e](\eta_0)=h$.
We rewrite the inclusion, with the definition of the extended Jordan projection $\beta_{[h]}(g,g^+)=\mathscr{L}_{[h]}(g)$
$$g h \in h \mathscr{L}_{[h]}(g) N .$$
Since $h^{-1}gh \in MA^{++},$ we deduce that $\mathscr{L}_{[h]}(g) \in MA^{++}$ and the $N$-coordinate is trivial, i.e. 
$$gh=h \mathscr{L}_{[h]}(g).$$

(2) The unipotent Bruhat section $[g^-]$ shares the same domain as $[h]$.
By Remark \ref{rem_unip-com-section}, these cross-sections are defined only up to their domain and by right multiplication by an element in $AM$.
Since $h_g$ the unique element in $hMA$ of Bruhat coordinates $(g^+,g^- \; ; \; e_{AM} )_{[g^-]}$,
then 
$$g (g^+,g^- \; ; \; e_{AM} )_{[g^-]} = (g^+, g^- \; ; \; \beta_{[g^-]}(g,g^+)  )_{[g^-]}.$$
Using properties of the right translation by $AM$ in Bruhat coordinates, we deduce that
$$gh_g=h_g \mathscr{L}_{[g^-]}(g).$$

(3) Follows first from the identity between transition functions and cocycle of Fact \ref{fait_cocycle-transfer}, using that $g^+$ is a fixed point for the action of $g$ on $\mathcal{F}$, then we apply Proposition \ref{prop_tranfer} (iii) to recognize the inverse term for the transition function.

(4) Follows from (3) because we are conjugating by an element in $AM$.
\end{proof}

\subsection{Dynamical action on the Furstenberg boundary}
We study the left action of loxodromic elements on the Furstenberg boundary.
We give an alternative proof that the basin of attraction is the Bruhat cell opposite to the repelling point.
This leads to a Definition \ref{defin-r-eps} of $(r,\varepsilon)$-loxodromic elements using the $K$-invariant distance on $\mathcal{F}$.
We give another proof that large iterates of loxodromic element are $(r,\varepsilon)$-loxodromic.

\begin{prop}\label{prop-lox-bassin}
\hypG
Let $g\in G$ be a loxodromic element.

Then $g^+$ is an attracting point for the action of $g$ on the Furstenberg boundary.
Furthermore, the \emph{basin of attraction} of $g^+$ is $\mathsf{b}(g^-)$, the Bruhat cell opposite to its repelling point.
\end{prop}
The classical proof uses the fundamental representations of $G$ introduced Tits (cf. \cite[Corollary 3.12]{sambarino2014hyperconvex}) and involves the notion of simultaneaous proximality in those representations ( cf. \cite{benoist1997proprietes}).
We only rely here on Bruhat decomposition and the convergence \eqref{equ_unipotents}.
\begin{proof}
Let us first assume that $g \in MA^{++}$. 
Then $g^+=\eta_0$ and $g^-= \check{\eta}_0$ and we are going to prove that its basin of attraction is $\mathsf{b}(\check{\eta}_0)$.

Since $MA$ normalizes $N^-$, we deduce that $g$ stabilizes the Bruhat cell $\mathsf{b}(\check{\eta}_0)$.
Indeed, for every $u_* \in N^-$, then $g u_* \eta_0=gu_* g^{-1}\eta_0 \in N^- \eta_0$.
Furthermore for any $u_* \in N^-$, then $ g^n u_* g^{-n} \rightarrow e_G$ when $n \rightarrow + \infty$.
This implies that $\mathsf{b}(\check{\eta}_0)$ is in the bassin of attraction of $\eta_0$.

Conversely, let $\xi \in \mathcal{F}$ be in the basin of attraction of $\eta_0$.
Choose for every element in the Weyl group $w \in W$ a representative $k_w \in N_K(A)$ and recall that $k_\iota \in N_K(A)$ denotes an element such that $ N^-= k_\iota N k_\iota^{-1}$.
Apply Bruhat decomposition $G= k_\iota \sqcup_{w\in W} B k_w B $ where $B=MAN$.
Then there exists $u_* \in N^-$ and $k_w \in N_K(A)$ such that $\xi=u_* k_w \eta_0$.
Now $g^n \xi = g^n u_* g^{-n} \; (g^n k_w \eta_0).$
Since $\xi$ is in the basin of attraction and $g^n u_* g^{-n}$ converges to $e_G$, we deduce that $g^n k_w \eta_0 \rightarrow \eta_0$.
Using that $k_w$ normalizes $MA$, we deduce that 
$g^n k_w \eta_0= k_w (k_w^{-1}g^n k_w) \eta_0 = k_w \eta_0.$
The sequence is stationary at $k_w \eta_0$,
by uniqueness of the limit $k_w \eta_0 = \eta_0$.
Hence $k_w \in M$ and $\xi \in N^- \eta_0$.

In the general case, let $g \in G$ be a loxodromic element. 
Consider the unique element $h_g \in G$ given by Fact \ref{fait_lox} such that 
$$ h_g^{-1} gh_g = \mathscr{L}_{[g^-]}(g) \in MA^{++}.$$
Hence, the attracting point of $g$ is $h_g \eta_0=g^+$ and its basin of attraction is $h_g \mathsf{b}(\check{\eta}_0)=\mathsf{b}(g^-)$.
\end{proof}

We give a definition of $(r,\varepsilon)$-loxodromic elements which is slightly different from what the reader may find in \cite{benoist1997proprietes} or \cite{benoist2000proprietes} because it does not use the notion of simultaneous proximality in the fundamental representations of $G$ given by Tits.
However, using our choice of distance on $\mathcal{F}$ and the intrinsic characterization of the basin of attraction of loxodromic elements, one can check that both definitions are equivalent.

\begin{defin}\label{defin-r-eps}
Let $r>0$ be a positive number and $\varepsilon \in (0,r]$.
An element $g \in G$ is $(r,\varepsilon)$-loxodromic if it satisfies the following conditions.
\begin{itemize}
\item[(i)] The element $g$ is loxodromic and $r \leq  \frac{1}{2} \mathrm{d}(g^+, \partial \mathsf{b}(g^-))$.
\item[(ii)] It maps the compact set $\mathcal{V}_{\varepsilon}(\partial \mathsf{b}(g^-))^\complement$ into the ball $B(g^+, \varepsilon)$.
\item[(iii)] The restriction of $g$ to $\mathcal{V}_{\varepsilon}(\partial \mathsf{b}(g^-))^\complement$ is an $\varepsilon$-Lipschitz map.
\end{itemize}
\end{defin}
These remarks follow from the previous definition.
	\begin{itemize}\label{remarque_prox}
	\item[1)] If an element is $(r,\varepsilon)$-loxodromic, then it is $(r',\varepsilon)$-loxodromic for every $\varepsilon \leq r'\leq r$.
	\item[2)] If an element is $(r,\varepsilon)$-loxodromic, then it is $(r,\varepsilon')$-loxodromic for every $r \geq \epsilon'\geq  \varepsilon$.
	\item[3)] If $g$ is is $(r,\varepsilon)$-loxodromic, then $g^n$ is also is $(r,\varepsilon)$-loxodromic for every $n\geq 1$.
	\end{itemize}		
Note that loxodromic elements that are not sufficiently contracting, for instance those too close to $e_G$, will never satisfy the second condition for being $(r,\varepsilon)$-loxodromic.
However, we give below another proof that every loxodromic element, iterated a large enough amount of times will be $(r,\varepsilon)$-loxodromic.

\begin{prop}\label{prop-lox-lips}
\hypG
Let $g \in G$ be a loxodromic element.

Then for all positive number $r \leq \frac{1}{2} \mathrm{d}(g^+, \partial \mathsf{b}(g^-))$ and  all $\varepsilon \in (0,r]$, there exists an integer $N_{r,\varepsilon} \geq 1$ such that for all $n \geq N_{r,\varepsilon}$, the element $g^n$ is $(r,\varepsilon)$-loxodromic.   
\end{prop}
\begin{proof}
Let $g \in G$ be a loxodromic element and fix $r \leq \frac{1}{2} \mathrm{d}(g^+,\partial \mathsf{b}(g^-) ) $ and $\varepsilon \in (0,r]$.
By choice of these parameters, condition (i) holds.

Note that $\mathcal{V}_\varepsilon (\partial \mathsf{b}(g^-))^\complement$ is a compact subset of $\mathsf{b}(g^-)$, which by Proposition \ref{prop-lox-bassin} is the basin of attraction of $g^+$.
Hence $\lbrace g^n \mathcal{V}_\varepsilon (\partial \mathsf{b}(g^-))^\complement\rbrace_{n\geq 1}$ is a sequence of compact sets in the basin of attraction shrinking towards $g^+$.
Consequently, condition (ii) holds for every $n\geq N_2$ sufficiently large.

Let us now prove that there exists an integer $N_{r,\varepsilon} $ such that for every $n \geq N_{r,\varepsilon}$ the restriction of $g^n$ to  $\mathcal{V}_{\varepsilon}(\partial \mathsf{b}(g^-))^\complement$ is an $\varepsilon$-Lipschitz map.
By Fact \ref{fait_lox}, we consider the element $h_g \in G$ such that $h_g^{-1}g h_g = \mathscr{L}_{[g^-]}(g)$.
Using that $AM$ normalises $N^-$, we express the action of $g$ on $\mathsf{b}(g^-)$ in the unipotent charts $[g^-](\mathsf{b}(g^-))=h_g N^-$  by 
\begin{align*}
c(g): h_g N^- &\longrightarrow  h_g N^- \\
h_g u_* &\longmapsto h_g  \big( \mathscr{L}_{[g^-]}(g) u_* \mathscr{L}_{[g^-]}(g)^{-1}  \big) .
\end{align*}
The chosen metric on $\mathcal{F}$ is induced by the identification $K/M \simeq \mathcal{F}$.
Furthermore, the compact Bruhat section $k(g^-):\mathsf{b}(g^-)\rightarrow K$ defined by $k_{\mathcal{I}}\circ [g^-]$ 
is a differentiable chart of $\mathsf{b}(g^-)$.
Therefore, any upper bound of the differential of the map
\begin{align*}
k(g^-)( \mathsf{b}(g^-) )  &\longrightarrow k(g^-)( \mathsf{b}(g^-) ) \\
k_{\mathcal{I}}(h_g u_*) &\longmapsto k_{\mathcal{I}}\circ c(g) ( h_g u_*) .
\end{align*}
restricted to $k(g^-) (\mathcal{V}_\varepsilon(\partial \mathsf{b}(g^-))^\complement ) $  provides a Lipschitz constant for the map 
\begin{align*}
\mathcal{V}_\varepsilon(\partial \mathsf{b}(g^-)^\complement &\longrightarrow \mathsf{b}(g^-) \\
\xi & \longmapsto g \xi.
\end{align*}
Set $C_{r,\varepsilon}:=\sup_{u_* \in [g^-] B(g^+,\varepsilon) } 
\Vert D_{u_*}k_{\mathcal{I}} \Vert \; \sup_{u_* \in [g^-]\mathcal{V}_\varepsilon(\partial \mathsf{b}( g^-))^\complement } 
\Vert D_{u_*}k_{\mathcal{I}} \Vert^{-1} 
 .$\\
At every point, the eigenvalues of the differential of $c(g)$ are $\lbrace e^{-\alpha(\lambda(g))} \rbrace_{\alpha \in \Sigma_+}$ where $\Sigma_+$ is the set of positive roots.
Denote by $\ell_g := \min_{\alpha \in \Sigma^+} \alpha(\lambda(g))$.
Since $g$ is loxodromic, $\ell_g$ is a positive number and we obtain the uniform exponential decay of the differential of $c(g)$ i.e. for every $n\geq 1$, 
$$ \sup_{u_* \in h_gN^-} \Vert D_{u_*}c(g^n) \Vert \leq e^{-n\ell_g} . $$
By hypothesis on $r$ and $\varepsilon$, we deduce that $B(g^+,\varepsilon) \subset \mathcal{V}_{\varepsilon} (\partial \mathsf{b}(g^-))^\complement$. 
Let $n \geq N_2$. 
Then by choice of $\ell_g$ and $C_{r,\varepsilon}$, we deduce that $C_{r,\varepsilon} e^{-n \ell_g}$ is a Lipschitz constant for the action of $g^n$ restricted to this compact subset of the basin of attraction.
Since this sequence decays exponentially fast, there exists $N_{r,\varepsilon}\geq N_2$ such that for every $n\geq N_{r,\varepsilon}$, then $C_{r,\varepsilon} e^{-n \ell_g} \leq \varepsilon $ and condition (iii) is satisfied.
\end{proof}

\subsection{Cocycle on the basin of attraction}
Let $g$ be a loxodromic element.
We prove in Lemma \ref{lem_magic-cocyle-lox} that the cocycle of $g$ and every point $\xi$ in its basin of attraction, with respect to $[g^-]$ is everywhere equal to its local extended Jordan projection.
Then we define a map which allow us to write in Proposition \ref{prop-magic-cocycle} an exact relation between the cocycle of $g$ and every $\xi$ in its basin of attraction and any local extended Jordan projection, with respect to any suitable choice of cross-sections.

\begin{lem}\label{lem_magic-cocyle-lox}
\hypG 
Then for all loxodromic element $g \in G$, all $n\geq 1$ and $\xi \in \mathsf{b}(g^-)$,
$$\beta_{[g^-]}(g^n,\xi)=
\mathscr{L}_{[g^-]}(g)^n.$$
\end{lem}
We give a different proof from Lee-Oh \cite{lee2020ergodic}.
\begin{proof}
Denote by $h_g$ the element of $G$ of Bruhat coordinates $(g^+,g^-\; ; \; e_{AM})_{[g^-]}$, then $gh_g = h_g \mathscr{L}_{[g^-]}(g)$ by Fact \ref{fait_lox} and the unipotent Bruhat sections $[h_g]$ and $[g^-]$ are equal.

By property of the unipotent Bruhat section, for every $\xi \in \mathsf{b}(g^-)=h_g N^- \eta_0$, there exists a unique $u_\xi \in N^-$ such that 
$\xi = h_g u_\xi \eta_0$ and
$h_g u_\xi$ reads in Bruhat coordinates 
$(\xi,g^- \; ; \; e_{AM})_{[g^-]}.$

On one hand, by definition of the cocycle and because $\mathsf{b}(g^-)$ is the basin of attraction of $g^+$, for all $n\geq 1,$ the element $g^n h_g u_\xi$ reads as
$$g^n(\xi, g^- \; ; \; e_{AM})_{[g^-]} =(g^n\xi, g^- \; ; \;  \beta_{[g^-]}(g^n,\xi) \; )_{[g^-]}.$$
Note that $\beta_{[g^-]}(g^n,\xi)$ is the unique element in $AM$ such that $g^n h_g u_\xi \in h_g N^- \beta_{[g^-]}(g^n,\xi) N .$
That the second coordinate remains equal to $g^-$ means that the part in $N$ is trivial.

On the other hand, using the definition of the signed Jordan projection,
$$g^n h_g u_\xi = h_g \mathscr{L}_{[g^-]}(g)^n u_\xi= h_g \big(\mathscr{L}_{[g^-]}(g)^n u_\xi \mathscr{L}_{[g^-]}(g)^{-n} \big) \mathscr{L}_{[g^-]}(g)^n  .$$
Since $AM$ normalizes $N^-$, we deduce that  $\mathscr{L}_{[g^-]}(g)^n u_\xi \mathscr{L}_{[g^-]}(g)^{-n} \in N^-$, hence
 $$g^n h_g u_\xi \in  h_g N^- \mathscr{L}_{[g^-]}(g)^n .$$
This allows us to deduce by uniqueness of the Bruhat decomposition in $h_g N^- MA N$ that $\beta_{[g^-]}(g^n,\xi)= \mathscr{L}_{[g^-]}(g)^n. $
\end{proof}

\begin{defin}\label{defin-ratio}
Given two cross-sections $s_1,s_2$ and a Bruhat cell $\mathsf{b}(\check{\xi})$, then for all $\xi_1 \in \mathcal{F}_{s_1} \cap \mathsf{b}(\check{\xi})$ and $\xi_2 \in \mathcal{F}_{s_2} \cap \mathsf{b}(\check{\xi})$
we define the \emph{Ratio}
$$ \mathscr{R}_{s_1,s_2}(\check{\xi}; \xi_1,\xi_2):= 
\mathscr{T}_{s_1,[\check{\xi}]}(\xi_1) 
\mathscr{T}_{[\check{\xi}],s_2}(\xi_2).$$
When $s_1=s_2$, we shorten the notation $\mathscr{R}_{s_1}:=\mathscr{R}_{s_1,s_1}$.
For every loxodromic element $g \in G$ such that $g^+ \in \mathcal{F}_{s_1}\cap \mathsf{b}(g^-)$, for all $\xi \in \mathcal{F}_{s_2} \cap \mathsf{b}(g^-)$, set
$$\mathscr{R}_{s_1,s_2}(g,\xi):= \mathscr{R}_{s_1,s_2}(g^- ; g^+, \xi).$$
\end{defin}
The regularity of the Ratio map depends on the regulatity of the transfer maps which in turn depend on that of the cross-sections.
Because the transition functions between the unipotent and compact Bruhat sections take value in $AM$, for any compact Bruhat sections $s_1,s_2$, the \emph{Ratio map} $\mathscr{R}_{s_1,s_2}$ is continuous on its domain and takes value in $AM$.

Using the ratio map, the following statement follows from Lemma \ref{lem_magic-cocyle-lox}.
\begin{prop}\label{prop-magic-cocycle}
\hypG

Then for all loxodromic element $g \in G^{lox}$, all integer $n\geq 1$ and $\xi \in \mathsf{b}(g^-)$, for any choice of compact (Bruhat) sections $s_0,s_1,s_2$ such that $(\xi,g^+,g^n\xi)\in \mathcal{F}_{s_0} \times \mathcal{F}_{s_1} \times \mathcal{F}_{s_2}$
\begin{equation}\label{eq-magic-cocycle}
\beta_{s_2,s_0}(g^n,\xi)= \mathscr{R}_{s_1,s_2}(g;g^n\xi)^{-1} \mathscr{L}_{s_1}(g)^n \mathscr{R}_{s_1,s_0}(g;\xi).
\end{equation}
\end{prop} 
\begin{proof}
Indeed, by using first the transition functions between $s_2,s_0$ and $[g^-]$, then applying Lemma \ref{lem_magic-cocyle-lox} on the middle term and finally using the transition function between $[g^-]$ and $s_1$ in the middle term yields  
\begin{align*}
\beta_{s_2,s_0}(g^n,\xi) &= \mathscr{T}_{s_2,[g^-]}(g^n \xi) \quad \beta_{[g^-]}(g^n ,\xi) \quad \mathscr{T}_{[g^-],s_0}( \xi) \\
&= \mathscr{T}_{s_2,[g^-]}(g^n \xi) \quad \mathscr{L}_{[g^-]}(g^n) \quad \;\;\mathscr{T}_{[g^-],s_0}( \xi) \\
&= \mathscr{T}_{s_2,[g^-]}(g^n \xi) 
\mathscr{T}_{[g^-],s_1}(g^+) \;
 \mathscr{L}_{s_1}(g^n) \;
 \mathscr{T}_{s_1,[g^-]}( g^+)
 \mathscr{T}_{[g^-],s_0}( \xi).
\end{align*}
Finally, using Definition \ref{defin-ratio} and properties of the transition functions we check that \\
$ \mathscr{T}_{s_2,[g^-]}(g^n \xi) 
\mathscr{T}_{[g^-],s_1}(g^+)= \mathscr{R}_{s_1,s_2}(g;g^n\xi)^{-1}$ and
$\mathscr{T}_{s_1,[g^-]}( g^+)
 \mathscr{T}_{[g^-],s_0}( \xi)=\mathscr{R}_{s_1,s_0}(g;\xi).$
\end{proof}

\subsection{Equicontinuity constants}
After constructing a distance of $AM_0$ that is symmetric and left and right invariant, 
we introduce for every $r>0$ and $\check{\xi} \in \mathcal{F}$, a family $(\delta_{r,\varepsilon}(\check{\xi}))_{\varepsilon \in (0,r]}$ of equicontinuity constants of a continuous fonction defined over a compact set.
These constants are thus positive and converge to zero when $\varepsilon$ goes to zero.
Furthermore, using the $K$-invariance of the distance on $\mathcal{F}$ and the action of $K$ on the compact and unipotent Bruhat sections, we show that these constants do not depend on the choice of $\check{\xi} \in \mathcal{F}$.

Let us now choose a distance between on $AM$ and which is left and right $AM$ invariant i.e. for every $x,y,z,w \in AM$,
$\mathrm{d}_{AM}(wxz,wyz)=\mathrm{d}_{AM} (x,y)$.
The maximal torus $A$ is an abelian group of finite dimension, any norm on its Lie algebra $\fa$ will induce by the exponential map a suitable distance on the former.
Since $M$ is the centralizer of $A$ in $K$, it is sufficient to construct an $M$ invariant distance between arc connected points and setting an infinite distance otherwise.
We will then endow $AM$ with the distance $\mathrm{d}_{AM}$ induced by the product group structure $A \times M$.

Now we construct an $M$-invariant norm on the Lie algebra of $M$.
Starting from an euclidean norm on $\mathfrak{m}$, we make it $Ad(M)$-invariant by taking its average with respect to the Haar measure on $M$.
Since $M $ is compact, its Haar measure is finite.
Therefore, the average is an $Ad(M)$-invariant norm on $\mathfrak{m}$.
It induces an $Ad(M)$-invariant scalar product on $T_{e_M}M$.
By transporting it on the tangent space over every point by left multiplication by $M$ we obtain a left invariant metric.
The induced riemannian distance on $M$ is only defined between arc connected points and is, by construction, left $M$-invariant and invariant by conjugation.
This suffices to deduce the $M$-invariance of such a distance.

Recall that for every $\check{\xi} \in \mathcal{F}$, then $[\check{\xi}]$ denotes a choice of unipotent Bruhat section of domain $\mathsf{b}(\check{\xi})$ and $k(\check{\xi}):= k_{\mathcal{I}} \circ [\check{\xi}]$ is the associated compact Bruhat section.
Therefore $k (\mathcal{F})$ denotes the family of such compact Bruhat sections.

\begin{prop}\label{prop_equicontinu}
\hypG
Let $r>0$.
Consider the compact, symmetric and invariant by conjugacy by $K$ subset
$$K_r := \lbrace h \in K \; \vert \;   h \mathcal{V}_r (\partial \mathsf{b}(\check{\eta}_0)) \subset \mathcal{V}_{2r}(\partial \mathsf{b}(\check{\eta}_0) )  \rbrace.$$
For every $\check{\xi} \in \mathcal{F}$, denote by
$$ \delta_{r,\varepsilon}(\check{\xi}):= \sup_{s \in K_r \cdot k(\check{\xi})} \Big\lbrace \mathrm{d}_{AM}( \mathscr{R}_{s}(\check{\xi}\;;\; \xi_1,\xi_2),e_{AM} ) \; \Big\vert \; \xi_1 \in  \mathcal{V}_{3r}(\partial \mathsf{b}(\check{\xi}))^\complement \text{ and } \xi_2 \in B(\xi_1,\varepsilon)  \Big\rbrace .$$
Then the following holds.
\begin{itemize}
\item[(a)] For every $r>0$ and every $\varepsilon \in (0,r]$, the constant $\delta_{r,\varepsilon}(\check{\eta}_0)$ is non-zero and  
$$ \delta_{r,\varepsilon}(\check{\eta}_0)  \underset{\varepsilon \rightarrow 0}{\longrightarrow} 0 .$$
\item[(b)] For every $\check{\xi} \in \mathcal{F}$, the equality holds $\delta_{r,\varepsilon}(\check{\eta}_0)=\delta_{r,\varepsilon}(\check{\xi})$.
\end{itemize}
\end{prop}

\begin{proof}
Note first that by $M$-invariance of the distance
$$\mathrm{d}_{AM}( \mathscr{R}_{s}(\check{\xi}\;;\; \xi_1,\xi_2),e_{AM} ) = \mathrm{d}_{AM}( \mathscr{T}_{s,[\check{\xi}]}(\xi_1) \; , \mathscr{T}_{s,[\check{\xi}]}(\xi_2)).$$

(a) 
Because the Bruhat sections and the Iwasawa decomposition are differentiable,
by Definition \ref{defin-transfer} of the transition maps, the map 
\begin{align*}
K_r \times \mathcal{V}_{2r}(\partial\mathsf{b}(\check{\eta}_0))^\complement & \longrightarrow AM \\
(c,\xi) & \longmapsto \mathscr{T}_{c\cdot k(\check{\eta}_0), [\check{\eta}_0]}(\xi)
\end{align*}
is continuous.
It is defined over a compact set and using the Definition \ref{defin-ratio} of the ratio map, we notice that $\delta_{r,\varepsilon}(\check{\eta}_0)$ is for every $\varepsilon >0$ bounded above by the equicontinuity constant of this map and bounded below by the equicontinuity constants for the restriction to $K_r \times \mathcal{V}_{3r}(\partial\mathsf{b}(\check{\eta}_0))^\complement$.
Hence the positivity and convergence to zero.

(b) Let $l \in K$ such that $\check{\xi}=l \check{\eta}_0$ and $k(\check{\xi})=l \cdot k(\check{\eta}_0)$. 
Since $l \mathcal{V}_{2r}(\partial\mathsf{b}(\check{\eta}_0))^\complement= \mathcal{V}_{2r}(\partial\mathsf{b}(\check{\xi}))^\complement$, note that $\delta_{r,\varepsilon}(\check{\xi})$ is associated to the continuous map
\begin{align*}
K_r \times \mathcal{V}_{2r}(\partial\mathsf{b}(\check{\eta}_0))^\complement & \longrightarrow AM \\
(c,\xi) & \longmapsto \mathscr{T}_{cl\cdot k(\check{\eta}_0), l \cdot [\check{\eta}_0]}(l\xi).
\end{align*}
Recall that the transition function is the unique element in $AM$ such that
$$ \big( cl \cdot k(\check{\eta}_0) \big)(l\xi) \in \big( l \cdot [\check{\eta}_0] \big) (l \xi) N  \mathscr{T}_{cl\cdot k(\check{\eta}_0), l \cdot [\check{\eta}_0]}(l\xi).$$
By first applying the definition of the translation $l \cdot s (l \xi) =l s(\xi)$, then  multiplying by $l^{-1}$ on the left, we obtain
$$  \big( l^{-1}cl \cdot k(\check{\eta}_0) \big)(\xi) \in  [\check{\eta}_0]  ( \xi) N  \mathscr{T}_{cl\cdot k(\check{\eta}_0), l \cdot [\check{\eta}_0]}(l\xi).  $$
Therefore 
$\mathscr{T}_{cl\cdot k(\check{\eta}_0), l \cdot [\check{\eta}_0]}(l\xi)= \mathscr{T}_{l^{-1}cl\cdot k(\check{\eta}_0), [\check{\eta}_0]}(\xi).$
Since $K_r$ is invariant by conjugation, in particular $l^{-1}K_r l = K_r$.
The continuous maps associated to $\delta_{r,\varepsilon}(\check{\xi})$ and $ \delta_{r,\varepsilon}(\eta_0) $ coincide, hence the constants are equal.
\end{proof}

\begin{defin}\label{defin_d-r-eps}
Let $r>0$. 
We define the family of \emph{equicontinuity constants}
$$\delta_{r,\varepsilon}:= \sup_{\check{\xi} \in \mathcal{F}} \quad \sup_{s \in K_r \cdot k( \check{\xi})} \Big\lbrace \mathrm{d}_{AM}( \mathscr{R}_{s}(\check{\xi}\;;\; \xi_1,\xi_2),e_{AM} ) \; \Big\vert \; \xi_1 \in  \mathcal{V}_{3r}(\partial \mathsf{b}(\check{\xi}))^\complement \text{ and  } \xi_2 \in B(\xi_1,\varepsilon)  \Big\rbrace .$$
\end{defin}

\subsection{Estimates for products of generic loxodromic elements}

Let $g_1,...,g_l \in G^{lox}$ be loxodromic elements. 
Taking the convention that $g_0=g_l$, we say the (ordered) family is \emph{generic} if $g_{i-1}^+,g_i^-$ are transverse for every $1 \leq i \leq l$ or in other words $g_{i-1}^+ \in \mathsf{b}(g_i^-)$.

\begin{prop}\label{prop_crucial}
\hypG
For all $r>0$ and $\varepsilon \in (0,r]$ and every generic family $g_1,...,g_l \in G$ of $(r,\varepsilon)$-loxodromic elements such that 
\begin{itemize}
\item[$\star$] $ r \leq \frac{1}{6} \mathrm{d}( \lbrace g_{i-1}^+, g_i^+ \rbrace, \partial \mathsf{b} (g_i^-))$ for all $1 \leq i \leq l$
with the convention $g_0=g_l$.  
\end{itemize}
Fix a choice of compact Bruhat sections $(s_i)_{0 \leq i \leq l}$ such that 
\begin{itemize}
\item[$\star \star$] $\mathcal{F}_{s_i} \supset \mathcal{V}_r(\partial \mathsf{b}(g_i^-))^\complement$ 
for every $1 \leq i \leq l$ and $\mathcal{F}_{s_0}\supset \mathcal{V}_{\varepsilon}(\partial \mathsf{b}(g_1^-))^\complement$.
\end{itemize}
Then for all $\xi_0 \in \mathcal{V}_{\varepsilon}(\partial \mathsf{b}(g_1^-))^\complement$ and for all integers  $n_1,...,n_l \geq 1$,  
$$\beta_{s_l,s_0}(g_l^{n_l}...g_1^{n_1},\xi_0)\in 
\mathscr{L}_{s_l}(g_l^{n_l} )
\mathscr{R}_{s_l,s_{l-1}}(g_l, g_{l-1}^+) ... 
\mathscr{L}_{s_1}(g_1^{n_1})
\mathscr{R}_{s_1,s_0}(g_1, \xi_0)
 B( e_{AM}, (2l-1) \delta_{r,\varepsilon} ). $$
Furthermore, $g_l^{n_l}...g_1^{n_1}$ is $(2r,2 \varepsilon)$-loxodromic with attracting (resp. repelling) point in $B(g_l^+,\varepsilon)$  (resp. $B(g_1^-,\varepsilon)$) and its extended Jordan projection satisfies
$$ \mathscr{L}_{s_l}(g_l^{n_l}... g_1^{n_1}) \in 
\mathscr{L}_{s_l}(g_l^{n_l} )
\mathscr{R}_{s_l,s_{l-1}}(g_l, g_{l-1}^+) ...
\mathscr{L}_{s_1}(g_1^{n_1})
\mathscr{R}_{s_1,s_l}(g_1, g_l^+)
 B( e_{AM}, 2l \delta_{r,\varepsilon} ) .$$
\end{prop}
One can find a proof in \cite[Lemma 3.6]{benoist2000proprietes} that $g_l^{n_l}...g_1^{n_1}$ is $(2r,2 \varepsilon)$-loxodromic with attracting (resp. repelling) point in $B(g_l^+,\varepsilon)$  (resp. $B(g_1^-,\varepsilon)$) and obtain an estimate for the Jordan projection $\lambda$.
We improve the result by giving an estimate for the elliptic part. 
\begin{proof}
Let us prove the estimate for the cocycle, the extended Jordan projection's estimate will follow from the Definition \ref{defin-mul-jordan} that for every loxodromic element $g$ and a suitable cross-section $s$ such that $g^+ \in \mathcal{F}_s$, then $\mathscr{L}_{s}(g)=\beta_{s}(g,g^+)$.

For all $1 \leq j \leq l$, we set $\xi_j:=
g_{j}^{n_{j}}...g_1^{n_1} \xi_{0}$.
At each step starting from $j=1$, the element $g_j$ is $(r,\varepsilon)$-loxodromic and $\xi_{j-1} \in \mathcal{V}_{\varepsilon} (\partial \mathsf{b}(g_j^-) )^\complement \cap \mathcal{F}_{s_{j-1}}$.
Hence $\xi_j:= g_j^{n_j} \xi_{j-1} \in B(g_j^+,\varepsilon)$, which by choice of $s_j$ and $r \leq \frac{1}{6} \mathrm{d}(g_j^+, \partial \mathsf{b}(g_j^-) \cup \partial \mathsf{b}(g_{j+1}^-) )$ is inside $\mathcal{V}_{r}(\partial \mathsf{b}(g_{j+1}^-))^\complement \cap \mathcal{V}_{r}(\partial \mathsf{b} (g_j^-) )^\complement \subset \mathcal{F}_{s_j}$.
We deduce, by induction, that $\xi_{j} \in B(g_{j}^+, \varepsilon) \subset \mathcal{F}_{s_j}$ for every $1 \leq j \leq l$.

By the cocycle relation and recognizing $\xi_{j}$ for every $1 \leq j \leq l-1$,
\begin{align*}
\beta_{s_l,s_0}(g_l^{n_l}...g_1^{n_1},\xi_0) 
 &= \beta_{s_l,s_{l-1}}(g_l^{n_l},
g_{l-1}^{n_{l-1}}...g_1^{n_1} \xi_{0}) \; \beta_{s_{l-1},s_0}(g_{l-1}^{n_{l-1}}...g_1^{n_1},\xi_0) \\
&= \beta_{s_l,s_{l-1}}(g_l^{n_l},
\xi_{l-1}) \hspace*{1.3 cm} \beta_{s_{l-1},s_0}(g_{l-1}^{n_{l-1}}...g_1^{n_1},\xi_0) \\
&= \beta_{s_l,s_{l-1}}(g_l^{n_l},
\xi_{l-1}) \; \cdots \; \beta_{s_j,s_{j-1}}(g_j^{n_j}, \xi_{j-1}) \; \cdots \;  \beta_{s_1,s_0}(g_1^{n_1},\xi_0) . 
\end{align*} 
We will first prove that the first term on the right hand side is in a $2 \delta_{r,\varepsilon}$ neighbourhood of $\mathscr{L}_{s_l}(g_l^{n_l} )
\mathscr{R}_{s_l,s_{l-1}}(g_l \; ; \; g_{l-1}^+)$. 
It will then follow by induction that every term of the form $\beta_{s_j,s_{j-1}}(g_j^{n_j}, \xi_{j-1})$ where $2 \leq j \leq l$ is $2 \delta_{r,\varepsilon}$ close to $\mathscr{L}_{s_j}(g_j^{n_j} )
\mathscr{R}_{s_j,s_{j-1}}(g_j \; ; \;  g_{j-1}^+)$.
Finally, we prove that the last term is in a $\delta_{r,\varepsilon}$ neighbourhood of $\mathscr{L}_{s_1}(g_1^{n_1} )
\mathscr{R}_{s_1,s_{0}}(g_1 ; \xi_0).$

Let us apply Proposition \ref{prop-magic-cocycle}, then replace $g_l^{n_l}\xi_{l-1}$ with $\xi_l$
\begin{align*}
\beta_{s_l,s_{l-1}}(g_l^{n_l},
\xi_{l-1}) &= \mathscr{R}_{s_l} (g_l \; ;\; g_l^{n_l}\xi_{l-1} )^{-1} \; \mathscr{L}_{s_l}(g_l)^{n_l} \; \mathscr{R}_{s_l,s_{l-1}} (g_l \; ; \; \xi_{l-1}) \\
&= \mathscr{R}_{s_l} (g_l \; ;\; \xi_l )^{-1} \; \mathscr{L}_{s_l}(g_l)^{n_l} \; \mathscr{R}_{s_l,s_{l-1}} (g_l \; ; \; \xi_{l-1}).
\end{align*}
By Definition \ref{defin-ratio} of the ratio $\mathscr{R}_{s_l} (g_l \; ;\; \xi_l )^{-1} = \mathscr{R}_{s_l} (g_l^- \; ;\; g_l^+ , \xi_l )^{-1}.$
Since $\xi_l \in B(g_l^+, \varepsilon)$ and by choice of $r \leq \frac{1}{6} \mathrm{d}(g_l^+,  \partial\mathsf{b}(g_l^-))$, we deduce by Definition \ref{defin_d-r-eps} of $\delta_{r\varepsilon}$ that 
$ \mathscr{R}_{s_l} (g_l \; ;\; \xi_l ) \in B(e_{AM},\delta_{r,\varepsilon}) .$
The first term is small, it remains to show that the third term is close to $\mathscr{R}_{s_l,s_{l-1}}(g_l^-; g_l^+, g_{l-1}^+).$
By definition of the ratio map,
\begin{align*}
\mathscr{R}_{s_l,s_{l-1}}(g_l^-; g_l^+, \xi_{l-1})
&= 
\mathscr{T}_{s_l,[h_l]}(g_l^+) \hspace*{2.4cm}   \mathscr{T}_{[h_l],s_{l-1}}(\xi_{l-1}) \\
&=
\mathscr{T}_{s_l,[h_l]}(g_l^+) \; 
\mathscr{T}_{[h_l],s_{l-1}}(g_{l-1}^+)
 \\
&\hspace*{2.3 cm}
\mathscr{T}_{s_{l-1},[h_l]}(g_{l-1}^+)
 \mathscr{T}_{[h_l],s_{l-1}}(\xi_{l-1}) \\
&= \mathscr{R}_{s_l,s_{l-1}}(g_l^-; g_l^+, g_{l-1}^+) 
\hspace*{.8 cm}
\mathscr{R}_{s_{l-1}}(g_l^-; g_{l-1}^+, \xi_{l-1}).
\end{align*}
Hence, the cocycle can be written as follows,
$$\beta_{s_l,s_{l-1}}(g_l^{n_l},\xi_{l-1}) 
= 
\mathscr{R}_{s_l}(g_l \; ; \; \xi_l)^{-1} \quad
\mathscr{L}_{s_l}(g_l)^{n_l} \;
\mathscr{R}_{s_l,s_{l-1}}(g_l \; ; \; g_{l-1}^+) \quad 
\mathscr{R}_{s_{l-1}}(g_l^-; g_{l-1}^+, \xi_{l-1}). 
$$
Finally, by choice of $r \leq \frac{1}{6} \mathrm{d}(g_{l-1}^+, \partial \mathsf{b}(g_l^-))$ and Definition \ref{defin_d-r-eps} of $\delta_{r,\varepsilon}$, the third term is small i.e.
$\mathscr{R}_{s_{l-1}}(g_l^-; g_{l-1}^+, \xi_{l-1}) \in  B(e_{AM}, \delta_{r,\varepsilon}).$
Given that the distance in $AM$ is symmetric and invariant by conjugation, we deduce that
$$ \beta_{s_l,s_{l-1}}(g_l^{n_l},\xi_{l-1}) 
\in 
\mathscr{L}_{s_l}(g_l)^{n_l} \;
\mathscr{R}_{s_l,s_{l-1}}(g_l \; ; \; g_{l-1}^+) \; B(e_{AM}, 2 \delta_{r,\varepsilon}) 
. $$
By induction, for every $2 \leq j \leq l$
$$ \beta_{s_j,s_{j-1}}(g_j^{n_j},\xi_{j-1}) 
\in 
\mathscr{L}_{s_j}(g_j)^{n_j} \;
\mathscr{R}_{s_j,s_{j-1}}(g_j \; ; \; g_{j-1}^+) \; B(e_{AM}, 2 \delta_{r,\varepsilon}) 
. $$
Now for $\beta_{s_1,s_0}(g_1^{n_1},\xi_0)$, by Proposition \ref{prop-magic-cocycle} and by replacing $g_1^{n_1}\xi_{0}$ with $\xi_1$
$$\beta_{s_1,s_{0}}(g_1^{n_1},
\xi_{0}) = \mathscr{R}_{s_1} (g_1 \; ;\; g_1^{n_1}\xi_{0} )^{-1} \; \mathscr{L}_{s_1}(g_1)^{n_1} \; \mathscr{R}_{s_1,s_{1}} (g_1 \; ; \; \xi_{0}).$$
Similarly, by choice of $r$ and definition of $\delta_{r,\varepsilon}$, we deduce that 
$$\beta_{s_1,s_{0}}(g_1^{n_1},
\xi_{0}) \in  \mathscr{L}_{s_1}(g_1)^{n_1} \; \mathscr{R}_{s_1,s_{1}} (g_1 \; ; \; \xi_{0}) \; B(e_{AM},\delta_{r,\varepsilon}).$$
Hence, 
$$\beta_{s_l,s_0}(g_l^{n_l}...g_1^{n_1},\xi_0)\in 
\mathscr{L}_{s_l}(g_l^{n_l} )
\mathscr{R}_{s_l,s_{l-1}}(g_l, g_{l-1}^+) ... 
\mathscr{L}_{s_1}(g_1^{n_1})
\mathscr{R}_{s_1,s_0}(g_1, \xi_0)
 B( e_{AM}, (2l-1) \delta_{r,\varepsilon} ). $$

Finally, for the extended Jordan projection, we apply the cocycle estimate for the attracting point $g^+$ of $g_l^{n_l}...g_1^{n_1}$ with cross-section $s_0=s_l$.
By \cite[Lemma 3.6]{benoist2000proprietes}, it is in $B(g_l^{n_l},\varepsilon)$, therefore by choice of $r\leq \frac{1}{6} \mathrm{d}( g_l^+,\partial \mathsf{b}(g_1^-))$, we deduce that
$ \mathscr{R}_{s_1,s_l}(g_1, g^+) \in \mathscr{R}_{s_1,s_l}(g_1, g_l^+) B(e_{AM},\delta_{r,\varepsilon}) .$

Hence $ \mathscr{L}_{s_l}(g_l^{n_l}... g_1^{n_1}) \in 
\mathscr{L}_{s_l}(g_l^{n_l} )
\mathscr{R}_{s_l,s_{l-1}}(g_l, g_{l-1}^+) ...
\mathscr{L}_{s_1}(g_1^{n_1})
\mathscr{R}_{s_1,s_l}(g_1, g_l^+)
 B( e_{AM}, 2l \delta_{r,\varepsilon} ) .$
\end{proof}

\begin{defin}\label{defin_schottky}
	Let $0<\varepsilon\leq r$. 
	A semigroup $\G \subset G$ is \emph{strongly $(r,\varepsilon)$-Schottky} if
	\begin{itemize}
	\item[(i)] every element is $(r,\varepsilon)$-loxodromic, 
	\item[(ii)] $d(h^+,\partial \mathsf{b}(h'^{-})  )\geq 6r$ for all $h,h'\in \G$. 
	\end{itemize} 
	We also write that $\G$ is a \emph{strong $(r,\varepsilon)$-Schottky semigroup}.
\end{defin}

\section{Invariant sets}
Let $\G$ is a Zariski dense subgroup of $G$. 
In the first paragraph, following \cite{dang_glorieux_2020}, we construct the non-wandering set $\Omega \subset \G \backslash G/M$ for regular Weyl chamber flows.
We notice that it is the smallest closed $A$-invariant subset of $\G \backslash G/M$ containing all the periodic orbits of the flows $\phi_{\lambda(\G^{lox})}^t$.

Denote by $\widetilde{\Omega}_G$ the preimage of $\Omega$ via the projection $G \rightarrow \G \backslash G /M$. 
Such a subset is closed, left $\G$-invariant and right $AM$-invariant.
Denote by $M_0$ the connected component of the identity of $M$.
In the second paragraph, following Guivarc'h-Raugi \cite{guivarc2007actions}, we introduce the sign group $M_\G$, a normal subgroup of finite index of $M$ containing $M_0$.
One can find another construction of the sign group in \cite{benoist2005convexes3}. 

Finally, using Guivarc'h-Raugi's classification of $\G$-invariant subsets of $K$ (cf. Theorem \ref{theo_GR_Gamma_orbite_K}) we construct a partition of left $\G$-invariant right $AM_\G$-invariant subsets of $\widetilde{\Omega}_G$.
We prove in Proposition \ref{prop_transitivite-G} that the topological dynamics of diagonal flows on these subsets are all conjugated.

\subsection{In the space of Weyl chamber}
\begin{defin}\label{defin_limit_set}
A point $\eta \in \mathcal{F}$ is a \emph{limit point} if there exists a sequence $(\g_n)_{n\geq 1}$ in $\G$ such that $\big( (\g_n)_* \mathrm{Haar}_{\mathcal{F}} \big)_{n\geq 1}$ converges weakly towards the Dirac measure in $\eta$.

The \emph{limit set} of $\G$, denoted by $L_+(\G)$, is the set of limit points of $\G$. 
It is a closed, $\G$-invariant subset of $\mathcal{F}$. 

Denote by $L_-(\G)$ the limit set of $\G^{-1}$ and finally let $L^{(2)}(\G)= \Big( L_+(\G)\times L_-(\G) \Big) \cap \mathcal{F}^{(2)} $.

\end{defin}
Note that when $\G$ is a subgroup, then $L_+(\G)=L_-(\G)$ and $L^{(2)}(\G)$ is the subset of pair of points of $L_+(\G)$ in general position. 
For the hyperbolic plane, we get the product of the usual limit set minus the diagonal.

By \cite{benoist1997proprietes} Lemma 3.6, the set of pairs of attracting and repelling points of loxodromic elements of $\G$ is dense in $L_+(\G)\times L_-(\G)$.
Therefore, using Hopf coordinates and the construction of attracting and repelling points of loxodromic elements, $L^{(2)}(\G)$ identifies with smallest closed $\G$-invariant subset of $G/AM$ containing
$$ \Big\lbrace h_\gamma AM \; \vert \; \exists \gamma \in \G^{lox} \text{ such that } h_\gamma^{-1} \gamma h_\gamma \in M A^{++}   \Big\rbrace .$$

\begin{theorem}[Theorem 4.5 \cite{dang_glorieux_2020}]\label{transitivity of Furstenberg boundary}
Let $G$ be a real linear, connected, semisimple Lie group of non-compact type (i.e. without compact factors) and $\Gamma$ be a Zariski dense subsemigroup of $G$.
 
Then the (diagonal) action of $\Gamma$ on $L^{(2)}(\Gamma)$ is topologically transitive, i.e. there are dense $\G$-orbits.
\end{theorem}

\begin{defin}
We denote by $\widetilde{\Omega}$ the subset of \emph{non-wandering Weyl chambers}, defined through the Hopf parametrization by: 
$$\widetilde{\Omega}:= \cH^{-1}(L^{(2)}(\G)\times \fa).$$
This is a $\G$-invariant and right $A$-invariant subset of $G/M$.
When $\G$ is a subgroup, we denote by $\Omega:=\G \backslash \widetilde{\Omega}$ the quotient space.
\end{defin}
By Theorem \ref{transitivity of Furstenberg boundary} the quotient $\Omega$ is the smallest closed $A$-invariant subset of $\G \backslash G/M$ containing the following subset
$$ \Big\lbrace \phi_{\lambda(\gamma)}^\mathbb{R} (h_\gamma M) \; \vert \; \gamma \in \G^{lox} \text{ and } h_\gamma^{-1} \gamma h_\gamma \in M A^{++}   \Big\rbrace .$$
Note that in rank one, the above set is the reunion of all periodic orbits for the geodesic flow. 
\subsection{The sign subgroup}
Denote by $M^{ab}:=M/[M,M]$ the abelianisation of the compact group $M$ and by $\pi_{ab}: M \rightarrow M^{ab}$ the projection.
Abusing notations, $\pi_{ab}$ also denotes the projection $AM \rightarrow AM^{ab}.$

\begin{fait}
\hypG
For all cross-sections $s$ and $s'$ the projection into $AM^{ab}$ of the signed translation maps $\mathscr{L}_s$ and $ \mathscr{L}_{s'}$ coincide on the intersection of their domains i.e. for every loxodromic element $g\in G$ such that $g^+ \in \mathcal{F}_s \cap \mathcal{F}_{s'}$, then
$$\pi_{ab} \big( \mathscr{L}_s(g) \big) = \pi_{ab} \big( \mathscr{L}_{s'}(g) \big) .$$
\end{fait}
\begin{proof}
By Definition \ref{defin-mul-jordan} of the signed translation map, we write 
$$ \mathscr{L}_{s'}(g)= \beta_{s'}(g,g^+). $$
Using first the Fact \ref{fait_cocycle-transfer} over signed cocycles and transition functions, then that $g^+$ is fixed by $g$, we deduce that
$$ \mathscr{L}_{s'}(g)= \mathscr{T}_{s',s}(g^+) \beta_{s}(g,g^+) \mathscr{T}_{s,s'}(g^+). $$
The middle term is an extended Jordan projection and  the first and last term are inverse (see Proposition \ref{prop_tranfer} (iii) on transition functions).
Hence 
$$ \mathscr{L}_{s'}(g)= \mathscr{T}_{s',s}(g^+) \mathscr{L}_{s}(g) \mathscr{T}_{s',s}(g^+)^{-1}. $$
The claim then follows by projecting the relation into the abelian group $AM^{ab}$. 
\end{proof}
Denote by $\mathscr{L}^{ab}$ the map that associates to every loxodromic element $g \in G$   the projection into $AM^{ab}$ of any signed Jordan projection.
We call this map the \emph{abelian signed Jordan projection}.
Denote by $\pi_{M^{ab}}$ the projection $AM^{ab} \rightarrow M^{ab}$.

\begin{defin}\label{defin-sign-group}
 Let $\G$ be a Zariski dense subgroup of $G$.
 Denote by $\G^{lox}$ the subset of loxodromic elements of $\G$.
We define the \emph{abelian sign group} of $\G$ by
$$M_\Gamma^{ab}:= \pi_{M^{ab}}\Big( \overline{ \langle \mathscr{L}^{ab}(\Gamma^{lox}) \rangle} \Big).$$
The \emph{sign group} of $\G$, denoted by $M_\G$ is given by
$M_\Gamma := \pi_{ab}^{-1} (M_\Gamma^{ab}).$
\end{defin}
The following Theorem will imply non-arithmeticity in $AM_\Gamma^{ab}$ of the abelian signed Jordan projections of $\G$. 
\begin{theorem}[Theorem 6.4 \cite{guivarc2007actions}]\label{thm_GR_non-arithm}
Let $G$ be a real linear, connected, semisimple Lie group of non-compact type (i.e. without compact factors).
Then for all Zariski dense subsemigroup $\Gamma \subset G$, the closed subgroup spanned by $\mathscr{L}^{ab}(\G^{lox})$ is of finite index in $AM^{ab}.$
\end{theorem}
\begin{corol}\label{corol_GR_non-arithm}
Let $G$ be a real linear, connected, semisimple Lie group of non-compact type (i.e. without compact factors) and $\G$ be a Zariski dense subsemigroup of $G$.
Then
$$ \overline{ \langle \mathscr{L}^{ab}(\Gamma^{lox}) \rangle} =AM_\Gamma^{ab}.$$
\end{corol}
\begin{proof}
Denote by 
$H:= \overline{ \langle \mathscr{L}^{ab}(\Gamma^{lox}) \rangle}. $
By definition, it is a closed subgroup of $AM^{ab}$.
In particular, $AM^{ab}/H$ is Hausdorff.
According to the previous Theorem \ref{thm_GR_non-arithm}, it is a finite group.
By endowing it with the discrete topology, we deduce that the morphism
\begin{align*}
\varphi : A & \longrightarrow AM^{ab}/H \\
a & \longmapsto aH
\end{align*}
is a continuous map that takes value in a finite group.
Since $A$ is connected, $\varphi$ is constant to $e_A H$, hence $\varphi(A)=AH = H$ and $A \subset H$. 

By Definition \ref{defin-sign-group} of the sign group, $H \subset AM_\G^{ab}$.
Conversely, for every $x \in A$ and $m\in M_\G^{ab}$, there exists $y \in A$ such that $ym \in H$.
We now write $xm$ as a product $xm=(xy^{-1}) ym.$
In the right hand side, the first term is in $A$ hence in $H$ and the second term is in $H$, hence $xm \in H$.
We thus conclude that $H= AM_\G^{ab}.$
\end{proof}

\begin{theorem}[Theorem 8.2 \cite{benoist2005convexes3},Theorem 1.9 \cite{guivarc2007actions}]\label{thm_GR_M_G}
Let $G$ be a real linear, connected, semisimple Lie group of non-compact type (i.e. without compact factors) and $\Gamma$ be a Zariski dense subsemigroup of $G$.
Then the following holds.
\begin{itemize}
\item[(a)] $M_\G$ is a closed normal subgroup of finite index of $M$ and contains the connected component of the identity $M_0$.
\item[(b)] There exists an integer $p_\Gamma \in [0,\dim \fa]$ such that $M_\G/M_0$ is isomorphic to $\big(\mathbb{Z}/2\mathbb{Z}\big)^{p_\Gamma}$. 
\item[(c)] $M_{\G^{-1}}= k_\iota M_\G k_\iota^{-1}$ where $k_\iota\in N_K(A)$ is an element such that $Ad(k_\iota)\fa^+=-\fa^+$.
\item[(d)] For all $g\in G$, the groups satisfy $M_{g\G g^{-1}}=M_\G$.
\end{itemize}
\end{theorem}

When $G$ is a split, real linear, algebraic group, Y. Benoist in \cite{benoist97positives} studies the following conditions:  
\begin{itemize}
\item[(C1)] There exists a Zariski dense subgroup $\G \subset G$ such that $M_\Gamma=M_0$.
\item[(C2)] There exists a Zariski dense subgroup $\G \subset G$ with $M_\G \supsetneq M_0$ such that the sign group of every Zariski dense subgroup of $\G$ strictly contains $M_0$.
\end{itemize}
In particular, he proves for $\mathrm{SL}(m, \R)$ that both conditions hold when $m$ is a multiple of $4$, in fact (C2) is true for all $m$.
.
However, when $m$ is even but not divisible by $4$, condition (C1) is false i.e. the sign group of every Zariski dense subgroup of $\mathrm{SL}(m,\R)$ is non trivial.

\subsection{$\G$-invariant subsets of $G$}
Recall the left action of $G$ on $K$ defined
for all $g \in G$ and $k \in K$ by $g.k = k_{\mathcal{I}}(gk).$
The projection 
$K \rightarrow \mathcal{F}$ is thus equivariant and endows $K$ with a fiber bundle structure of fiber $M$ over the Furstenberg boundary.
Apply a result of Guivarc'h-Raugi \cite{guivarc2007actions} on the left action of $G$ on $K$.
Denote by $L_G(\G)$ the preimage in $K$ of the limit set $L(\G)\subset \mathcal{F}$.
Then the closed right $M$-invariant and left $\G$-invariant subset $L_G(\G) \subset K$ partitions 
into $\vert M/M_\G \vert$ closed, $\G$-invariant, minimal subsets.
Furthermore, these invariant subsets are right $M_\G$-invariant.
Lastly, using Iwasawa decomposition, we partition $\widetilde{\Omega}_G$ into left $\G$-invariant and right $AM_\G$-invariant subsets of $G$.

\begin{theorem}[Theorem 2 \cite{guivarc2007actions} ]\label{theo_GR_Gamma_orbite_K}
\hypG
Let $\G$ be a Zariski dense, discrete subsemigroup of $G$.

Then the following holds.
\begin{itemize}
\item[1)] $L_G(\G)\subset K$ partitions into $\vert M/M_\G\vert$ closed, minimal $\G$-invariant subsets i.e. in each partition, every $\G$-orbit is dense.
\item[2)] There is an indexation of this partition by $M/M_\G$ i.e. $L_G (\G) = \sqcup_{ [m] \in M/M_\G} L_{[m]}(\G)$ such that for every $m\in M$,
$$ L_{[m]}(\G)=L_{[e_M]}(\G)m.$$
\item[3)] Every element of the partition turns out to be right $M_\G$-invariant.
\end{itemize}
\end{theorem}
Recall that for every compact Bruhat section $s$,
the map $g \in s(\mathcal{F}_s)NAM \mapsto k_{\mathcal{I}}(g) \in K$ reads in $\mathcal{B}_s$ coordinates for the source and target as
\begin{align*}
\mathcal{F}_s^{(2)} \times AM & \longrightarrow \mathcal{F}_s \times M \\
(\xi,\eta \; ; \; x)_s & \longmapsto (\xi \; ; \; x_M)_s.
\end{align*}
Bruhat-Hopf coordinates make the following diagram commutative and equivariant for every compact Bruhat section $s$.
$$\xymatrix{
\mathcal{F}_s^{(2)} \times AM \simeq s(\mathcal{F}_s)NAM \subset G \ar[d]_M \ar[r] 
	&\mathcal{F}_s \times M \simeq s( \mathcal{F}_s) M \subset K  \ar[d]^M \\
\mathcal{H}^{-1}\big(\mathcal{F}_s^{(2)} \times \fa \big) \subset G/M \ar[r] & \mathcal{F}_s
 } $$
Let us now translate Theorem \ref{theo_GR_Gamma_orbite_K} using the restriction to $K$ of the Bruhat-Hopf coordinates given by the right side of the diagram.
\begin{corol} 
\hypG
Let $\G$ be a Zariski dense, discrete subsemigroup of $G$.
Then for every compact Bruhat section $s$, the following holds.
\begin{itemize}
\item[1)] For every element in $L_G(\G)$ of coordinates $(\xi \; ; \; x)_s \in \big( \mathcal{F}_s \cap L(\G) \big)  \times M $, 
there exists a unique element $[m] \in M/M_\G$ such that the element of coordinate $(\xi \; ; \; x)_s$ is in $L_{[m]}(\G)$.
Furthermore, the $\G$-orbit of this element $\G  (\xi \; ; \; x)_s$ is dense in $L_{[m]}(\G)$.
\item[2)] For every element in $ L_{[e_M]}(\G)$ of coordinate $(\xi \; ; \; x)_s$ and for all $m \in M$, the translate of coordinate $(\xi \; ; \; xm)_s$ is in $L_{[m]}(\G)$.
\item[3)] For every element in $L_{[m]}(\G)$ of coordinate $(\xi \; ; \; x)_s $ and for all $c \in M_\G$ then the element of coordinate $(\xi \; ; \; xc)_s$ remains in $L_{[m]}(\G)$.
\end{itemize} 
\end{corol}

Denote by $\widetilde{\Omega}_G$ the preimage in $G$ of $\widetilde{\Omega}\subset G/M$ by the projection $G \rightarrow G/M$. 
It is a closed, left $\G$-invariant and right $AM$-invariant subset of $G$.
For every compact Bruhat section $s \in k(\mathcal{F})$, the intersection $\widetilde{\Omega}_G \cap s(\mathcal{F}_s)NAM$ reads in Bruhat-Hopf coordinates as
$$ \mathcal{B}_s \big(\widetilde{\Omega}_G \cap s(\mathcal{F}_s)NAM \big) = \big( L^{(2)}(\G) \cap \mathcal{F}_s^{(2)} \big) \times   AM .$$
In other words, every element of coordinate $(\xi,\eta \; ; \; x)_s  \in L^{(2)}(\G) \times AM$ with $\xi \in \mathcal{F}_s$ is in $\widetilde{\Omega}_G$.
The previous Theorem and left side of the diagram allow us to partition it into closed left $\G$-invariant and right $AM_\G$-invariant subsets.
To simplify notations, for every $x \in AM$, we denote by $x_M$ its projection in $M$.
\begin{defin}\label{defin_invariant-G}
For every $m\in M$, we denote by
$\widetilde{\Omega}_{[m]}:= L_{[m]}(\G)AN \cap \widetilde{\Omega}_G$ and $\Omega_{[m]}:=  \G \backslash \widetilde{\Omega}_{[m]}.$
In other words, $\widetilde{\Omega}_{[m]}$ is the subset of elements of coordinate $(\xi, \eta \; ; \; x)_s \in L^{(2)}(\G) \times AM$ whose compact Iwasawa projection of coordinate $(\xi\; ; \; x_M)_s$ is in $L_{[m]}(\G)$, for every suitable compact Bruhat section $s$. 
\end{defin}

\begin{prop}\label{prop_transitivite-G}
\hypG
Let $\G$ be a Zariski dense, discrete subgroup of $G$.
Then the following holds.
\begin{itemize}
\item[(a)] The left $\G$-invariant and
right $AM_\G$-invariant subsets $\widetilde{\Omega}_{[m]}$ form a partition of $ \widetilde{\Omega}_G$, i.e.
$$ \widetilde{\Omega}_G = \bigsqcup_{[m] \in M/M_\G} \widetilde{\Omega}_{[m]}.$$
\item[(b)] For every $m \in M$, then $\widetilde{\Omega}_{[m]}= \widetilde{\Omega}_{[e_M]}m$.
\item[(c)] For all $[m]\in M/M_\G$, the dynamical systems $(\Omega_{[m]},\phi_\theta^t)$ and $(\Omega_{[e_M]},\phi_\theta^t)$ are conjugated.
\end{itemize}
\end{prop}
\begin{proof}
The left $\G$-invariance in (a) is a consequence of the first point of Theorem \ref{theo_GR_Gamma_orbite_K} and of the left $\G$-invariance of $\widetilde{\Omega}_G$.
It also follows from the same point that the subsets $\big(\widetilde{\Omega}_{[m]} \big)_{[m] \in M/M_\G}$ form a partition of $\widetilde{\Omega}_G$.
The right $AM_\G$-invariance is due to the right $M_\G$-invariance of $L_{[m]}(\G)$ and the properties of the Bruhat-Hopf coordinates given by Proposition \ref{prop_cocycle} and Proposition \ref{prop_BH_K}.

Point (b) is a direct consequence of the second point of Theorem \ref{theo_GR_Gamma_orbite_K} and the compatibility of the Bruhat-Hopf coordinates with the compact Iwasawa projection.

Point (c) follows from the commutativity of the right action by multiplication by $M$ with that of $A$, because every element of $M$ commute with every element of $A$.
\end{proof}

\section{Decorrelation}
Let $\G$ be a Zariski dense subgroup of $G$.
We construct a pair of points $(\xi_1,\check{\xi}_1) \in L^{(2)}(\G)$ and show that there exists $(r,\varepsilon)$-loxodromic elements in $\G$ of attracting and repelling points in an $\varepsilon$-neighbourhood of these points and whose signed cocycle are dense in an $M_\G$-invariant set.

Consider the family of equicontinuity constants $\delta_{r,\varepsilon}$ of Definition \ref{defin_d-r-eps}.
To simplify notations, we introduce the family of constants 
$$\tilde{\delta}_{r,\varepsilon}:= (8 \dim \fa + 4 \dim M_0 +5 ) \delta_{r,\varepsilon}.$$
In this section, we prove the following Proposition.
\begin{prop}\label{prop_decorrelation}
Let $G$ be a real linear, connected, semisimple Lie group of non-compact type (i.e. without compact factors.)
\emph{Assume that $M_0$ is abelian}.
Let $\Gamma$ be a Zariski dense subsemigroup of $G$.
Then there exists
\begin{itemize}
\item[1)] $( \xi_1,\check{\xi}_1) \in L^{(2)}(\Gamma)$,
\item[2)] a real positive number
$$0 < r_1 \leq \frac{1}{6} \mathrm{d}(\xi_1, \partial \mathsf{b}(\check{\xi}_1) ),$$ 
\end{itemize}
such that for all $r \in (0,r_1]$ and $\varepsilon \in (0,r]$, for any choice of compact Bruhat sections $c_1,\check{c}_1$ with
$$ B(\xi_1, r)\subset \mathcal{F}_{c_1}  \text{ and  } \mathcal{V}_{6r}\big( \partial \mathsf{b}(\check{\xi}_1) \big)^\complement \subset \mathcal{F}_{\check{c}_1}$$
there exists a finite family $(g_i)_{i\in I} \subset \G$ and a point $a_{r,\varepsilon} \in A$ that satisfy the following conditions.
\begin{itemize}
\item[$\dagger$] For all $i\in I$, the element $g_i$ is $(2r,2\varepsilon)$-loxodromic with
$$ (g_i^+,g_i^-) \in B(\xi_1, \varepsilon) \times B(\check{\xi}_1,\varepsilon) .$$
\item[$\ddagger$] For all $\eta \in \mathcal{V}_{6r}\big( \partial \mathsf{b}(\check{\xi}_1) \big)^\complement$ and $(\eta_i)_{i\in I} \subset B(\eta, \varepsilon)$, the family 
$ \lbrace \beta_{c_1,\check{c}_1}(g_i, \eta_i) \rbrace_{i\in I} $
is $\tilde{\delta}_{r,\varepsilon}$-dense in $a_{r,\varepsilon}\mathscr{R}_{c_1,\check{c}_1}(\check{\xi}_1; \xi_1,\eta) M_\Gamma$ i.e. 
$$a_{r,\varepsilon} \mathscr{R}_{c_1,\check{c}_1}(\check{\xi}_1; \xi_1,\eta)  M_\Gamma \subset  \cup_{i\in I} B( \beta_{c_1,\check{c}_1}(g_i, \eta_i), \tilde{\delta}_{r,\varepsilon} ) .$$
\end{itemize}
\end{prop}

In the first paragraph, we construct in Lemma \ref{lem_decorrelation_discret} infinitely many elements (as products of loxodromic elements) whose cocycle will hit all the connected components of $AM_\G$.

In the second paragraph, we decorrelate in an $M_0$ orbit of $AM_\G$ that projects into a convex cone of non-empty interior of $\mathfrak{a}^{++}$.
More specifically, we construct in Lemma \ref{lem_decorM}:
\begin{itemize}
\item[(a)] a convex cone of non-empty interior $\mathcal{C}_0 \subset \fa^{++}$,
\item[(b)] a pair of transverse points $(\xi_0,\check{\xi}_0) \in L^{(2)}(\G)$,
\item[(c)]  a real positive number $r_0>0$.
\end{itemize}
for which there exists, for all $0 < \varepsilon \leq r \leq r_0$ an $(r,\varepsilon)$-Schottky generating family $F_{r,\varepsilon}=(\g_1,...,\g_l)$, containing at most $4 \dim \fa + 2 \dim M_0$ loxodromic elements, such that the generalized Jordan projection of the elements of the form 
$$ \lbrace \g_1^{n_1}... \g_l^{n_l} \; \vert \; n_1, ... n_l \geq 1 \rbrace$$
are $l\delta_{r,\varepsilon}$-dense in a translate in $AM_\G$ of $\exp (\mathcal{C}_0) M_0$.
Note that the constants $\delta_{r,\varepsilon}$ become as small as we want as $\varepsilon$ goes to $0$.

In the third paragraph, we prove the Proposition \ref{prop_decorrelation} by combining the previous Lemmata with an overlapping cone argument.

\subsection{The connected components of $AM_\G$}
Denote by $p$ the integer such that $M_\Gamma / M_0$ is isomorphic to $(\mathbb{Z}/2 \mathbb{Z})^{p}$.
Since $M/M_0$ is abelian, the projection in $M/M_0$ of every signed Jordan projection does not depend on the choice of the cross-section.
The following Lemma does not require that $M_0$ is abelian.
\begin{lem}\label{lem_decorrelation_discret}
Let $G$ be a real linear, connected, semisimple Lie group of non-compact type (i.e. without compact factors).
Let $\Gamma$ be a Zariski dense subgroup of $G$.
Denote by $p$ the integer such that $M_\Gamma / M_0$ is isomorphic to $(\mathbb{Z}/2 \mathbb{Z})^{p}$ and by $\pi_{M/M_0}: AM \rightarrow M/M_0$ the projection.

Then for all $\xi_0 \in L_+(\Gamma),$ there exists $h_1,...,h_{p} \in \Gamma^{lox}$ such that taking the notation $h_0^+ := \xi_0,$ the following holds.
\begin{itemize}
\item[(i)] For every choice of cross-sections $s_1,...,s_p$ such that $h_i^+ \in \mathcal{F}_{s_i}$ for all $1 \leq i \leq p$,
the set $\lbrace \pi_{M/M_0}(\mathscr{L}_{s_i}(h_i) ) \rbrace_{1 \leq i \leq p}$ forms a basis of the vector space $M_\Gamma /M_0$.
\item[(ii)] For all $1 \leq i \leq p$, the pair $(h_{i-1}^+,h_i^-)\in L^{(2)}(\Gamma)$ is transverse.
\item[(iii)] Assume now that $s_0, s_p$ are compact Bruhat sections of respective domains $\mathsf{b}(h_1^-)$ and $\mathsf{b}(h_p^-)$, then there exists $m_p \in M$ and a large integer $N \in \mathbb{N}$ such that for all $\nu \in \lbrace 0,1 \rbrace^p$, for all $n\geq N$, 
$$\pi_{M/M_0} \Big(  \beta_{s_p m_p, s_0 } (h_p^{2n+\nu_p} ... h_1^{2n+\nu_1},\xi_0) \Big)=\nu.$$
\end{itemize}
\end{lem}

For the first step of the proof of this Lemma we use the non-arithmeticity Corollary \ref{corol_GR_non-arithm} to choose $p$ loxodromic elements in $\G$.
We order them. 
For the second step, since the repelling point of the $i$th element is not necessarily transverse to the attracting point of the $i-1$th term, we conjugate inductively these elements.
Thanks to the Fact below, the abelianised Jordan projection of the conjugated element will remain in the same connected component of $AM_{\G}$. 
To obtain the third point, we use the explicit formula of the cocycle given by Proposition \ref{prop-magic-cocycle} and the cocycle relation and combine it with a Ping-Pong argument.
Finally, the corrective term $m_p \in M$ of the cross-section is chosen using the Definition \ref{defin-ratio} of the ratio maps.

\begin{fait}\label{fait_jordan-sign-ab}
\hypG
Then for all $u\in G$ and all loxodromic element $g \in G$, the conjugate $ugu^{-1}$ is loxodromic 
of attracting point $ug^+$ and basin of attraction $u\mathsf{b}(g^-)= \mathsf{b}(ug^-)$.
Furthermore,
$$ \mathscr{L}^{ab}(ugu^{-1})= \mathscr{L}^{ab}(g).$$ 
\end{fait}
\begin{proof}
By Proposition \ref{prop-lox-bassin} a loxodromic element $g$ has attracting point $g^+$ in $\mathcal{F}$ and its basin of attraction is the Bruhat cell opposite to its repelling point $\mathsf{b}(g^-)$.
By Fact \ref{fait_lox}, consider $h_g \in G$ such that $\mathscr{L}_{[g^-]}(g)=h_g^{-1} gh_g$.
Since the Jordan projection is invariant by conjugation, $ugu^{-1}$ is also loxodromic and diagonalised by $uh_g$. 
Therefore, its attracting point is $ug^+$ of basin of attraction $\mathsf{b}(ug^-)$.
The abelian signed Jordan projection relation comes from Fact \ref{fait_lox} by choosing to compute $\mathscr{L}_{u \cdot[g^-]}(ugu^{-1})=\mathscr{L}_{[g^-]}(g)$ and then using Fact \ref{fait_jordan-sign-ab} to argue that the abelian signed Jordan projection does not depend on the choice of the cross-sections.
\end{proof}

\begin{proof}[Proof of Lemma \ref{lem_decorrelation_discret}]
Since $\mathbb{Z}/2 \mathbb{Z}$ is a field, $M_\G/M_0$ is a vector field over it. 
By Corollary \ref{corol_GR_non-arithm}, the abelian signed Jordan projection of $\G^{lox}$ spans $AM_\G^{ab}$ i.e. $AM_\G^{ab}= \overline{ \langle \mathscr{L}^{ab}(\G^{lox}) \rangle }.$ 
Furthermore, because $M_0$ is a closed normal subgroup of $M$ containing the latter's commutator subgroup, we deduce $M_\G^{ab}/M_0^{ab}=M_\G/M_0$.
Using that this is a discrete vector space and projecting in $M^{ab}$ the abelian signed Jordan projection, we write
$$M_\G/M_0=M^{ab}_\G/M^{ab}_0= \Big\langle
  \pi_{M^{ab}/M_0^{ab}}  \big( \mathscr{L}^{ab}(\G^{lox}) \big)  
  \Big\rangle.$$
The left and middle sides are $\mathbb{Z}/2\mathbb{Z}$ vector space of dimension $p$. 
The right hand side provides us with a generating set of the vector space, we extract a basis from it.
Hence there exists $g_1,...,g_p \in \G^{lox}$ such that 
$$ M_\G/M_0= \Big\langle
\pi_{M^{ab}/M_0^{ab}}  \big( \mathscr{L}^{ab}(g_p) \big), ... , 
\pi_{M^{ab}/M_0^{ab}}  \big( \mathscr{L}^{ab}(g_1) \big)
\Big\rangle .$$
Now using that $M^{ab}_\G/M^{ab}_0=M_\G /M_0$, we deduce for every suitable choice of compact Bruhat sections $b_1,...,b_p$, that
$$ M_\G/M_0= \Big\langle
\pi_{M/M_0}  \big( \mathscr{L}_{b_p}(g_p) \big), ... , 
\pi_{M/M_0}  \big( \mathscr{L}_{b_1}(g_1) \big) 
\Big\rangle .$$
By the above Fact \ref{fait_jordan-sign-ab}, condition (i) holds for every family $h_1,...,h_p$ such that for every $i=1,...,p$ the element $h_i$ is a conjugate of $g_i$.  

Let us now construct $h_1,...,h_p$.
Set $u_0:=e_G$ and $g_0^+:= \xi_0$.
We are going to choose by induction $u_1,...,u_p \in \G$ such that for every $i=1,...,p$, 
$$\big(  u_{i-1}^{-1}g_{i-1}^+ , u_i^{-1} g_i^- \big) \in \Big( L_+(\G) \times L_-(\G) \Big) \cap \mathcal{F}^{(2)}. $$ 
Repelling points of loxodromic elements lie in $L_-(\G)$ i.e. for every $i=1,...,p$, 
$$g_i^- \in L_-(\G).$$
By minimality of the action of $\G^{-1}$ on $L_-(\G)$ and because there are no isolated points in this subset, we choose $u_1\in \G$ such that
$u_1^{-1} g_1^-$ also lies in the Bruhat cell opposite to $\xi_0$, meaning that $u_1^{-1}g_1^- \in L_-(\G) \cap \mathsf{b}(\xi_0).$
By Proposition \ref{prop_transverse_pairs}, we deduce the first step $(\xi_0, u_1^{-1}g_1^-) \in L^{(2)}(\G).$
Using the same minimality arguments on the action of $\G^{-1}$ on $L_-(\G)$, we proceed as such to construct $u_i$ given $u_1,...,u_{i-1}$ such that $\big(  u_{i-1}^{-1}g_{i-1}^+ , u_i^{-1} g_i^- \big) \in L^{(2)}(\G)$.
Now that $u_1,...,u_p \in \G$ are chosen, we set for every $i=1,...,p$
$$h_i:=u_i^{-1}g_iu_i.$$ 
By the above Fact \ref{fait_jordan-sign-ab}, condition (i) holds. 
Furthermore, because $\G$ is a subgroup, every $h_i$ is a loxodromic element of $\G$ with
$$(h_i^+, h_i^-)=(u_i^{-1}g_i^+, u_i^{-1}g_i^-).$$ 
The family $h_1,...,h_p$ verify condition (ii) by construction of the $u_i$.

Let us now check condition (iii).
Choose $s_1,...,s_p$ compact Bruhat sections of respective domains $\mathsf{b}(h_1^-),...,\mathsf{b}(h_p^-)$ and set $s_1=s_0$.
For all $n_1,...,n_p \geq 1$, denote by $\underline{n}:=(n_1,...,n_p)$ and for all $i=1,...,p$ we set 
$$\xi_{i,\underline{n}}:=h_i^{n_i}...h_1^{n_1} \xi_0 .$$

Let us compute the cocycle $\beta_{s_p,s_0}(h_p^{n_p}...h_1^{n_1},\xi_0)$. 
We want to understand in which connected component of $AM$ it takes value.
Apply the cocycle relation, condition (ii) ensures the maps are well defined, then recognize the $\xi_{i,\underline{n}}$.
Finally apply the identity \eqref{eq-magic-cocycle} of Proposition \ref{prop-magic-cocycle} between the cocycle of loxodromic elements and ratio maps on each term using that the domain of $s_i$ is $\mathsf{b}(h_i^-)$ for every $i=1,...,p$.
\begin{align*}
\beta_{s_p,s_0}(h_p^{n_p}...h_1^{n_1},\xi_0)
&= 
\beta_{s_p,s_{p-1}}(h_p^{n_p}, h_{p-1}^{n_{p-1}}...h_1^{n_1}\xi_0)
\; ... \;
\beta_{s_1,s_0}(h_1^{n_1},\xi_0)
 \\
 &= 
\beta_{s_p,s_{p-1}}(h_p^{n_p}, \xi_{p-1,\underline{n}})
\; ... \;
\beta_{s_2,s_1}(h_2^{n_2}, \xi_{1,\underline{n}})
\beta_{s_1,s_0}(h_1^{n_1},\xi_0) \\
&= \hspace*{1cm} \mathscr{R}_{s_p}(h_p;\xi_{p, \underline{n}})^{-1} \quad \mathscr{L}_{s_p}(h_p)^{n_p} \quad
\mathscr{R}_{s_p,s_{p-1}}(h_p;\xi_{p-1, \underline{n}}) \; ...\\
& \hspace*{1cm} ... \quad
\mathscr{R}_{s_2}(h_2;\xi_{2, \underline{n}})^{-1} \quad \mathscr{L}_{s_2}(h_2)^{n_2} \quad
\mathscr{R}_{s_2,s_{1}}(h_2;\xi_{1, \underline{n}}) \\
&\hspace*{1.8cm}
\mathscr{R}_{s_1}(h_1;\xi_{1, \underline{n}})^{-1}
\quad \mathscr{L}_{s_1}(h_1)^{n_1} \quad
\mathscr{R}_{s_1,s_0}(h_1;\xi_0).
\end{align*}

Condition (i) allow us to deduce that the products of the middle terms $\mathscr{L}_{s_p}(h_p)^{n_p}... \mathscr{L}_{s_1}(h_1)^{n_1}$ take value in the connected component of $AM_\G$ corresponding to the projection of $\underline{n}$ in $\big(\mathbb{Z}/2\mathbb{Z}\big)^p$ that we denote by $\nu$. Then
\begin{equation}\label{eq_decor-discret}
\pi_{M/M_0}\big( \mathscr{L}_{s_p}(h_p)^{n_p}... \mathscr{L}_{s_1}(h_1)^{n_1} \big)=\nu.
\end{equation}
Note that this equation does not depend on the choice of compact Bruhat section $s_1,...,s_p$ of same domains.

It remains to control the connected components of $AM$ in which the ratio terms take value.
First, by a Ping-Pong argument, we choose a large integer $N$ which will allows us to control the sequence $(\xi_{i,\underline{n}})_{1 \leq i \leq p}$.
Then we slightly modify the choice of $s_1,...,s_p$ while preserving their domains.
Lastly, we check that under these modifications  the ratio terms are $AM_0$ valued.

Let us start by the Ping-Pong argument.
For all $i=2,...,p$ denote by $\mathsf{b}^0(h_i^-,h_{i-1}^-)$ the connected component of $\mathsf{b}(h_i^-) \cap \mathsf{b}(h_{i-1}^-) $ containing $h_{i-1}^+$.
By condition (ii) then $\xi_0 \in \mathsf{b}(h_1^-)$ and $h_1^+ \in \mathsf{b}(h_2^-)$.
By Proposition \ref{prop-lox-bassin} applied on  the loxodromic element $h_1$, there exists a large integer $N_1 \geq 1$ such that for every $n_1\geq N_1$, the element $h_1^{n_1}\xi_0$ is sufficiently close to $h_1^+$ to satisfy
$$ h_1^{n_1}\xi_0 =\xi_{1,n_1} \in \mathsf{b}^0(h_2^-, h_1^-).$$
Assume for any $i=1,...,p$ the following induction hypothesis, that there exists a large integer $N_{i-1}$ such that for every $\underline{n} \in \big( [N_{i-1},+\infty)\cap \mathbb{N} \big)^p $ and every $j=1,...,i-1$ 
$$ \xi_{j,\underline{n}} \in \mathsf{b}^0(h_{j+1}^-,h_j^-) .$$
In particular $\xi_{i-1,\underline{n}} \in \mathsf{b}(h_i^-)$. 
Also, by condition (ii) then $h_i^+ \in \mathsf{b}(h_{i+1}^-)$.
As before, we apply Proposition \ref{prop-lox-bassin} on $h_i$ to choose a large integer $N_i \geq N_{i-1}$ such that for all $\underline{n} \in \big( [N_{i},+\infty)\cap \mathbb{N} \big)^p $,
$$ h_i^{n_i}\xi_{i-1,\underline{n}}=\xi_{i,\underline{n}} \in \mathsf{b}^0(h_{i+1}^-,h_i^-).  $$  
Since $N_i$ is larger that $N_{i-1}$, the induction hypothesis is inherited for every $j=1,...,i$ i.e. $\xi_{j,\underline{n}}\in \mathsf{b}^0(h_{j+1}^-, h_j^-).$
Hence, by induction, there exists a large integer $N \geq 1$ such that for all $\underline{n} \in \big( [N,+\infty)\cap \mathbb{N} \big)^p $ and all $i=1,...,p$
\begin{equation}\label{eq_decorping-pong}
\xi_{i,\underline{n}} \in \mathsf{b}^0(h_{i+1}^-,h_i^-).
\end{equation}

Now that the large integer $N$ is chosen, assume that $\underline{n} \in \big( [2N,+\infty) \cap \mathbb{N} \big)^p$.
Let us now modify the sections by right multiplication by elements of $M$ and prove that the ratio terms for the new family of compact Bruhat section take value in $AM_0$.
Recall the Definition \ref{defin-ratio} of the ratio map.
$$\mathscr{R}_{s_i,s_{i-1}}(h_i;\xi_{i-1, \underline{n}})= \mathscr{T}_{s_i,[h_i^-]}(h_i^+)\mathscr{T}_{[h_i^-],s_{i-1}}(\xi_{i-1, \underline{n}}).$$
By Definition \ref{defin-transfer} of the transition functions, the domain of $\mathscr{R}_{s_i,s_{i-1}}(h_i;.)$ is $\mathsf{b}(h_i^-)\cap \mathsf{b}(h_{i-1}^-)$.
Set $m_0=e_M$.
By induction, we multiply $s_1,...,s_p$ on the right by elements $m_1,...,m_p\in M$ such that for every $i=1,...,p$ the restriction to the connected component containing $h_{i-1}^+$ of the map
\begin{align*}
\mathsf{b}^0(h_i^-,h_{i-1}^-) &\longrightarrow AM \\
\xi_{i-1} &\longmapsto \mathscr{R}_{s_i. m_i, s_{i-1}.m_{i-1}}(h_i;\xi_{i-1})
\end{align*}
takes value in $AM_0$.
In particular, by choice of $N$ such that condition \eqref{eq_decorping-pong} holds, we deduce that all $\mathscr{R}_{s_i.m_i,s_{i-1}.m_{i-1}}(h_i;\xi_{i-1, \underline{n}})$ term take value in $AM_0$.
Replacing them in the cocycle expression, we write
\begin{align*}
\beta_{s_p.m_p,s_0}(h_p^{n_p}...h_1^{n_1},\xi_0)
&= \mathscr{R}_{s_p.m_p}(h_p;\xi_{p, \underline{n}})^{-1} \; \mathscr{L}_{s_p.m_p}(h_p)^{n_p} \;
\mathscr{R}_{s_p.m_p,s_{p-1}.m_{p-1}}(h_p;\xi_{p-1, \underline{n}}) \; ...\\
& \hspace*{1.5cm} ...\;
\mathscr{R}_{s_1.m_1}(h_1;\xi_{1, \underline{n}})^{-1}
\; \mathscr{L}_{s_1.m_1}(h_1)^{n_1} \;
\mathscr{R}_{s_1.m_1,s_0}(h_1;\xi_0).
\end{align*}
Let us now prove that the left hand terms of the form $\mathscr{R}_{s_i.m_i}(h_i; \xi_{i,\underline{n}})^{-1}$ take value in $AM_0$.
Recall that, 
$$\mathscr{R}_{s_i.m_i} (h_i; \xi_{i,\underline{n}})= \mathscr{T}_{s_i.m_i,[h_i^-]}(h_i^+)\mathscr{T}_{[h_i^-],s_i.m_i}(\xi_{i,\underline{n}}). $$
Using that the domain of $s_i.m_i$ is $\mathsf{b}(h_i^-)$, we deduce that $\mathscr{R}_{s_i.m_i}(h_i;.)$ is well defined on it.
Furthermore, by Proposition \ref{prop_tranfer}, (iii)  
$$\mathscr{T}_{s_i.m_i,[h_i^-]}^{-1}=\mathscr{T}_{[h_i^-],s_i.m_i}.$$ 
Hence by continuity of the transition functions defined in a connected set, we deduce that the continuous maps
$\xi_i \in \mathsf{b}(h_i^-) \mapsto \mathscr{R}_{s_i.m_i}(h_i;\xi_i)^{-1}$
take value in $AM_0$.

Finally, since all ratio terms take value in $AM_0$ and by equation \eqref{eq_decor-discret}, we deduce condition (iii) that for all $\nu \in \lbrace 0,1 \rbrace^p$, all $\underline{n}$ of the form $(2n+\nu_i)_{1 \leq i \leq p}$ such that $n\geq N$, 
$$\pi_{M/M_0} \big( \beta_{s_p.m_p,s_0}(h_p^{n_p}...h_1^{n_1},\xi_0) \big) = \pi_{M/M_0} \big( \mathscr{L}_{s_p}(h_p)^{n_p}... \mathscr{L}_{s_1}(h_1)^{n_1} \big)=\nu. $$
\end{proof}

\subsection{The connected component $AM_0$}
\begin{lem}\label{lem_decorM}
Let $G$ be a real linear, connected, semisimple Lie group of non-compact type (i.e. without compact factors.)
\emph{Assume that $M_0$ is abelian}.
Let $\Gamma$ be a Zariski dense subsemigroup of $G$.
Then there exists 
\begin{itemize}
\item[(a)] a convex cone of non empty interior $\mathcal{C}_0$,
\item[(b)] a pair of transverse points $(\xi_0,\check{\xi}_0) \in L^{(2)}(\G),$
\item[(c)] a real positive number $r_0>0$,
\end{itemize}
such that for all $r \in (0,r_0]$ and $\varepsilon \in (0,r]$
 and any Bruhat section $s_0$ of domain $\mathsf{b}(\check{\xi}_0),$ \\
there exists $F_{r,\varepsilon}\subset \Gamma$ and $x_{r,\varepsilon} \in AM_\Gamma$ such that the following holds.
\begin{itemize}
\item[$\heartsuit$] $F_{r,\varepsilon}$ is a finite subset of at most $4 \dim \fa + 2 \dim M_0$ elements.
\item[$\clubsuit$] $F_{r,\varepsilon}$ is a subset of a strong $(r,\varepsilon)$-Schottky Zariski dense subsemigroup.
\item[$\diamondsuit$] There exists an ordering of $F_{r,\varepsilon}=(\gamma_1,...,\gamma_l)$ such that $\gamma_1^-= \check{\xi}_0$ and $\gamma_l^+=\xi_0$, for which every element of the form $w=\gamma_l^{n_l}... \gamma_1^{n_1}$ with $n_1,...,n_l \geq 1$, satisfies
$$ (w^+,w^-) \in B(\xi_0, \varepsilon) \times B(\check{\xi}_0,\varepsilon).$$
\item[$\spadesuit$] For such an ordering, the set
$$ \mathscr{L}_{s_0} \big(  \lbrace \gamma_l^{n_l}... \gamma_1^{n_1} \; \vert \; n_1,...,n_l \geq 1 \rbrace \big) $$
is $l_{AM_0}\delta_{r,\varepsilon} $-dense in 
$\exp(\mathcal{C}_0)x_{r,\varepsilon}M_0$, where $l_{AM_0}:=8 \dim \fa +4 \dim M_0 +1$.
\end{itemize}
The family of constants
$\delta_{r,\varepsilon}$ is given in Definition \ref{defin_d-r-eps}, for every $r>0$, they converge to $0$ when $\varepsilon$ goes to $0$.
\end{lem}

The first step of the proof is given by the following Lemma, which is a consequence of \cite[Proposition 4.3]{benoist1997proprietes}.
We give a reference for a proof.
The last steps involve the non-arithmeticity of Corollary \ref{corol_GR_non-arithm} and density Lemmata \ref{lemme_densite2}, \ref{lemme_densite3_cone} of the appendix.
These statements require that $M_0$ is abelian.
	\begin{lem}[Lemme 5.6 \cite{dang_glorieux_2020}]\label{lemme_schottky_theta}
	\hypG Let $\Gamma\subset G$ be Zariski dense subsemigroup. 
	For all $\theta$ in the interior of the limit cone $\mathcal{C}(\G)$, there exists a finite set $S\subset \G$, a positive number $r_0>0$ such that
		\begin{itemize}
		\item[(i)] $\theta$ is in the interior of the convex cone $C(\lambda (S)):= \sum_{g\in S}\R_+ \lambda(g)$,
		\item[(ii)] the elements of $\lambda(S)$ form a basis of $\fa$, 
		\item[(iii)] for all $r\in (0,r_0]$ and $\varepsilon \in (0,r]$, there exists an integer $N>0$ such that for all $n\geq N$, the family $S_n:= (g^n)_{g\in S}$ spans a Zariski-dense strong $(r,\varepsilon)$-Schottky semigroup of $\G$. 
		\end{itemize}
	\end{lem}

\begin{proof}[Proof of Lemma \ref{lem_decorM} ]
First fix $\theta$ in the interior of the limit cone.
Apply now Lemma \ref{lemme_schottky_theta}.
Set 
$ \mathcal{C}_0:= \sum_{g \in S} \R_+ \lambda(g) .$
By (i), it is indeed a convex cone of non-empty interior.
Let us now order the elements $S:=(g_1,..., g_{r_G})$ where $r_G=\dim \fa$ by (ii).
By (iii), for any integer $n$ sufficiently large, $S_n$ spans a strong $(r,\varepsilon)$-Schottky Zariski dense subsemigroup.
We deduce that $g_{r_G}^+$ is in the basin of attraction of $g_1$, meaning that $g_{r_G}^+\in \mathsf{b}(g_1^-)$, which by Proposition \ref{prop_transverse_pairs} is the same as $(g_{r_G}^+,g_1^-)\in L^{(2)}(\G)$.

Let $r \in (0,r_0]$ and $\varepsilon \in (0,r]$,
 fix a compact Bruhat section $s_0$ of domain $\mathsf{b}(\check{\xi}_0)$. 
Let us choose $F_{r,\varepsilon} \subset \G$. 
Consider a large integer $N$ such that for every $n\geq N$, the subset $S_n$ spans a Zariski dense, strong $(r,\varepsilon)$-Schottky subsemigroup.

By Theorem \ref{thm_GR_M_G}, the group $M_\G/M_0$ is isomorphic to $(\mathbb{Z}/2\mathbb{Z})^p$.
Consequently, for every element $m\in M_\G$, its square $m^2$ is in $M_0$. 
In particular, for every loxodromic element $\g\in \G$ and any suitable compact Bruhat section $s$ such that $\g^+ \in \mathcal{F}_s$, 
$$\mathscr{L}_s(\g^2)= \big( \mathscr{L}_s(\g) \big)^2 \in AM_0.$$ 
Since $M_0$ is abelian and a normal subgroup of $M_\G$, we deduce that the multiplicative Jordan projection of squares does not depend on the choice of $s$ and coincides with $\mathscr{L}^{ab}$.
We therefore remove the subscript. 
Denote by $\G_n$ the Zariski dense subsemigroup generated by $S_{2n}.$
By Corollary \ref{corol_GR_non-arithm} and using that $M_0$ is abelian,
$$ \overline{\langle \mathscr{L}^{ab}(\G_n) \rangle}=AM_{\G_n}^{ab}\supset AM_0. $$
Let us prove that the subset of squares $\mathscr{L}(\G_n)^2$ spans a dense subgroup of $AM_0$. 
Every element $x\in AM_0$ admits a square root that we denote by $\sqrt{x} \in AM_0.$
Now we approximate it in $\langle \mathscr{L}^{ab}(\G_n) \rangle$.
For all $\delta >0 $, there exists a finite number of integers $(k_j)_{j\in J} \subset \mathbb{Z}$ and a finite number of elements $(\g_j)_{j\in J}$ such that
$$ \sqrt{x} \in B\Big( \prod_{j\in J} \mathscr{L}^{ab}(\g_j)^{k_j} ,\sqrt{\delta} \Big). $$
Taking the squares, we obtain the approximation by squares, 
$$x \in B\Big( \prod_{j\in J} \mathscr{L}^{}(\g_j)^{2k_j} ,\delta \Big).$$
Hence 
$$\overline{ \langle \mathscr{L}(\G_n)^2  \rangle }=AM_0.$$
Apply density Lemma \ref{lemme_densite2} in $AM_0$ for the family of squares $\mathscr{L}(\G_n)^2$. 
Consider $F_{r,\varepsilon}'$ of at most $3 \dim \fa + 2 \dim M_0$ elements such that the subgroup spanned by squares $\mathscr{L}(F_{r,\varepsilon}')^2$ is $ \delta_{r,\varepsilon}$-dense in $AM_0$.
Denote by 
$$F_{r,\varepsilon}:= S_{2n} \cup \lbrace \g^2 \; \vert \; \g \in F_{r,\varepsilon}' \rbrace.$$
The subgroup spanned by $\mathscr{L}(F_{r,\varepsilon})$ is $\delta_{r,\varepsilon}$-dense in $AM_0$.
Apply now density Lemma \ref{lemme_densite3_cone} to such a family.
There exists $v_{r,\varepsilon} \in \fa$ such that the subsemigroup generated by $\mathscr{L}(F_{r,\varepsilon})$ is $\delta_{r,\varepsilon}$-dense in 
$$  \exp \bigg(v_{r,\varepsilon}+ \sum_{\g \in F_{r,\varepsilon}} \R_+ \lambda(\g) \bigg) M_0 .$$
Now since $F_{r,\varepsilon}$ contains $S_{2n}$ and by choice of $\mathcal{C}_0$,
$$ \mathcal{C}_0 \subset  \sum_{\g \in F_{r,\varepsilon}} \R_+ \lambda(\g), $$
we deduce $\delta_{r,\varepsilon}$-density of the subsemigroup generated by $\mathscr{L}(F_{r,\varepsilon})$ in $\exp(v_{r,\varepsilon}+\mathcal{C}_0) M_0$.

Let us now compute $x_{r,\varepsilon}\in AM$.
We order $F_{r,\varepsilon}=(\g_1,...,\g_l)$ such that $\g_1:=g_1^{2n}$ and $\g_l:=g_{r_G}^{2n}$.
Fix compact Bruhat sections $s_1,...,s_l$ 
such that for every $i=1,...,l$
$$B(\g_i^+,\varepsilon) \subset \mathcal{F}_{s_i}.$$
We assume that for every $i=2,...,l$ then $\mathscr{R}_{s_i,s_{i-1}}(\g_i, \g_{i-1}^+)\in AM_0$.  
Since $\mathscr{R}_{s_i,s_{i-1}}(\g_i,.)$ restricted to the connected component of $\mathsf{b}(\g_i^-) \cap \mathcal{F}_{s_{i-1}}$ containing $\g_{i-1}^+$ takes value in a connected component of $AM$, by multiplying $s_i$ on the right by an element of $M$ one can always assume that this restricted map takes value in $AM_0$. 
Set $\mathscr{R}_i:= \mathscr{R}_{s_i,s_{i-1}}(\g_i, \g_{i-1}^+)$ with convention that $s_0=s_l$ and $\g_0=\g_l$ and
$$ x_{r,\varepsilon}:= \exp(v_{r,\varepsilon}) \mathscr{R}_l...
\mathscr{R}_2 \mathscr{R}_1. $$

Let us now check the card conditions.
Since $S_{2n}$ contains $\dim \fa$ elements and generates the strong $(r,\varepsilon)$-Schottky  $\G_n$ and by choice of $F_{r,\varepsilon}'$, the subset $F_{r,\varepsilon}$ satisfies both $\heartsuit$ and $\clubsuit$.

By choice of ordering $\g_1^-=g_1^-=\check{\xi}_0$ and $\g_l^+=g_{r_G}^+=\xi_0$.
Apply Proposition \ref{prop_crucial}, for all $n_1,...,n_l \geq 1$ the element $w=\g_l^{n_l}...\g_1^{n_1}$ is loxodromic and satisfies  $ (w^+,w^-) \in B(\xi_0, \varepsilon) \times B(\check{\xi}_0,\varepsilon).$
Hence $\diamondsuit$ is satisfied.

Furthermore, by Proposition \ref{prop_crucial} we estimate the Jordan projection
$$ \mathscr{L}_{s_l}(w) \in  \mathscr{L}_{s_l}(\g_l)^{n_l} \mathscr {R}_l ... \mathscr{L}_{s_1}(\g_1)^{n_1} \mathscr{R}_1 \; B\big(e_{AM}, 2l \delta_{r,\varepsilon}\big) .$$
Note that $\mathscr{R}_l,...,\mathscr{R}_2$ take value in the abelian group $AM_0$ by choice of $s_2,..,s_l$.
Furthermore, the $\g_i$ are squares, hence integer powers of $\mathscr{L}_{s_i}(\g_i)$ take value in $AM_0$ and we can remove the subscript.
Hence by reordering the terms in $AM_0$,
 $$ \mathscr{L}_{s_l}(w) \in \mathscr{L}(\g_l)^{n_l} ... \mathscr{L}(\g_1)^{n_1}  \mathscr{R}_l...\mathscr{R}_1 \; B\big(e_{AM}, 2l \delta_{r,\varepsilon}\big) .$$
The first part of the left hand side 
$$ \lbrace \mathscr{L}(\g_l)^{n_l} ... \mathscr{L}(\g_1)^{n_1} \; \vert \; n_1,...,n_l \geq 1 \rbrace $$ coincides with the subsemigroup of $AM_0$ generated by $\mathscr{L}(F_{r,\varepsilon})$ which is $\delta_{r,\varepsilon}$-dense in 
$$\exp(\mathcal{C}_0) \exp(v_{r,\varepsilon}) M_0.$$
We deduce that
$$ \lbrace \mathscr{L}_{s_l}(\g_l^{n_l} ...\g_1^{n_1} )\; \vert \; n_1,...,n_l \geq 1 \rbrace $$
is $(2l+1)\delta_{r,\varepsilon}$-dense in 
$$ \exp(\mathcal{C}_0) \exp(v_{r,\varepsilon})M_0 \; \mathscr{R}_l...\mathscr{R}_1=\exp(\mathcal{C}_0)  M_0 x_{r,\varepsilon},$$
using that $M_0$ centralises $A$.
Since $M_0$ is a normal subgroup, we deduce $(2l+1)\delta_{r,\varepsilon}$-density in $\exp(\mathcal{C}_0)  M_0 x_{r,\varepsilon}=\exp(\mathcal{C}_0)  x_{r,\varepsilon} M_0 .$
By $\heartsuit$, then $l\leq 4\dim \fa + 2 \dim M_0$, hence 
$$(2l+1)\delta_{r,\varepsilon}\leq (8 \dim \fa +4 \dim M_0 +1) \delta_{r,\varepsilon} $$ and condition $\spadesuit$ is satisfied.
\end{proof}

\subsection{Proof of Proposition \ref{prop_decorrelation} }

Let us first find the pair $(\xi_1,\check{\xi}_1) \in L^{(2)}(\G)$ using the previous Lemmas of this section.
Consider the pair of transverse points $(\xi_0,\check{\xi}_0)\in L^{(2)}(\G)$ given by the decorrelation in $AM_0$ Lemma \ref{lem_decorM} (b).
Apply Lemma \ref{lem_decorrelation_discret} to $\xi_0$ to reach every connected component of $AM_\G$.
There exists loxodromic elements $h_1,...,h_p \in \G^{lox}$ such that taking the notation $h_0^+ := \xi_0,$ the following holds.
\begin{itemize}
\item[(i)] For every choice of sections $s_1,...,s_p$ such that $h_i^+ \in \mathcal{F}_{s_i}$ for all $1 \leq i \leq p$,
the set $\lbrace \pi_{M/M_0}(\mathscr{L}_{s_i}(h_i) ) \rbrace_{1 \leq i \leq p}$ forms a basis of the vector space $M_\Gamma /M_0$.
\item[(ii)] For all $1 \leq i \leq p$, the pair $(h_{i-1}^+,h_i^-)\in L^{(2)}(\Gamma)$ is transverse.
\item[(iii)] Assume now that $s_0, s_p$ are compact Bruhat sections of respective domains $\mathsf{b}(h_1^-)$ and $\mathsf{b}(h_p^-)$, then there exists $m_p \in M$ and a large integer $N \in \mathbb{N}$ such that for all $\nu \in \lbrace 0,1 \rbrace^p$, for all $n\geq N$, 
$$\pi_{M/M_0} \Big(  \beta_{s_p m_p, s_0 } (h_p^{2n+\nu_p} ... h_1^{2n+\nu_1},\xi_0) \Big)=\nu.$$
\end{itemize}
Since $h_p^+$ has no reason to be transverse to $\check{\xi}_0$, we need the following choice.
By density of attracting and repelling points of loxodromic elements in $L^{(2)}(\G)$, there exists a loxodromic element $h_{p+1} \in \G^{lox}$ such that
$$\left\lbrace  
\begin{array}{c}
(h_{p+1}^+,\check{\xi}_0)  \in  L^{(2)}(\G) \\
(h_p^+, h_{p+1}^-)  \in  L^{(2)}(\G)
\end{array} \right.$$
Such a choice is always possible because there are no isolated points in the limit sets $L_\pm(\G)$.
Set now
\begin{equation}\label{eq_prop-decor1}
(\xi_1,\check{\xi}_1):=(h_{p+1}^+,\check{\xi}_0).
\end{equation}

Let us now find the positive number $r_1$. 
Consider the real number $r_0$ given by Lemma \ref{lem_decorM} (c).
We set 
$$r_0':= \inf_{1 \leq i \leq p+1} 
\left\lbrace 
\frac{1}{6} \mathrm{d}(h_{i-1}^+, \partial\mathsf{b}(h_i^-))
,\frac{1}{2} \mathrm{d}(h_i^+, \partial \mathsf{b}(h_i^-))
\right\rbrace .$$
By (ii), choice of $h_{p+1}$ and using that $h_1,...,h_{p+1}$ are loxodromic, we deduce that
both $r_0$ and $r_0'$ are positive real numbers.
This leads us to define the positive real number
\begin{equation}\label{eq_prop-decor2}
r_1:= \inf(r_0,r_0').
\end{equation}

Let $r \in (0,r_1]$ and $\varepsilon \in (0,r]$.
Fix a choice of compact Bruhat sections $c_1, \check{c}_1$ such that
$$ B(\xi_1,r)\subset \mathcal{F}_{c_1} \text{ and } \mathcal{V}_{6r}(\partial \mathsf{b}(\check{\xi}_1))^\complement \subset \mathcal{F}_{\check{c}_1} .$$

\textbf{Reaching every connected component of $AM_\G$}\\
By Proposition \ref{prop-lox-lips} on loxodromic elements $h_1,...,h_{p+1}$, there exists a large integer $N_{r,\varepsilon} \geq 1$ such that for every $n\geq N_{r,\varepsilon}$, each $h_i^n$ are $(r,\varepsilon)$-loxodromic. 

Since $\xi_0$ is in the basin of attraction of $h_1$, then by Proposition \ref{prop-lox-bassin}, we choose another integer $N_1 \geq 1$ such that for all $n\geq N_1$ large enough, $h_1^n \xi_0 \in B(h_1^+,\varepsilon)$.
Set $N_2:= \sup(N_1,N_{r_\varepsilon})$.
By a Ping-Pong argument using the dynamical properties of $(r,\varepsilon)$-loxodromic elements, this implies that for all $n_1,...,n_p \geq N_2 $, then $h_{p}^{n_p}...h_1^{n_1}\xi_0 \in B(h_p^+, \varepsilon)$.
For all $\underline{n}:=( n_{p},...,n_1)$ family of positive integers, denote by $\xi_{p,\underline{n}}:=h_p^{n_p}...h_1^{n_1}\xi_0$.
By Proposition \ref{prop-magic-cocycle}, for all $\xi_p\in B(h_p^+,\varepsilon)$, then 
$$\beta_{c_1,s_p m_p}(h_{p+1}^{2n}, \xi_p) = \mathscr{R}_{c_1}(h_{p+1}; h_{p+1}^{2n} \xi_p)^{-1} \mathscr{L}(h_{p+1})^{2n} \mathscr{R}_{c_1,s_p m_p}(h_{p+1}; \xi_p) .$$ 
Note that by choice of $\varepsilon\leq r_1$ the balls $B(h_p^+,\varepsilon)$ resp. $B(h_{p+1}^+,\varepsilon)$  are included in connected components of $\mathsf{b}(h_{p+1}^-) \cap \mathcal{F}_{s_p}$ resp. $\mathcal{F}_{c_1} \cap \mathsf{b}(h_{p+1}^-)$.
Therefore, using that $h_{p+1}^{2n}$ is $(r,\varepsilon)$-loxodromic when $n \geq N_2$, 
the restriction to $B(h_p^+,\varepsilon)$ of $\beta_{c_1,s_p m_p}(h_{p+1}^{2n},. )$ to $B(h_p^+,\varepsilon)$ is constant mod $AM_0$.
For all $n \geq N_2$, the map
\begin{align*}
\big([N_2,\infty)\cap \mathbb{N}\big)^p & \longrightarrow AM \\
\underline{n} & \longmapsto \beta_{c_1,s_0}(h_{p+1}^{2n} h_p^{n_p}...h_1^{n_1}, \xi_0 )
\end{align*}
reaches every connected components of $AM_\G$.
Indeed, by the cocycle relation 
$$\beta_{c_1,s_0}(h_{p+1}^{2n} h_p^{n_p}...h_1^{n_1}, \xi_0 )=\beta_{c_1,s_p m_p}(h_{p+1}^{2n},\xi_{p,\underline{n}}) \beta_{s_pm_p,s_0}(h_p^{n_p}...h_1^{n_1},\xi_0), $$
and by (iii), we control which connected component of $AM_\G$ the right term hits, the left term being constant mod $AM_0$ as discussed above.
Thus, for all $\nu \in \lbrace 0,1 \rbrace^p \simeq M_\G/M_0$ there exists and we choose $n_p(\nu),...,n_1(\nu) \geq N_2$ such that if we denote by 
$$ \left\lbrace \begin{array}{l}
h_{[\nu]}:=h_{p+1}^{2n} h_p^{n_p(\nu)}...h_1^{n_1(\nu)} \\
x_{[\nu]}:= \beta_{c_1,s_0}(h_{[\nu]},\xi_0)
\end{array} \right. $$
then $\pi_{AM/AM_0}(x_{[\nu]})=\nu.$

\textbf{A particular subset of loxodromic elements of $\G$ }\\
Consider now the subset $F_{r,\varepsilon}$, the point $x_{r,\varepsilon} \in AM_\G$ and the convex cone of non-empty interior $\mathcal{C}_0$ given by Lemma \ref{lem_decorM}.
They satisfy 
\begin{itemize}
\item[$\heartsuit$] $F_{r,\varepsilon}$ is a finite subset of at most $4 \dim \fa + 2 \dim M_0$ elements.
\item[$\clubsuit$] $F_{r,\varepsilon}$ is a subset of a strong $(r,\varepsilon)$-Schottky Zariski dense subsemigroup.
\item[$\diamondsuit$] There exists an ordering of $F_{r,\varepsilon}=(\gamma_1,...,\gamma_l)$ such that $\gamma_1^-= \check{\xi}_0$ and $\gamma_l^+=\xi_0$, for which every element of the form $w=\gamma_l^{n_l}... \gamma_1^{n_1}$ with $n_1,...,n_l \geq 1$, satisfies
$$ (w^+,w^-) \in B(\xi_0, \varepsilon) \times B(\check{\xi}_0,\varepsilon).$$
\item[$\spadesuit$] For such an ordering, the set
$$ \mathscr{L}_{s_0} \big(  \lbrace \gamma_l^{n_l}... \gamma_1^{n_1} \; \vert \; n_1,...,n_l \geq 1 \rbrace \big) $$
is $l_{AM_0}\delta_{r,\varepsilon} $-dense in 
$\exp(\mathcal{C}_0)x_{r,\varepsilon}M_0$ where $l_{AM_0}:=8 \dim \fa +4 \dim M_0 +1$.
\end{itemize}
We are going to choose $(g_i)_{i\in I}$ among elements of the form $h_{[\nu]}\g_l^{n_l}...\g_1^{n_1}$, where $n_l,...,n_1 \geq 1$ are integers and $\nu \in \lbrace 0,1 \rbrace^p$.

By choice of $r_1$, we deduce $\dagger$ for all elements of 
$$ \lbrace h_{[\nu]} \g_l^{n_l}...\g_1^{n_l} \; \vert \; n_1,...,n_l \geq 1 \; \text{and} \; \nu \in (\mathbb{Z}/2\mathbb{Z})^p \rbrace .$$
Meaning that all elements of the set above are $(2r,2\varepsilon)$-loxodromic with attracting and repelling points in $B(\xi_1,\varepsilon) \times B(\check{\xi}_1,\varepsilon)$.

\textbf{Cocycle estimates}\\
By equation \eqref{eq_prop-decor1} recall that $\check{\xi}_0=\check{\xi}_1=\g_1^-$.
Let $n_l,...,n_1 \geq 1$ be integers. 
Then by choice of $r\leq r_1$ and $\varepsilon\leq r$, the element $\g = \g_l^{n_l}...\g_1^{n_1}$ is $(r,\varepsilon)$-loxodromic, of attracting point in $B(\g_l^+,\varepsilon)$ and repelling point in $B(\g_1^-,\varepsilon)$.
By Proposition \ref{prop-magic-cocycle} on loxodromic element $\g$ and $\eta \in \mathcal{V}_{6r}(\partial \mathsf{b}(\check{\xi}_1))^\complement$, by  Definition \ref{defin_d-r-eps} of the equicontinuity constant $\delta_{r,\varepsilon}$,
we deduce
$$\beta_{s_0,\check{c}_1}(\g, \eta) \in \mathscr{L}_{s_0}(\g) \mathscr{R}_{s_0,\check{c}_1}(\g ; \eta) B(e_{AM},\delta_{r,\varepsilon}).$$
Now, $\mathscr{R}_{s_0,\check{c}_1}(\g ; \eta)$ is  $\delta_{r,\varepsilon}$ close to $\mathscr{R}_{s_0,\check{c}_1}(\g_1^-; \g_l^+,\eta)$.
Hence using $\xi_0=\g_l^+$ and $\check{\xi}_1=\g_1^-$, we deduce
$$\beta_{s_0,\check{c}_1}(\g, \eta) \in \mathscr{L}_{s_0}(\g) \mathscr{R}_{s_0,\check{c}_1}(\check{\xi}_1 ; \xi_0, \eta) B(e_{AM},2\delta_{r,\varepsilon}).$$
For all $\nu \in \lbrace 0,1 \rbrace^{p}$, by the cocycle relation,
$$\beta_{c_1,\check{c}_1}(h_{[\nu]}\g,\eta )=
\beta_{c_1,s_0}(h_{[\nu]}, \g \eta ) \beta_{s_0,\check{c}_1}(\g,\eta). $$
By a Ping-Pong argument on $\g_1,..., \g_l$ we deduce that $\g \eta \in B(\g_l^+,\varepsilon)$.
Similarly, the same type of argument on the $(r,\varepsilon)$-loxodromic elements $h_1^{n_1(\nu)},...,h_p^{n_p(\nu)},h_{p+1}^{2n}$ yields that 
$$\beta_{c_1,s_0}(h_{[\nu]}, \g \eta ) \in \beta_{c_1,s_0}(h_{[\nu]} ,\g_l^+ ) B(e_{AM}, \delta_{r,\varepsilon}) .$$
Using $\g_l^+=\xi_0$ and the definition of $x_{[\nu]}$, we deduce the following estimate
$$ \beta_{c_1,\check{c}_1}(h_{[\nu]}\g,\eta  )
\in x_{[\nu]} \mathscr{L}_{s_0}(\g) \mathscr{R}_{s_0,\check{c}_1}(\check{\xi}_1 ; \xi_0, \eta) B(e_{AM},3 \delta_{r,\varepsilon}) .$$
To recover the term $\mathscr{R}_{c_1,\check{c}_1}(\check{\xi}_1; \xi_1,\eta)$ as in $\ddagger$,  one can check using the definition of the Ratio maps that $ \mathscr{R}_{s_0,\check{c}_1}(\check{\xi}_1; \xi_0,\eta) = \mathscr{R}_{c_1,s_0}(\check{\xi}_1; \xi_1,\xi_0)^{-1} \mathscr{R}_{c_1,\check{c}_1}(\check{\xi}_1,\xi_1,\eta) .$
Denote by $y_0:= \mathscr{R}_{c_1,s_0}(\check{\xi}_1; \xi_1,\xi_0)^{-1}$. 
Then for all $\nu \in (\mathbb{Z}/2\mathbb{Z})^p$, all $\g \in \lbrace \g_l^{n_l}...\g_1^{n_1} \; \vert \; n_1,...,n_l\geq 1 \rbrace$ and all $\eta \in \mathcal{V}_{6r}(\partial \mathsf{b}(\check{\xi}_1))^\complement$,
\begin{equation}\label{eq-decor-idcocycle}
\beta_{c_1,\check{c}_1}(h_{[\nu]}\g,\eta  )
\in x_{[\nu]} \mathscr{L}_{s_0}(\g) y_0  \mathscr{R}_{c_1,\check{c}_1}(\check{\xi}_1 ; \xi_1, \eta) B(e_{AM},3 \delta_{r,\varepsilon}).
\end{equation}

\textbf{Overlapping cone argument}\\
Using $\spadesuit$ on the Jordan term, we deduce that for every $\eta \in \mathcal{V}_{6r}(\partial \mathsf{b}(\check{\xi}_1))^\complement$, the subset of cocycles
$$ \lbrace \beta_{c_1,\check{c}_1} (h_{[\nu]} \g_l^{n_l}...\g_1^{n_1},\eta ) \; \vert \;
n_1,...,n_l\geq 1
 \rbrace $$
is $(3 + l_{AM_0}) \delta_{r,\varepsilon}$-dense in the translated cone $x_{[\nu]}  \exp(\mathcal{C}_0) x_{r,\varepsilon} M_0 \; y_0 \mathscr{R}_{c_1,\check{c}_1}(\check{\xi}_1 ; \xi_1, \eta) .$ 
The left terms $x_{[\nu]}$ ensures that when $\nu$ varies in $(\mathbb{Z}/2\mathbb{Z})^p$, all connected components of $AM_\G$ are reached.
Denote by $\pi_A: AM \rightarrow A$ the projection.
Using that $\mathcal{C}_0$ is convex of non-empty interior, we deduce that there exists $a_{r,\varepsilon} \in A$ such that the intersection, over the number of connected components of $AM_\G$, of the projection in $A$ of these translated cones, contains $a_{r,\varepsilon} \exp (\mathcal{C}_0)$, i.e.
$$ a_{r,\varepsilon} \exp (\mathcal{C}_0) \; \subset \bigcap_{\nu \in (\mathbb{Z}/2\mathbb{Z})^p}  
\pi_A\big( x_{[\nu]}  \exp(\mathcal{C}_0) x_{r,\varepsilon} M_0 y_0 \big) .$$
Hence the disjoint union of translated cones contains $a_{r,\varepsilon} \exp(\mathcal{C}_0)M_\G$ i.e.
$$a_{r,\varepsilon} \exp(\mathcal{C}_0)M_\G \; \subset  \bigsqcup_{\nu \in (\mathbb{Z}/2\mathbb{Z})^p}  
 x_{[\nu]}  \exp(\mathcal{C}_0) x_{r,\varepsilon} M_0 \; y_0 .$$
Hence by right multiplication by $\mathscr{R}_{c_1,\check{c}_1}(\check{\xi}_1; \xi_1,\eta)$, we deduce that
$$a_{r,\varepsilon} \exp(\mathcal{C}_0)M_\G
\mathscr{R}_{c_1,\check{c}_1}(\check{\xi}_1; \xi_1,\eta)
 \; \subset  \bigsqcup_{\nu \in (\mathbb{Z}/2\mathbb{Z})^p}  
 x_{[\nu]}  \exp(\mathcal{C}_0) x_{r,\varepsilon} M_0 \; y_0  \mathscr{R}_{c_1,\check{c}_1}(\check{\xi}_1; \xi_1,\eta).$$
Using the $(3+l_{AM_0}) \delta_{r,\varepsilon}$ density of cocycles in the disjoint union on the right yields
$$a_{r,\varepsilon} \exp(\mathcal{C}_0)M_\G
\mathscr{R}_{c_1,\check{c}_1}(\check{\xi}_1; \xi_1,\eta)
 \; \subset 
\underset{n_1,...,n_l \geq 1}{ \bigcup_{\nu \in (\mathbb{Z}/2\mathbb{Z})^p}} \beta_{c_1,\check{c}_1}(h_{[\nu]}\g_l^{n_l}...\g_1^{n_1} ,\eta )B(e_{AM}, (3+l_{AM_0}) \delta_{r,\varepsilon}).$$
By compacity, we choose a finite family 
$$ (g_i)_{i\in I} \subset \lbrace h_{[\nu]} \g_l^{n_l}...\g_1^{n_l} \; \vert \; n_1,...,n_l \geq 1 \; \text{and} \; \nu \in (\mathbb{Z}/2\mathbb{Z})^p \rbrace $$
such that for all $\eta \in \mathcal{V}_{6r}(\partial \mathsf{b}(\check{\xi}_1))^\complement$,
$$a_{r,\varepsilon} M_\G \mathscr{R}_{c_1,\check{c}_1}(\check{\xi}_1; \xi_1,\eta)
 \; \subset \bigcup_{i\in I} B\big( \beta_{c_1,\check{c}_1}(g_i,\eta) ,(3+l_{AM_0}) \delta_{r,\varepsilon}\big) .$$
Set $\tilde{\delta}_{r,\varepsilon}:= (8 \dim \fa + 4 \dim M_0 +5 ) \delta_{r,\varepsilon}=(l_{AM_0}+4)\delta_{r,\varepsilon}$.
Apply Proposition \ref{prop-magic-cocycle} with Definition \ref{defin_d-r-eps} of the equicontinuity constants $\delta_{r,\varepsilon}$ for every family $(\eta_i)_{i\in I} \subset B(\eta,\varepsilon)$ to deduce $\ddagger$,
$$a_{r,\varepsilon} M_\G \mathscr{R}_{c_1,\check{c}_1}(\check{\xi}_1; \xi_1,\eta)
 \; \subset \bigcup_{i\in I} B\big( \beta_{c_1,\check{c}_1}(g_i,\eta_i) , \tilde{\delta}_{r,\varepsilon}\big) .$$
\section{Conditions for topological mixing}
We prove the following necessary and sufficient conditions.
\begin{theorem}\label{theo_mixing}
Let $G$ be a real linear, connected, semisimple Lie group of non-compact type (i.e. without compact factors) and $\Gamma$ be a Zariski dense subgroup of $G$.
For all $\theta \in \fa^{++}$, the following topological mixing conditions occur. 
\begin{itemize}
\item[(NC)] If the dynamical system $(\Omega_{[e_M]},\phi_\theta^t)$ is topologically mixing then $\theta \in \overset{\circ}{ \mathcal{C}}(\Gamma)$.
\item[(SC)] Assume that the connected component of the identity $M_0$ of $M$ is \emph{abelian}.
Then the converse is true i.e. if $\theta$ is in the interior of the Benoist cone, then the dynamical system $(\Omega_{[e_M]},\phi_\theta^t)$ is topologically mixing.
\end{itemize}
\end{theorem}

\subsection{Necessary condition}
Let $\theta \in \fa^{++}$.
We prove that if the dynamical system $(\Omega_{[e_M]}, \phi_\theta^t)$ is topologically mixing, 
then $\theta$ is in the interior of the limit cone $\mathcal{C}(\G)$.

Since this dynamical system factors $(\Omega,\phi_\theta^t)$, we deduce topological mixing of the regular Weyl chamber flow. 
Using now $\theta\in \fa^{++}$ and the necessary and sufficient condition for mixing \cite{dang_glorieux_2020}, we deduce that 
$$\theta \in \overset{\circ}{\mathcal{C}}(\G).$$ 

\subsection{Sufficient condition}

The key arguments are given by Theorem \ref{transitivity of Furstenberg boundary}, decorrelation Proposition \ref{prop_decorrelation} and the Proposition \ref{DG-prop} below.

Let $\theta \in \fa^{++}$ be in the interior of the limit cone.
We want to prove that 
for all non-empty open sets $ U, V\subset  \Omega_{[e_M]},$ there exists $T >0$ such that for every $t\geq T$, 
$$ \phi_\theta^t \big( U \big) \cap V \neq \emptyset. $$

It is equivalent to prove that 
for all non-empty open sets $ \mathcal{U}, \mathcal{V}\subset  \widetilde{\Omega}_{[e_M]},$ there exists $T >0$ such that for every $t\geq T$, 
$$ \mathcal{U} e^{t \theta} \cap \G \mathcal{V} \neq \emptyset. $$

By Theorem \ref{transitivity of Furstenberg boundary}, the action of $\G$ on $L^{(2)}(\G)$ has dense orbits.
The latter are the first and second Bruhat-Hopf coordinates of $\widetilde{\Omega}_{[e_M]}$.
Using that left and right actions commute, we align $\mathcal{U}$ and $\mathcal{V}$ in the same $AM$ orbit as a right $AM$-invariant subsets given by Proposition \ref{prop_decorrelation}: of first and second Bruhat-Hopf coordinates in a neighbourhood of $(\xi_1, \check{\xi}_1)$.

The Proposition \ref{DG-prop}, applied to $\theta$ and the neighbourhood of $(\xi_1,\check{\xi}_1)$, proves the density in affine half lines of direction $\theta$, of the Jordan projection of loxodromic element whose attracting and repelling points in that neighbourhood.
Using it, we construct elements in $\G$ that will satisfy the mixing statement up to right multiplication by $M_\G$.
Finally decorrelation Proposition \ref{prop_decorrelation} allows to choose very contracting loxodromic elements in $\G$ whose attracting and repelling points are in a neighbourhood of $(\xi_1,\check{\xi}_1)$, of signed Jordan projection dense in an $M_\G$-invariant set.

\begin{prop}[Proposition 5.4 \cite{dang_glorieux_2020}]\label{DG-prop}
Let $G$ be a real linear, connected, semisimple Lie group of non-compact type (i.e. without compact factors) and $\Gamma$ be a Zariski dense subsemigroup of $G$.
Fix $\theta\in \fa^{++}$ of norm $1$ in the interior of the limit cone $\mathcal{C}(\G)$.

Then for every nonempty open subset $\mathcal{O}^{(2)}\subset L^{(2)}(\G)$, for all $x_0\in A$ and $\delta_0 >0$ there exists $T_0>0$ such that for all $t\geq T_0$ there exists a loxodromic element $\gamma_t \in \G$ with
\begin{equation}\label{eq-DG-prop}
 \left\{
      \begin{aligned}
        (\gamma_t^+,\gamma_t^-) &\in \mathcal{O}^{(2)} \\
       \exp\big( \lambda(\gamma_t) \big)&\in B\big(x_0 e^{t \theta},\delta_0\big)  \\
      \end{aligned}
    \right.
\end{equation}

\end{prop}

%
Recall that for every non-trivial section $s$ of the $MAN$-bundle $G \rightarrow \mathcal{F}$, we denote by $\mathcal{F}_s$ its domain and by $\mathcal{F}_s^{(2)}:= (\mathcal{F}_s \times \mathcal{F})\cap \mathcal{F}^{(2)}$ the subset of ordered transverse pairs of first coordinate in $\mathcal{F}_s$. 
Denote by $\pi_A$ the projection $AM \rightarrow A$.
\begin{proof}[Proof of Theorem \ref{theo_mixing} (SC)]
Let $\theta \in \fa^{++}$ be in the interior of the limit cone.

We want to prove the following statement:
for all non-empty open sets 
$\mathcal{U}^{(2)} \subset L^{(2)}(\G)$ and $ \mathcal{V}^{(2)} \subset L^{(2)}(\G)$,
for all $u,v \in AM_\G$ and $\delta>0$, there exists $T_1>0$ such that for every $t\geq T_1$,
for all compact Bruhat sections $c_U, c_V$ such that $\mathcal{U}^{(2)} \subset \mathcal{F}_{c_U}^{(2)}$ and $\mathcal{V}^{(2)}\subset \mathcal{F}_{c_V}^{(2)}$, in Bruhat-Hopf coordinates 
$$ \phi_\theta^t \big( \mathcal{U}^{(2)} \times B(u,\delta)  \big)_{c_U} \cap \G \big( \mathcal{V}^{(2)} \times B(v,\delta)  \big)_{c_V} \neq \emptyset, $$
meaning that there exists $h_t \in \G$ such that
$$ \big( \mathcal{U}^{(2)} \times B(ue^{t \theta},\delta)  \big)_{c_U} \cap h_t \big( \mathcal{V}^{(2)} \times B(v,\delta)  \big)_{c_V} \neq \emptyset .$$

 Consider the pair $(\xi_1,\check{\xi}_1)\in L^{(2)}(\G)$, the real positive number $r_1>0$  given by Proposition \ref{prop_decorrelation} and the associated compact Bruhat sections $c_1, \check{c}_1$.

\textbf{Step 1:} Apply topological transitivity of the action of $\G$ on $L^{(2)}(\G)$ given by Theorem \ref{transitivity of Furstenberg boundary}.
Then there exists $h_U, h_V \in \G$ such that
$$
 \left\{ \begin{aligned}
h_U \mathcal{U}^{(2)} &\ni (\xi_1, \check{\xi}_1) \\
h_V\mathcal{V}^{(2)} &\ni (\xi_1, \check{\xi}_1).
      \end{aligned}
    \right.$$
By left $\G$ invariance and right $AM_\G$ invariance of $\widetilde{\Omega}_{[e_M]}$,  there exists $u_1,v_1 \in AM_\G$ such that in Bruhat-Hopf coordinates,
 $$
 \left\{ \begin{aligned}
h_U \big( \mathcal{U}^{(2)} \times B(u,\delta) \big)_{c_U} &\ni (\xi_1, \check{\xi}_1 \; ; \; u_1)_{c_1} \\
h_V \big( \mathcal{V}^{(2)} \times B(v,\delta) \big)_{c_V} &\ni (\xi_1, \check{\xi}_1\; ; \; v_1)_{\check{c}_1}.
      \end{aligned}
    \right.$$ 
Choose $r \in (0,r_1]$ and $\delta_1>0$ small enough such that in Bruhat-Hopf coordinates
 $$
 \left\{ \begin{aligned}
h_U \big( \mathcal{U}^{(2)} \times B(u,\delta) \big)_{c_U} &\supset \big(B( \xi_1,r) \times B ( \check{\xi}_1, r) \times B (u_1,\delta_1)\big)_{c_1} \\
h_V \big( \mathcal{V}^{(2)} \times B(v,\delta) \big)_{c_V} &\supset \big( B( \xi_1,r) \times B ( \check{\xi}_1, r) \times B( v_1, \delta_1)  \big)_{\check{c}_1}.
      \end{aligned}
    \right.$$  
By Proposition \ref{prop_decorrelation}, for all $\varepsilon \in (0, r]$, there exists a finite family $(g_i)_{i\in I} \subset \G$ and a point $a_{r,\varepsilon} \in A$ satisfying the following conditions.
\begin{itemize}
\item[$\dagger$] For all $i\in I$, the element $g_i$ is $(2r,2\varepsilon)$-loxodromic with
$$ (g_i^+,g_i^-) \in B(\xi_1, \varepsilon) \times B(\check{\xi}_1,\varepsilon) .$$
\item[$\ddagger$] For all $\eta \in \mathcal{V}_{6r}\big( \partial \mathsf{b}(\check{\xi}_1) \big)^\complement$ and $(\eta_i)_{i\in I} \subset B(\eta, \varepsilon)$, the family 
$ \lbrace \beta_{c_1,\check{c}_1}(g_i, \eta_i) \rbrace_{i\in I} $
is $\tilde{\delta}_{r,\varepsilon}$-dense in $a_{r,\varepsilon}\mathscr{R}_{c_1,\check{c}_1}(\check{\xi}_1; \xi_1,\eta) M_\Gamma$ i.e. 
$$a_{r,\varepsilon} \mathscr{R}_{c_1,\check{c}_1}(\check{\xi}_1; \xi_1,\eta)  M_\Gamma \subset  \cup_{i\in I} B( \beta_{c_1,\check{c}_1}(g_i, \eta_i), \tilde{\delta}_{r,\varepsilon} ) .$$
\end{itemize}

\textbf{Step 2:} Choose $\varepsilon \in (0,r]$ such that $\tilde{\delta}_{r,\varepsilon}\leq \delta_1/2.$
Denote by $\mathcal{O}^{(2)}:=B(\xi_1, \varepsilon) \times B(\check{\xi}_1,\varepsilon)$. 
We are going to prove the topological mixing statement for $u_1, v_1 \in AM_\G$, small $\delta_1>0$, when $\mathcal{U}^{(2)}=\mathcal{V}^{(2)}=\mathcal{O}^{(2)}$. 

Let us apply Proposition \ref{DG-prop} to $\theta$ which is in the interior of the limit cone, the above open subset of $L^{(2)}(\G)$, for $x_0:= \pi_A \big(a_{r,\varepsilon}^{-1} u_1 v_1^{-1}  \big) $ and $\delta_1/2$.
We thus consider $T_0>0$ and a family of loxodromic elements $(\g_t)_{t \geq T_0}$ satisfying the system \eqref{eq-DG-prop}.
Apply $\dagger$, since $g_i^-$ is the attracting point of $g_i^{-1}$, we deduce for all $i\in I$ 
$$ \g_t^{-1}g_i^{-1} B(\check{\xi}_1, \varepsilon) \subset B(\check{\xi}_1,\varepsilon) .$$ 
Hence for all $i \in I$ and every $\check{\xi} \in \g_t^{-1}g_i^{-1} B(\check{\xi}_1, \varepsilon)$, 
\begin{align*}
g_i \g_t \big(\g_t^+, \check{\xi} \; ; \; v_1\big)_{\check{c}_1} 
&= \big(g_i \g_t^+, g_i \g_t \check{ \xi}\; ;\; \beta_{c_1,\check{c}_1}(g_i \g_t, \g_t^+) v_1 \big)_{c_1} \\
&= \big(g_i \g_t^+, g_i \g_t \check{ \xi}\; ;\; \beta_{c_1,\check{c}_1}(g_i , \g_t^+) \mathscr{L}_{\check{c}_1}(\g_t)  v_1  \big)_{c_1}\\
&\in  \mathcal{O}^{(2)} \times 
\lbrace \beta_{c_1, \check{c}_1} (g_i , \g_t^+) \mathscr{L}_{\check{c}_1}(\g_t) v_1 \rbrace.
\end{align*}

We discuss the cocycle terms using the decorrelation. 
By $\ddagger$, the set 
\begin{equation}\label{eq_set-cocycle}
\lbrace \beta_{c_1, \check{c}_1} (g_i , \g_t^+) \mathscr{L}_{\check{c}_1}(\g_t) v_1 \; \vert \; i\in I \rbrace
\end{equation}
is $\delta_1/2$-dense in 
$$ \Big( a_{r,\varepsilon} \mathscr{R}_{c_1, \check{c}_1}(\check{\xi}_1; \xi_1 , \xi_1)M_\G \Big) \mathscr{L}_{\check{c}_1}(\g_t) v_1 .$$
Since the ratio $\mathscr{R}_{c_1, \check{c}_1}(\check{\xi}_1; \xi_1 , \xi_1)$ is trivial and $M_\G$ is a normal subgroup of $M$, we deduce that the above subset of cocycles \eqref{eq_set-cocycle} is $\delta_1/2$-dense in 
$a_{r,\varepsilon} \mathscr{L}_{\check{c}_1}(\g_t) u_1 M_\G .$
Furthermore, by equation \eqref{eq-DG-prop} it is $\delta_1 $- dense in $ a_{r,\varepsilon}x_0 e^{t \theta}u_1 M_\G .$
By choice of $x_0$, remark $\pi_A( a_{r,\varepsilon}x_0 e^{ 
t \theta}v_1 )= \pi_A( u_1 e^{t\theta}).$
Hence
$$ \pi_A( u_1 e^{t \theta})M_\G \subset \bigcup_{i\in I} B( \beta_{c_1,\check{c}_1}(g_i \g_t, \g_t^+) v_1 ,  \delta_1). $$
Since $u_1 e^{t\theta} \in  \pi_A( u_1 e^{t \theta})M_\G$, we choose for all $t \geq T_0$ an element $h_t\in \lbrace g_i\g_t\rbrace_{i\in I}$ such that
$$u_1 e^{t \theta} \in B( \beta_{c_1,\check{c}_1}(h_t, \g_t^+) v_1 , \delta_1).$$
Consider $w\in B(v_1,\delta_1)$ such that $\beta_{c_1,\check{c}_1}(h_t, \g_t^+) w= u_1 e^{t \theta}$.
Then for all $\check{\xi} \in h_t^{-1} B(\check{\xi}_1,\varepsilon),$  
\begin{align*}
h_t \big(\g_t^+, \check{\xi} \; ; \; w\big)_{\check{c}_1} 
&= \big(h_t \g_t^+, h_t \check{ \xi}\; ;\; \beta_{c_1,\check{c}_1}(h_t, \g_t^+) w \big)_{c_1} \\
&= \big(h_t \g_t^+, h_t \check{ \xi}\; ;\; u_1 e^{t\theta} \big)_{c_1} \in \phi_\theta^t \big( \mathcal{O}^{(2)} \times B(u_1,\delta_1) \big)_{c_1}.
\end{align*}
Therefore, all points of such coordinates are in
$ \phi_\theta^t \big( \mathcal{O}^{(2)} \times B(u_1,\delta_1)  \big)_{c_1} \cap h_t \big( \mathcal{O}^{(2)} \times B(v_1,\delta_1) \big)_{\check{c}_1}.$
Hence for all $t \geq T_0$, there exists $h_t \in \G$ such that
\begin{equation}\label{eq-top-mix}
 \phi_\theta^t \big( \mathcal{O}^{(2)} \times B(u_1, \delta_1)  \big)_{c_1} \cap h_t \big( \mathcal{O}^{(2)} \times B(v_1, \delta_1) \big)_{\check{c}_1} \neq \emptyset. 
\end{equation}
By choice of $\varepsilon>0$, remark that 
     $$
 \left\{ \begin{aligned}
h_U \big( \mathcal{U}^{(2)} \times B(u,\delta) \big)_{c_U} &\supset \big( \mathcal{O}^{(2)} \times B (u_1,\delta_1)\big)_{c_1} \\
h_V \big( \mathcal{V}^{(2)} \times B(v,\delta) \big)_{c_V} &\supset \big( \mathcal{O}^{(2)} \times B( v_1, \delta_1)  \big)_{\check{c}_1}.
      \end{aligned}
    \right.$$ 
Note that relation \eqref{eq-top-mix} ensures 
$ \phi_\theta^t \big( h_U\big( \mathcal{U}^{(2)} \times B(u, \delta)  \big)_{c_U} \big) \cap h_t h_V\big( \mathcal{V}^{(2)} \times B(v, \delta) \big)_{\check{c}_V} \neq \emptyset. $
Since the flow commutes with left multiplication by $\G$, we deduce that for all $t \geq T_0$,
$$ \phi_\theta^t \big( \mathcal{U}^{(2)} \times B(u, \delta)  \big)_{c_U}  \cap h_U^{-1} h_t h_V\big( \mathcal{V}^{(2)} \times B(v, \delta) \big)_{\check{c}_V} \neq \emptyset. $$
\end{proof}

\section*{Appendix : density Lemmata}

\begin{lem}\label{lemme_densite2}
Let $C$ be a compact connected abelian real linear Lie group and $V$ be a finite dimensionnal real vector space.

Then for all subset $E \subset V \times C$ that span a dense subgroup in $V \times C$, for all small real number $\delta >0$, there exists a finite subset $F_{\delta} \subset E$ of at most $3 \dim V + 2 \dim C$ elements such that the subgroup generated by $F_\delta$ is $\delta$-dense in $V \times C$.
\end{lem}
It a consequence of the following Lemma.
\begin{lem}[Lemma 6.1 \cite{dang_glorieux_2020}]\label{lemme_densite1}
Let $V$ be a finite dimensionnal real vector space.

Then for all subset $E \subset V $ that span a dense subgroup in $V $, for all small real number $\delta >0$ and all basis $B \subset E$ of $V$, there exists a finite subset $F_{\delta} \subset E$ of at most $2 \dim V $ elements such that the subgroup generated by $B \cup F_\delta$ is $\delta$-dense in $V $.
\end{lem}

\begin{proof}[ Proof of Lemma \ref{lemme_densite2} ] 
By Corollary 3.7 of 
\cite{brocker1985representations}, the group $C$ is isomorphic to a torus.
Consequently, its universal cover $\tilde{C}$
is a real vector space of dimension $\dim (C)$.

Fix a small real number $\delta>0.$
Denote by $\widetilde{V}=V \times \tilde{C}$ the universal cover of $V\times C$.
Then $\widetilde{V}$ is a real vector space of dimension $\tilde{d}= d+\dim C$.

We want to apply Lemma \ref{lemme_densite1} on this vector space. 
Let us first construct out of $E$ a subset that spans a dense additive subgroup.
Denote by $p: \widetilde{V} \rightarrow V \times C$ the covering map.
Fix a basis $(b_1,...,b_d, b_{d+1},...,b_{\tilde{d}})$ of $\widetilde{V}$ such that $(p(b_1),...,p(b_d))$ is a basis of $V\times \lbrace 0\rbrace$ and the additive subgroup generated by $( b_{d+1},...,b_{\tilde{d}} ) $ is the kernel of the covering map $\ker(p)$.
With such a basis, we explicit the isomorphism between $\widetilde{V}/\langle b_{d+1},...,b_{\tilde{d}} \rangle$ and $ V \times C$.
Then the following subset 
$$\widetilde{\mathcal{D}}:= Vect(b_1,...,b_d)\times \bigg( \overset{\dim C}{\underset{j=1}{ \prod}} [0,1[b_{d+j} \bigg),$$ 
is a fundamental domain of the covering.
Consider now the subset of elements of this fundamental domain that project into elements of $E$,
$$\widetilde{E}:= p^{-1}(E) \cap \widetilde{\mathcal{D}}.$$
We deduce that $\widetilde{E}\cup  \lbrace  b_{d+1},...,b_{\tilde{d}} \rbrace$ spans a dense additive subgroup of $\widetilde{V}$.
Fix now a subset $B' \subset \widetilde{E}$ such that $\pi_V(p(B'))$ is a basis of $V$.  

Apply now density Lemma \ref{lemme_densite1} on $\widetilde{V}$, for the subset $\widetilde{E}\cup  \lbrace  b_{d+1},...,b_{\tilde{d}} \rbrace$ and choice of basis $B'\cup  \lbrace  b_{d+1},...,b_{\tilde{d}} \rbrace$.
There exists and we choose a finite subset $\widetilde{F} \subset \widetilde{E}$ of at most $2\tilde{d}$ 
elements, such that
$\widetilde{F} \cup B'\cup  \lbrace  b_{d+1},...,b_{\tilde{d}} \rbrace  $ 
spans a $\delta$-dense additive subgroup of $\widetilde{V}$.

Finally, we project $\widetilde{F} \cup B'\cup  \lbrace  b_{d+1},...,b_{\tilde{d}} \rbrace  $  to $V \times C$ using the covering map.
Then $p(\widetilde{F} \cup B') \subset E$ is a 
finite subset of at most $3d+2\dim C$ elements 
that spans a $\delta$-dense additive subgroup of  $V\times C$.
\end{proof}

\begin{lem}\label{lemme_densite3_cone}
Let $C$ be a compact connected abelian real linear Lie group and $V$ be a finite dimensionnal real vector space.
Fix $\delta >0$ a small real number.

Then for all finite subset $F \subset V \times C$
that spans a $\delta$-dense subset of $V \times C$, there exist an element $v_F \in V$ such that the \emph{semigroup} genenerated by
 $F$ is $\delta-$dense in
$$\bigg( v_F + \underset{f\in F}{\sum} \R_+ \pi_V(f) \bigg) \times C .$$
\end{lem}

\begin{proof}
We adapt a proof of Y. Benoist of Lemma 6.2 \cite{benoist2000proprietes}.

Consider the compact subset of $V$
$$\widetilde{D}:= \bigg\lbrace \underset{ f\in F}{\sum} t_f \pi_V(f) \;  \bigg\vert \; 0\leq t_f \leq 1 \bigg\rbrace.$$
Then $\widetilde{D}\times C$ is a compact subset of $V\times C$.
By hypothesis, the additive subgroup generated by 
$F$ is $\delta-$dense in $V\times C$.
Then, applying compacity, we choose a finite subset  
$X \subset \langle F\rangle$ that is $\delta$-dense in $\widetilde{D}\times C$, i.e. such that 
$$\widetilde{D}\times C \subset \bigcup_{x \in X}B(x,\delta).$$

Denote by $\langle F\rangle_+$ the subsemigroup generated by $F$.
Choose an element of the additive sub-semigroup $h \in \langle F\rangle_+$ such that $hX \subset \langle F\rangle_+$.
Such a choice is possible because $V \times C$ is \emph{abelian}. 

Then the translate $h (\widetilde{D}\times C)$ is $\delta$-covered by $hX \subset \langle F \rangle_+$, i.e.
\begin{equation} \label{eq_lem1}
h (\widetilde{D}\times C) \; \subset \; \bigcup_{x \in X}B(hx,\delta) \; \subset \; \bigcup_{x \in \langle F \rangle_+}B(x,\delta) .
\end{equation}
Remark now that
$$h (\widetilde{D}\times C)=(\pi_V(h)+\widetilde{D})\times \pi_C(h)C = (\pi_V(h)+\widetilde{D})\times C .$$
Denote now by $L$ the close convex cone generated by $\pi_V(F)$, i.e.  $L:=\sum_{f\in F} \R_+ \pi_V(f)$.
Then, by translating on the left by $\langle F \rangle_+$ in the previous equality, a translate of $L$ appears on the right hand side i.e.
$$\langle F \rangle_+\big( (\pi_V(h)+\widetilde{D})\times C  \big) = \big( (\pi_V(h)+L) \times C\big).$$
Finally, combining with \eqref{eq_lem1}, we deduce that $\langle F\rangle_+$ is 
$\delta$-dense in $\big( (\pi_V(h)+L) \times C\big)$ i.e.
$$\big( (\pi_V(h)+L) \times C\big) \subset \bigcup_{x \in \langle F \rangle_+}B(x,\delta).$$
\end{proof}
\bibliographystyle{alpha}

\end{document}